\documentclass[preprint,10pt]{elsarticle}
\usepackage[showframe = false, headheight = 2cm, margin = 2cm, bottom = 2cm]{geometry}

\usepackage{amssymb}
\usepackage{bigints}
\usepackage{mathrsfs}
\usepackage{calc}
\usepackage{tikz}
\usepackage{xstring}
\usepackage{tkz-euclide}
\usepackage{graphics}
\usepackage{multicol}
\usepackage{adjustbox}
\usepackage{subcaption}
\usepackage{multirow}
\usepackage{etoolbox}
\usepackage{datetime}
\usepackage{amsmath, accents, bm}
\usepackage{amsthm}
\usepackage{amssymb}
\usepackage{afterpage}
\usepackage{acro}
\usepackage{float}
\usepackage{relsize}
\usepackage{array}
\usepackage{textcomp}
\usepackage{comment}
\usepackage{diagbox}
\usepackage{listings}
\usepackage{tabularx}
\usepackage{ragged2e}
\usepackage{lscape}
\usepackage{enumitem}
\usepackage{eurosym}
\usepackage{footmisc}
\usepackage{todonotes}
\usepackage{mwe}
\usepackage{gensymb}
\usepackage{wasysym}
\usepackage{graphicx}
\usepackage{blindtext}
\usepackage{appendix}
\usepackage{multicol}
\usepackage{hyperref}
\usepackage{cleveref}
\usepackage{cancel}
\usepackage{bbm}
\usepackage{mathtools}
\usepackage{nomencl}
\usepackage{accents}
\usepackage{calligra}
\usepackage{miama}
\usepackage{empheq}
\usepackage{xcolor}
\usepackage{soul}
\usepackage[normalem]{ulem}
\definecolor{codegreen}{rgb}{0,0.6,0}
\usetikzlibrary{patterns.meta}

\journal{}

\theoremstyle{definition}
\newtheorem{definition}{Definition}[section]
\DeclareMathOperator*{\argmin}{arg\,min}

\theoremstyle{remark}

\begin{document}

\newcommand{\parddx}[2]{\dfrac{\partial #1}{\partial #2}}
\newcommand{\parddxdx}[2]{\dfrac{\partial^2 #1}{\partial #2^2}}
\newcommand{\E}[2]{\mathbb{E}^{#1, #2}}
\newcommand{\tE}[2]{\widetilde{\mathbb{E}}^{#1, #2}}
\newcommand{\tETranspose}[2]{\widetilde{\mathbb{E}}^{#1, #2^{T}}}
\newcommand{\ETranspose}[2]{\mathbb{E}^{#1, #2^{T}}}
\newcommand{\MassM}[1]{\mathbb{M}^{(#1)}}
\newcommand{\MassMP}[2]{\mathbb{M}^{(#1, #2)}}
\newcommand{\MassP}[3]{\mathbb{M}^{(#1, #2)}_{#3}}
\newcommand{\MassMinv}[1]{\mathbb{M}^{(#1)^{-1}}}
\newcommand{\weakForm}[2]{\int_{\Omega} #1 \wedge \star \: #2}
\newcommand{\boldCal}[1]{\boldsymbol{\mathcal{#1}}}
\newcommand{\bilnear}[2]{\left(#1, #2\right)_{L^2(\Omega)}}
\usetikzlibrary{math}
\begin{frontmatter}



\title{Optimal solutions employing an algebraic Variational Multiscale approach \\ Part I: Steady Linear Problems}

\author[1,2]{Suyash Shrestha}
\ead{s.shrestha@upm.es}
\author[2]{Marc Gerritsma}
\author[1,3]{Gonzalo Rubio}
\author[2]{Steven Hulshoff}
\author[1,3]{Esteban Ferrer}
\affiliation[1]{organization={ETSIAE-UPM-School of Aeronautics, Universidad Politécnica de Madrid},
            addressline={Plaza Cardenal Cisneros 3},
            city={Madrid},
            postcode={E-28040},
            country={Spain}}

\affiliation[2]{organization={Delft University of Technology, Faculty of Aerospace Engineering},
            addressline={Kluyverweg 1},
            city={Delft},
            postcode={2629 HS},
            country={The Netherlands}}
\affiliation[3]{organization={Center for Computational Simulation, Universidad Politécnica de Madrid, Campus de Montegancedo},
            addressline={Boadilla del Monte},
            city={Madrid},
            postcode={E-28660},
            country={Spain}}




\begin{abstract}
This work extends our previous study from S. Shrestha et al. (2024) by introducing a new abstract framework for Variational Multiscale (VMS) methods at the discrete level. We introduce the concept of what we define as the optimal projector and present a discretisation approach that yields a numerical solution closely approximating the optimal projection of the infinite-dimensional continuous solution. In this approach, the infinite-dimensional unresolved scales are approximated in a finite-dimensional subspace using the numerically computed Fine-Scale Greens' function of the underlying symmetric problem. The proposed approach involves solving the VMS problem on two separate meshes: a coarse mesh for the full PDE and a fine mesh for the symmetric part of the continuous differential operator. We consider the 1D and 2D steady advection-diffusion problems in both direct and mixed formulations as the test cases in this paper. We first present an error analysis of the proposed approach and show that the projected solution is achieved as the approximate Greens' function converges to the exact one. Subsequently, we demonstrate the working of this method where we show that it can exponentially converge to the chosen optimal projection. We note that the implementation of the present work employs the Mimetic Spectral Element Method (MSEM), although, it may be applied to other Finite/Spectral Element or Isogeometric frameworks. Furthermore, we propose that VMS should not be viewed as a stabilisation technique; instead, the base scheme should be inherently stable, with VMS enhancing the solution quality by supplementing the base scheme.

\end{abstract}




\end{frontmatter}

\tableofcontents


\section{Introduction}
The present work serves as a continuation of \cite{Shrestha2024ConstructionScales}, where we formulate a Variational Multiscale (VMS) \cite{Hughes1998TheMechanics, ShakibFarzin1989FiniteEquations} approach using the Fine-Scale Green's function to obtain optimal solutions for advection-diffusion problems. In \cite{Shrestha2024ConstructionScales} we introduced a generalised method for explicitly computing the Fine-Scale Greens' function for any projector. Possessing this Fine-Scale Greens' function allows one to tackle the fine-scale problem within the VMS framework. This also presents challenges, however, such as the need to analytically derive and numerically integrate the classic Greens' function of the associated PDE. The numerical integration issue can be mitigated with advanced quadrature rules \cite{Mori2006NumericalTransformation}, but the analytical derivation of Greens' functions is only feasible for simple problems, see \cite[\S 8]{RichardHaberman2019HabermanBoundary}. In \cite{Shrestha2024ConstructionScales}, we highlighted that, in fact, only the Greens' function for the \emph{symmetric} part of the differential operator is needed, which we leverage in this paper.

Our proposed solution to the analytical derivation challenge involves solving the VMS problem using two separate meshes: a \emph{coarse mesh} for the full PDE and a \emph{fine mesh} for computing the classic Greens' function for the \emph{symmetric} part of the differential operator. These meshes differ in refinement levels, adjusted by selecting different polynomial degrees. This paper specifically uses the Mimetic Spectral Element Method (MSEM) \cite{Jain2021ConstructionMeshes}, but the concepts can be applied to other Finite/Spectral Element or Isogeometric frameworks.

The VMS approach considers the computation of a solution equal to a chosen projection of the exact solution onto the finite-dimensional space where the discretisation is performed. A key aspect of the VMS approach is thus the choice of \emph{projector}. While any projector can theoretically be chosen,
it is not trivial as to which one can be regarded as the ``optimal" projector. Optimality depends on the desired solution characteristics. For instance, one might prefer a nodally exact solution and consider a nodal projector to be the optimal projector, while another might aim to minimise the $L^2$ error and consider the $L^2$ projector as the optimal. We offer our view of optimality by defining the projector based on the PDE in question. We demonstrate the existence of a unique projector for any symmetric PDE that allows direct computation of the projected solution using only the input data, i.e., the PDE source term and boundary conditions. This projector, constructed using the energy norm of the symmetric PDE, is what we define as the \emph{optimal} projector. 

It is important to clarify the use of the terminology ``optimal projector" in this manuscript. Here, ``optimal projector" refers specifically to the particular choice of projector we have made based on its suitability for the governing PDE. This terminology does not imply that our approach is universally superior or the ultimate solution among all possible projectors. Instead, it reflects a choice tailored to the PDE's structure and this study's objectives. While other projectors may be optimal in their own contexts or norms, the goal of this study is not to explore all such possibilities. Rather, we aim to demonstrate the efficacy of our framework in recovering optimal solutions based on the defined projector.

We begin in \Cref{sec:optmal_sym} by reviewing the concept of obtaining the optimal solution for symmetric PDEs. Initially, we discuss the optimal projector in a general context for any generic symmetric differential operator, followed by specific examples involving the Poisson equation in both direct and mixed formulations. Then in \Cref{sec:optmal_skew_sym}, we derive the VMS formulation for the steady advection-diffusion problem using the optimal projector defined for the symmetric problem. In \Cref{sec:greens_func}, we explain the numerical computation of the (Fine-Scale) Greens' functions for the symmetric operator, summarise the VMS approach in an abstract framework, and present an error analysis. Additionally, we present our argument and motivation on the particular choice of the projector based on the symmetric operator. Next, we present numerical tests in \Cref{sec:numerical_ex} where we demonstrate the linear theory by applying the VMS approach to the 1D steady advection-diffusion equation in both direct and mixed formulations and the 2D steady advection-diffusion equation in the mixed formulation. Finally, in \Cref{sec:summary}, we conclude with a summary of our work and outline future steps.
\section{Optimal projections: Symmetric operators}
\label{sec:optmal_sym}
We consider the problem statement defined as follows. Let $\Omega \subset \mathbb{R}^d$ be an open bounded domain with a sufficiently smooth Lipschitz continuous boundary $\partial \Omega$. Consider $V$ to be a Hilbert space in $\Omega$ endowed with a norm $||\cdot||_V$ and a scalar product $(\cdot,\cdot)_V$. We take $V^*$ to be the dual space of $V$ with $\prescript{}{V^*}{\langle\cdot,\cdot\rangle}_{V}$ being the duality pairing between the two. Let $\mathcal{L}: V \to V^*$ be a linear isomorphism between $V$ and $V^*$, the problem can then be stated as: given $f \in V^*$, find $\phi \in V$ such that
\begin{gather}
    \mathcal{L} \phi = f, \quad \text{in } \Omega \label{eq:pde} \\
    \phi = g, \quad \text{on } \Gamma_d \\
    \nabla \phi \cdot \hat{n} = h, \quad \text{on } \Gamma_n \\
    \text{with } \Gamma_d \cup \Gamma_n = \partial \Omega,
\end{gather}
where $\hat{n}$ is the outward unit normal vector on $\partial \Omega$ and $g$ and $h$ are prescribed functions. For a symmetric operator $\mathcal{L}$, we may derive a weak form of \eqref{eq:pde} by considering the energy minimisation principle in a suitably chosen norm, typically the energy norm, see \cite{Arnold2010FINITESTABILITY, Boffi2013MixedApplications, Marsden2008TheMethods}. The energy functional associated with the PDE is expressed in a generic form as 
\begin{gather}
    J(\phi; f, \{g, h\}) := \frac{1}{2} \left\lVert \phi \right\rVert_{X}^2 + \mathcal{D}(\phi, f) + \mathcal{B}(\phi, \{g, h\}), \label{eq:energy_norm}
\end{gather}
where $\left\lVert \cdot \right\rVert_{X}$ is a norm defined in $\Omega$ and $\mathcal{D}(\cdot, \cdot)$ and $\mathcal{B}(\cdot, \cdot)$ are inner products in $\Omega$ and $\partial \Omega$ respectively. 
Taking variations with respect to $\phi$ yields a variational (weak) formulation of \eqref{eq:pde} which reads: given $f \in V^*$, $g$, and $h$, find $\phi \in V$ such that
\begin{gather}
    a(v, \phi) = (v, f)_{L^2(\Omega)} + b(v, \{g, h\}),  \quad \forall v \in V, \label{eq:weak_bilinear}
\end{gather}
where $a(\cdot, \cdot)$ is a bilinear form in $\Omega$ and $b(\cdot, \cdot)$ is a bilinear form in $\partial \Omega$, both resulting from $\mathcal{L}$. Since we want to solve \eqref{eq:pde} numerically we work in $\bar{V}$ which we define to be a closed finite-dimensional subspace of $V$. The traditional \emph{Galerkin} approach to arrive at a numerically implementable formulation is to restrict both the \emph{test space} and \emph{trial space} to the finite-dimensional subspace $\bar{V}$ which gives
\begin{gather}
    a(v^h, {\phi}^h) = (v^h, f)_{L^2(\Omega)} + b(v^h, \{g, h\}),  \quad \forall v^h \in \bar{V}.
    \label{eq:weak_glk}
\end{gather}
The restriction of the trial space implies that the solution $\phi$ is approximated by ${\phi}^h \in \bar{V}$. If the problem is well-posed, one can show that this weak formulation satisfies the original problem in \eqref{eq:pde} if the solution is sufficiently smooth, by invoking integration by parts on the bilinear form $a(\cdot, \cdot)$. We can further show that the solution of \eqref{eq:weak_glk} is in fact a \emph{projection}
of the infinite-dimensional (strong) solution $\phi \in V$ onto the finite-dimensional subspace $\bar{V}$. The projector $\mathcal{P}$ in question is based on the energy norm $J(\phi; 0, 0)$, which is a linear mapping $\mathcal{P} :  V \rightarrow \Bar{V}$ with $\mathcal{P}^2 = \mathcal{P}$, $\text{Range}(\mathcal{P}) = \Bar{V}$, and is stated as a minimisation problem as
\begin{gather}
    \mathcal{P} \phi := \argmin_{\bar{\phi} \in \bar{V}} \left\{\frac{1}{2} \left\lVert \bar{\phi} - \phi \right\rVert_{X}^2\right\}, \label{eq:proj_Y}
\end{gather}
where $\bar{\phi} \in \bar{V}$ is the projected quantity being sought and $\phi \in V$ is the exact (strong) solution to \eqref{eq:pde}. 
Taking variations of the above expression yields the variational form of the problem expressed as
\begin{gather}
    a(v^h, \bar{\phi}) = a(v^h, \phi), \quad \forall v^h \in \bar{V},
    \label{eq:proj_bilinear_fd}
\end{gather}
with $a(\cdot, \cdot)$ being the aforementioned bilinear form. The infinite-dimensional term on the right-hand side of \eqref{eq:proj_bilinear_fd} can be recognised to be the continuous level weak form in \eqref{eq:weak_bilinear} with \emph{just} the test space restricted to $\bar{V}$
\begin{gather}
    a(v^h, {\phi}) = (v^h, f)_{L^2(\Omega)} + b(v^h, \{g, h\}), \quad \forall v^h \in \bar{V}.
\end{gather}
Substituting this into \eqref{eq:proj_bilinear_fd} gives
\begin{gather}
    a(v^h, \bar{\phi}) = (v^h, f)_{L^2(\Omega)} + b(v^h, \{g, h\}), \quad \forall v^h \in \bar{V}.
    \label{eq:proj_weak_fd}
\end{gather}
We end up with the same finite-dimensional variational (weak) formulation as in \eqref{eq:weak_glk}, however, we arrive at it by employing the notion of a projector and the continuous level weak form. Hence, the solution to \eqref{eq:weak_glk} is the projection of the strong solution i.e. ${\phi}^h = \bar{\phi} = \mathcal{P}\phi$. Note that the specific projector based on the energy norm of the PDE allows the bilinear form on the right-hand side of \eqref{eq:proj_bilinear_fd} to be fully characterised using only the problem input data i.e. the PDE source term and boundary conditions. Noting this, we formulate a definition of what we call the \emph{optimal} projector.
\begin{definition}
    \textit{We define the \emph{optimal} projector as the one based on the energy norm of the PDE which yields a bilinear form $a(\cdot, \cdot)$ for which the infinite-dimensional term $a(v^h, \phi)$ can be completely characterised using only the problem data, i.e. the source term of the PDE and the boundary conditions.}
\end{definition}
In the context of the Variational Multiscale (VMS) method, we refer to $\bar{\phi}$ as the so-called resolved (coarse) scales which are representable in our finite-dimensional subspace. Furthermore, we may define the so-called unresolved (fine) scales as $\phi' := \bar{\phi} - \phi$, which are all the components truncated away by the projector that cannot be represented in our finite-dimensional setting. We note that the unresolved scales live in an infinite-dimensional subspace $V'$ which is the null space of the projector $V' = \text{Ker}(\mathcal{P})$. Given this definition of the unresolved scales, the projection statement for the optimal projector reads
\begin{gather}
    \mathcal{P} \phi := \argmin_{\bar{\phi} \in \bar{V}} \left\{\frac{1}{2} \left\lVert \bar{\phi} - \phi \right\rVert_{X}^2\right\} = \argmin_{{\phi}^{\prime} \in {V}^{\prime}} \left\{\frac{1}{2} \left\lVert \phi^{\prime} \right\rVert_{X}^2\right\}, \label{eq:proj}
\end{gather}
which we may interpret as the minimisation of the energy of the fine scales. Hence the optimal projector minimises the energy of the fine scales which cannot be represented in our finite-dimensional subspace. Moreover, these fine scales are orthogonal to the elements of the resolved space in the energy norm i.e.
\begin{gather}
    a(\bar{\phi}, \phi') = 0, \quad \forall \bar{\phi} \in \bar{V}, \: \forall \phi' \in V'.
\end{gather}
We shall demonstrate this concept by considering a Poisson problem.

\subsection{Poisson problem with a direct formulation}
Consider a Poisson problem stated as follows
\begin{gather}
    -\nabla \cdot \nabla \phi = f, \quad  \text{in } \Omega \label{eq:poisson} \\
    \phi = g, \quad \text{on} \: \Gamma_d \label{eq:poisson_Dbc}\\
    \nabla \phi \cdot \hat{n} = h, \quad \text{on} \: \Gamma_n \label{eq:poissonNbc}\\
    \text{with } \Gamma_d \cup \Gamma_n = \partial \Omega.
\end{gather}
We start by defining the Sobolev space $H^1_0(\Omega) := \{\varphi \in H^1 \: | \: \varphi = 0 \: \text{on } \Gamma_d\}$ and $H^{-1}(\Omega)$ as its dual space. We seek $\phi \in V = H^1_0(\Omega)$ given $f \in V^* = H^{-1}(\Omega)$, and for our finite-dimensional setting we take $\bar{V}$ to be a closed (polynomial) subspace of $H^1_0(\Omega)$. The energy functional associated with this PDE reads as follows
\begin{gather}
    J(\phi; f, h) := \frac{1}{2} \left\lVert \phi \right\rVert_{H^1_0}^2 + \mathcal{D}(\phi, f) + \mathcal{B}(\phi, h) = \frac{1}{2} \int_{\Omega} |\nabla \phi |^2 \: \mathrm{d}\Omega - \int_{\Omega} \phi f \: \mathrm{d}\Omega - \int_{\partial \Omega} \phi h \: \mathrm{d}\Gamma_n. \label{eq:poisson_energy}
\end{gather}
In physical terms, if we consider \eqref{eq:poisson} with $f = 0$, as describing a potential flow problem with $\phi$ as the potential, then $\nabla \phi$ is the (irrotational) velocity and the first term in \eqref{eq:poisson_energy} defines the kinetic energy. The remaining terms define the work done along the boundary. 
Taking variations with respect to $\phi$ gives the following variational form
\begin{gather}
    \left(\nabla v, \nabla \phi \right)_{L^2(\Omega)} = \left(v, f \right)_{L^2(\Omega)} + \int_{\partial \Omega} v h \: \mathrm{d} \Gamma_n, \quad \forall v \in V, \label{eq:weak_poisson}
\end{gather}
which is a weak restatement of \eqref{eq:poisson} with $\phi, v \in V$ (not $\bar{V}$).


If we consider the projector based on the energy norm $J(\phi; 0, 0)$, we get the $H_0^1$ projector for which the minimisation problem reads
\begin{gather}
    \mathcal{P}_{H_0^1} \phi = \argmin_{\bar{\phi} \in \bar{V}} \left\{\frac{1}{2} \left\lVert \nabla \bar{\phi} - \nabla \phi \right\rVert_{L^2}^2\right\}. \label{eq:H01_min}
\end{gather}
Taking variations of \eqref{eq:H01_min} leads to
\begin{gather}
    \left(\nabla v^h, \nabla \bar{\phi}\right)_{L^2(\Omega)} = \left(\nabla v^h, \nabla 
 \phi \right)_{L^2(\Omega)}, \quad \forall v^h \in \bar{V}. \label{eq:H01_opt}
\end{gather}
The right-hand side of \eqref{eq:H01_opt} can be identified as being the weak form in \eqref{eq:weak_poisson} but with the test space restricted to $\bar{V}$. Hence we have
\begin{gather}
    \left(\nabla v^h, \nabla \bar{\phi}\right)_{L^2(\Omega)} = \left(v^h, f \right)_{L^2(\Omega)} + \int_{\partial \Omega} v^h h \: \mathrm{d} \Gamma_n, \quad \forall v^h \in \bar{V}. \label{eq:poisson_fem}
\end{gather}
We can strongly impose the Dirichlet boundary conditions on the above system and we arrive at a formulation which we can implement to solve \eqref{eq:poisson} numerically provided that $\phi$ is sufficiently smooth. We recognise that \eqref{eq:poisson_fem} is just the weak form of \eqref{eq:weak_poisson} with both the test and trial space restricted to $\bar{V}$ i.e. the Galerkin formulation. However, we arrive to it by considering the minimisation problem associated with the $H_0^1$ projector. 
Moreover, in the context of the VMS method, we can note that the fine scales are orthogonal with the resolved (projector) solution in the energy norm i.e.
\begin{gather}
    (\nabla \bar{\phi}, \nabla \phi')_{L^2(\Omega)} = 0, \quad \forall \bar{\phi} \in \bar{V}, \: \forall \phi' \in V'.
\end{gather}
\subsection{Poisson problem in a mixed formulation}
\label{subsec:poisson_mix}
Having considered the direct formulation for the Poisson problem, we may alternatively consider the mixed formulation of the PDE. For example, the Poisson problem in \eqref{eq:poisson} may be expressed in a mixed formulation as
\begin{gather}
    \underline{u} - \nabla \phi = 0, \quad \text{in } \Omega \\
    \nabla \cdot \underline{u} = -f, \quad \text{in } \Omega \\
    \phi = g, \quad \text{on} \: \Gamma_d \\
    \underline{u} \cdot \hat{n} = h, \quad \text{on} \: \Gamma_n \\
    \text{with } \Gamma_d \cup \Gamma_n = \partial \Omega.
\end{gather}
We define the Sobolev space $H(\mathrm{div}, \Omega) := \{ \underline{\varphi} \in [L^2(\Omega)]^d \: | \: \nabla \cdot \underline{\varphi} \in L^2(\Omega) \}$ and we seek $\underline{u} \in H(\mathrm{div}, \Omega)$, $\phi \in L^2(\Omega)$ given $f \in L^2(\Omega)$, and we take two sets of finite-dimensional spaces, $\bar{V}$ as the closed subspace of $H(\mathrm{div}, \Omega)$ and $\bar{W}$ as the closed subspace of $L^2(\Omega)$. We still have an energy functional for the mixed formulation expressed using the new variable $\underline{u}$, however, we no longer have a minimisation problem.
Instead, the pair $(\underline{u}, \phi) \in H(\mathrm{div}, \Omega) \times L^2(\Omega)$ are variationally characterised as the saddle point of the energy functional which reads \cite{Arnold2010FINITESTABILITY}
\begin{gather}
    I(\underline{u}, \phi; f, g) := \frac{1}{2} \int_{\Omega} |\underline{u}|^2 \: \mathrm{d}\Omega + \int_{\Omega} \phi \nabla \cdot \underline{u} \: \mathrm{d}\Omega - \int_{\Omega} \phi  f \: \mathrm{d}\Omega - \int_{\partial \Omega} g \underline{u} \cdot \hat{n} \: \mathrm{d}\Gamma_d.
\end{gather}
Taking variations of the above functional we get
\begin{gather}
    \left(\underline{v}, \underline{u} \right)_{L^2(\Omega)} + \left(\nabla \cdot \underline{v}, \phi \right)_{L^2(\Omega)} = \int_{\partial \Omega} g \underline{v} \cdot \hat{n} \: \mathrm{d} \Gamma_d,  \quad \forall \underline{v} \in H(\mathrm{div}, \Omega) \label{eq:mix_poi_1}\\
    \left(\eta, \nabla \cdot \underline{u}\right)_{L^2(\Omega)} = -\left(\eta, f \right)_{L^2(\Omega)}, \quad \forall \eta \in L^2(\Omega). \label{eq:mix_poi_2}
\end{gather}
As done before, we can express this problem in terms of a projector $\mathcal{P}$. More specifically, the projection in question is a constrained minimisation problem which reads
\begin{gather}
    \mathcal{P} \underline{u} := \argmin_{\underline{\bar{u}} \in \bar{V}} \left\{\frac{1}{2} \left\lVert \underline{\bar{u}} - \underline{u} \right\rVert_{L^2}^2 \right\}, \quad \quad s.t \: \nabla \cdot \underline{\bar{u}} = \mathcal{P}_{L^2} f,
    \label{eq:mix_proj}
\end{gather}
where $\mathcal{P}_{L^2} f$ is the $L^2$ projection of the source term of the PDE. This constrained minimisation problem can be uniquely characterised as the saddle point of the energy functional
\begin{gather}
    {I}(\underline{\bar{u}} - \underline{u}, \bar{\phi} - \phi; 0, 0) = \frac{1}{2} \int_{\Omega} |\underline{\bar{u}} - \underline{u}|^2 \: \mathrm{d}\Omega + \int_{\Omega} (\bar{\phi} - \phi) \nabla \cdot (\underline{\bar{u}} - \underline{u}) \: \mathrm{d}\Omega.
\end{gather}
If we take variations with respect to $(\underline{\bar{u}} - \underline{u})$ we get
\begin{gather}
    (\underline{v}^h, \underline{\bar{u}})_{L^2(\Omega)} + (\nabla \cdot \underline{v}^h, \bar{\phi})_{L^2(\Omega)} = (\underline{v}^h, \underline{u})_{L^2(\Omega)} + (\nabla \cdot \underline{v}^h, \phi)_{L^2(\Omega)}, \quad \forall \underline{v}^h \in \bar{V}.
\end{gather}
Once again, we can completely characterise the right-hand side of the above equation by simply restricting the test space to $\bar{V}$ in \eqref{eq:mix_poi_1}
\begin{gather}
    (\underline{v}^h, \underline{\bar{u}})_{L^2(\Omega)} + (\nabla \cdot \underline{v}^h, \bar{\phi}) = \int_{\partial \Omega} g \underline{v}^h \cdot \hat{n} \: \mathrm{d} \Gamma_d,  \quad \forall \underline{v}^h \in \bar{V}.
\end{gather}
Similarly, if we take variations with respect to $\bar{\phi} - \phi$ we get
\begin{gather}
    (\eta^h, \nabla \cdot \underline{\bar{u}})_{L^2(\Omega)} = (\eta^h, \nabla \cdot \underline{u})_{L^2(\Omega)}, \quad \forall \eta^h \in \bar{W} \\
    (\eta^h, \nabla \cdot \underline{\bar{u}})_{L^2(\Omega)} = -(\eta^h, f)_{L^2(\Omega)}, \quad \forall \eta^h \in \bar{W}.
\end{gather}
An important remark about the choice of function spaces is that they satisfy the De Rham sequence, i.e. applying the divergence to an element of $H(\mathrm{div}, \Omega)$ maps the element to $L^2(\Omega)$ which implies that the divergence operator is a surjective mapping between $H(\mathrm{div}, \Omega)$ and $L^2(\Omega)$. This should also hold at the finite-dimensional setting to guarantee the well-posedness of the mixed formulation. See \cite[\S 5.3]{Jain2021ConstructionMeshes} regarding the discrete inf-sup condition.
\begin{equation*}
\begin{array}{ccc}
    H(\mathrm{div}, \Omega) & \xrightarrow[\hspace{2cm}]{\displaystyle \nabla \cdot} & L^2(\Omega)\\ \\
    \bar{V} & \xrightarrow[\displaystyle \nabla_h \cdot]{\hspace{2cm}} & \bar{W}
\end{array}
\end{equation*}

We note that the difference between the direct and mixed formulations relates to whether the divergence or gradient operator is applied weakly or strongly. In the direct formulation for the Poisson equation, we apply the strong gradient and the weak divergence which yields natural Neumann boundary conditions. On the other hand, we apply the weak gradient and the strong divergence in the mixed formulations which gives natural Dirichlet boundaries. 
\section{Optimal projections: skew-symmetric operators}
\label{sec:optmal_skew_sym}
We now turn our attention to a more general class of problems involving skew-symmetric operators with the problem stated as
\begin{gather}
    \mathcal{C}\phi + \mathcal{L} \phi = f, \quad \text{in } \Omega \label{eq:pde_non_sym} \\
    \phi = g, \quad \text{on } \Gamma_d \\
    \nabla \phi \cdot \hat{n} = h, \quad \text{on } \Gamma_n \\
    \text{with } \Gamma_d \cup \Gamma_n = \partial \Omega,
\end{gather} 
with $\mathcal{C}$ being a \emph{linear} skew-symmetric operator (an advection term for example) and $\mathcal{L}$ being a symmetric operator. We still work in function spaces dictated by the symmetric operator, i.e. a Hilbert space $V$ defined in $\Omega$ where we seek $\phi \in V$ at the continuous level and $\bar{\phi} \in \bar{V}$ for the discrete setting with the optimal projector $\mathcal{P} : V \xrightarrow{} \bar{V}$. However, unlike the case of the symmetric problem, it is not possible to derive a weak form of \eqref{eq:pde_non_sym} from an energy minimisation principle unless we move to a Lagrangian formulation where the skew-symmetric part vanishes. We choose to stick to an Eulerian formulation, hence, we derive a weak form of \eqref{eq:pde_non_sym} by testing it with a test function $v \in V$ which yields
\begin{gather}
    a(v, \phi) + c(v, \phi) = (v, f)_{L^2(\Omega)} + b(v, \{g, h\}), \quad \forall v \in V.
    \label{eq:weak_non_symm}
\end{gather}
The resulting weak form is similar to \eqref{eq:weak_bilinear} with the addition of the bilinear form $c(v, \phi)$ emerging from the skew-symmetric term. The inability to derive the weak form from an energy minimisation principle also implies that we cannot directly derive a finite-dimensional variational form of \eqref{eq:weak_non_symm} using the notion of an optimal projector as we had demonstrated for the symmetric case. In short, this means that if we were to directly derive a finite-dimensional system from \eqref{eq:weak_non_symm} by restricting the test and trial space to $\bar{V}$ (classic Galerkin formulation), the discrete solution we would obtain by solving the system would \emph{not} be the optimal projection of the exact solution! 

In order to ensure we get the optimal projection from our discretisation using \eqref{eq:weak_non_symm}, we employ Variational Multiscale (VMS) analysis wherein the unresolved components of the solution are explicitly accounted for in the variational equation. We know that invoking the optimal projector associated with the symmetric operator yields the following
\begin{gather}
    a(v^h, \bar{\phi}) = a(v^h, \phi).
\end{gather}
The projector allows us to separate the resolved ($\bar{\phi}$) and unresolved (fine) scales ($\phi^{\prime}$) where the continuous $\phi$ can be uniquely expressed as 
\begin{gather}
    \phi := \bar{\phi} + \phi^{\prime}, \quad \bar{\phi} \in \bar{V}, \: \phi^{\prime} \in V^{\prime}, \label{eq:phi_split} 
\end{gather}
and further implies orthogonality between the resolved and unresolved (fine) scales
\begin{gather}
    a(\bar{\phi}, \phi^{\prime}) = 0. \label{eq:orth_subscale}
\end{gather}
If we substitute the split form from \eqref{eq:phi_split} into \eqref{eq:weak_non_symm} we get
\begin{gather}
    a(v, \bar{\phi}) + a(v, \phi^{\prime}) + c(v, \bar{\phi}) + c(v, \phi^{\prime}) = (v, f)_{L^2(\Omega)} + b(v, \{g, h\}), \quad \forall v \in V.
\end{gather}
If we restrict the test space to $\bar{V}$, we get
\begin{gather}
    a(v^h, \bar{\phi}) + \cancel{a(v^h, \phi^{\prime})} + c(v^h, \bar{\phi}) + c(v^h, \phi^{\prime}) = (v^h, f)_{L^2(\Omega)} + b(v^h, \{g, h\}), \quad \forall v^h \in \bar{V} \\
    a(v^h, \bar{\phi}) + c(v^h, \bar{\phi}) + c(v^h, \phi^{\prime}) = (v^h, f)_{L^2(\Omega)} + b(v^h, \{g, h\}), \quad \forall v^h \in \bar{V}
\end{gather}
where the term $a(v^h, \phi^{\prime}) = 0$ due to the condition in \eqref{eq:orth_subscale}. We thus find that we need information on the fine-scales $\phi^{\prime}$ to find the optimal solution $\bar{\phi}$. Naturally, $\phi^{\prime}$ only appears in the bilinear form associated with the skew-symmetric operator. In the absence of the skew-symmetric term we have the formulation to find the optimal solution of the symmetric problem that was previously discussed. 

The unresolved scales $\phi^{\prime}$ are naturally unknowns and thus must be solved simultaneously with $\bar{\phi}$. We employ the Fine-Scale Greens' function established in \cite{Shrestha2024ConstructionScales} to compute $\phi^{\prime}$ and its gradients. An important point to note is that we only use Greens' function of the \emph{symmetric problem} and not the full skew-symmetric problem. The latter case could be argued as `cheating' as the Greens' function of the full PDE would directly yield the exact solution to the original problem. Additionally, finding the Greens' function for the full PDE is as difficult as solving the original problem, making the VMS approach obsolete. Hence, we proceed with the Greens' function of only the symmetric operator which overcomes the aforementioned flaws/shortcomings. We express the fine-scale problem in a generic setting by substituting the split form from \eqref{eq:phi_split} into the original PDE in \eqref{eq:pde_non_sym} and constraining it as follows
\begin{gather}
    \mathcal{L} \phi' + \mathcal{P}^T \bar{\lambda} = f - \mathcal{L} \bar{\phi} - (\mathcal{C} \bar{\phi} + \mathcal{C} \phi^{\prime}) \label{eq:fine_scale_con0}\\
    \mathcal{P} \phi' = 0, \label{eq:fine_scale_con1}
\end{gather}
which is the same format used in \cite{Shrestha2024ConstructionScales, Hughes2007VariationalMethods} where $\bar{\lambda} \in \bar{V}$ is the Lagrange multiplier used to constraint $\phi^{\prime}$ to $V'$. Given that $\mathcal{L}$ in invertible ($\mathcal{L}^{-1} = \mathcal{G}$ i.e. the Greens' function), we may write
\begin{gather}
    \phi' = \mathcal{G}(r - \mathcal{P}^T \bar{\lambda}), \quad \text{with } r := f - \mathcal{L} \bar{\phi} - (\mathcal{C} \bar{\phi} + \mathcal{C} \phi^{\prime}).
\end{gather}
When substituted into \eqref{eq:fine_scale_con1} and rearranged for $\bar{\lambda}$, we get
\begin{gather}
    \bar{\lambda} = \left(\mathcal{P} \mathcal{G} \mathcal{P}^T\right)^{-1} \mathcal{P} \mathcal{G} r.
\end{gather}
Hence, we have 
\begin{gather}
    \phi' = (\mathcal{G} - \mathcal{GP}^T (\mathcal{PGP}^T)^{-1}\mathcal{PG}) r = \mathcal{G}^{\prime} r,
\end{gather}
where $\mathcal{G}^{\prime}$ is the Fine-Scale Greens' function, see \cite{Shrestha2024ConstructionScales, Hughes2007VariationalMethods} for more details. Expanding the full expression for $r$ gives
\begin{gather}
    \phi' = \mathcal{G}^{\prime}f - \cancel{\mathcal{G}^{\prime}\mathcal{L} \bar{\phi}} - \mathcal{G}^{\prime}\mathcal{C} \bar{\phi} - \mathcal{G}^{\prime} \mathcal{C} \phi^{\prime}) = \mathcal{G}^{\prime}(\mathscr{R}\bar{\phi} - \mathcal{C} \phi^{\prime}),
    \label{eq:phi_prime}
\end{gather}
where we define $\mathscr{R}\bar{\phi} := f - \mathcal{C} \bar{\phi}$ as the coarse scale residual with $\mathcal{L} \bar{\phi}$ being eliminated as it lives exactly in the finite-dimensional dual space $\bar{V}^*$ which is orthogonal to $\mathcal{G}^{\prime}$ given the property $\mathcal{G}^{\prime} \mathcal{P}^T = 0$. With a closed form expression for $\phi^{\prime}$ we now have a coupled system we can solve to obtain $\bar{\phi}$, namely
\begin{gather}
    a(v^h, \bar{\phi}) + c(v^h, \bar{\phi}) + c(v^h, \phi^{\prime}) = (v^h, f)_{L^2(\Omega)} + b(v^h, \{g, h\}), \quad \forall v^h \in \bar{V} \\
    \phi' = \mathcal{G}^{\prime}(\mathscr{R}\bar{\phi} - \mathcal{C} \phi^{\prime}). \label{eq:phi_prime_2}
\end{gather}
However, due to the fact that $\phi^{\prime}$ appears on both sides of \eqref{eq:phi_prime_2} we cannot directly solve the coupled system for an arbitrary operator $\mathcal{C}$. 
Fortunately, we can simplify \eqref{eq:phi_prime_2} using the linearity of $\mathcal{C}$. We first realise that the bilinear form $c(\cdot, \cdot)$ is given by
\begin{gather}
    c(\alpha, \beta) := (\alpha, \mathcal{C}\beta)_{L^2(\Omega)},
\end{gather}
hence we are interested in finding $\mathcal{C}\phi^{\prime}$ which we can compute by simply applying $\mathcal{C}$ to \eqref{eq:phi_prime_2}. 
\begin{gather}
    \mathcal{C}\phi' = \mathcal{C}\mathcal{G}^{\prime}(\mathscr{R}\bar{\phi} - \mathcal{C} \phi^{\prime}).
\end{gather}
We can then rearrange the equation to get
\begin{gather}
    \mathcal{C} \phi' + \mathcal{C} \mathcal{G}^{\prime} \mathcal{C} \phi^{\prime} = \mathcal{C} \mathcal{G}^{\prime}\mathscr{R}\bar{\phi} \\
    (\mathbb{I} + \mathcal{C} \mathcal{G}^{\prime}) \mathcal{C} \phi^{\prime} = \mathcal{C} \mathcal{G}^{\prime}\mathscr{R}\bar{\phi} \\
    \mathcal{C} \phi^{\prime} = (\mathbb{I} + \mathcal{C} \mathcal{G}^{\prime})^{-1} \mathcal{C} \mathcal{G}^{\prime}\mathscr{R}\bar{\phi} \\
    \mathcal{C} \phi^{\prime} = \sigma^{\mathcal{C}}_{SG} \mathscr{R}\bar{\phi},
\end{gather}
where we have defined a new operator $\sigma^{\mathcal{C}}_{SG}$ which we refer to as the Suyash-Greens' operator. A common theme in the literature on the algebraic VMS approach is to express the unresolved scales as $\phi^{\prime} = \tau \mathscr{R}\bar{\phi}$ where $\tau$ is a parameter determined through scaling arguments and is usually based on some a-priori expected/desired behaviour of the fine scales, see \cite{ShakibFarzin1989FiniteEquations, Bazilevs2007VariationalFlows, Brezzi1997BG, Holmen2004SensitivityFlow, tenEikelder2018CorrectContext, Koobus2004AShedding} for examples. However, the Suyash-Greens' operator is derived directly from the variational form of the PDE and the Fine-Scale Greens' function of the symmetric operator, hence it does not rely on any ad-hoc derivation.

\subsection{Steady advection-diffusion with a direct formulation}
Having derived the general formulation for a generic linear skew-symmetric operator, we now focus on a particular example problem in the form of a steady advection-diffusion equation. Consider the advection-diffusion problem and its underlying symmetric Poisson problem described as follows

\begin{minipage}{0.49\linewidth}
    \begin{gather}
        \underline{c} \nabla \phi - \nabla \cdot (\kappa \nabla \phi) = f, \quad \text{in } \Omega \label{eq:adv_diff} \\
        \phi = g, \quad \text{on } \Gamma_d \\
        \nabla \phi \cdot \hat{n} = h, \quad \text{on } \Gamma_n \\
        \text{with } \Gamma_d \cup \Gamma_n = \partial \Omega,
    \end{gather}
\end{minipage}
\begin{minipage}{0.49\linewidth}
    \begin{gather}
        -\nabla \cdot (\kappa \nabla \phi) = f, \quad \text{in } \Omega \label{eq:poisson_ad} \\
        \phi = g, \quad \text{on } \Gamma_d \\
        \nabla \phi \cdot \hat{n} = h, \quad \text{on } \Gamma_n \\
        \text{with } \Gamma_d \cup \Gamma_n = \partial \Omega,
    \end{gather}
\end{minipage}
\vspace{0.5em}

\noindent with $\underline{c}$ being a prescribed uniformly continuous vector field and $\kappa$ being a symmetric positive definite matrix. We motivate the choice of the function spaces based on the symmetric problem. Hence, we seek $\phi \in V = H^1_0(\Omega)$ given $f \in V^* = H^{-1}(\Omega)$, and for our finite-dimensional setting, we take $\bar{V}$ to be a closed (polynomial) subspace of $H^1_0(\Omega)$. The \emph{optimal} projector for the symmetric problem reads
\begin{gather}
    \mathcal{P} \phi = \argmin_{\bar{\phi} \in \bar{V}} \left\{\frac{1}{2} \left\lVert \sqrt{\kappa} (\nabla \bar{\phi} - \nabla \phi) \right\rVert_{L^2}^2\right\}, \quad \text{with } \left\lVert \sqrt{\kappa} (\nabla \varphi) \right\rVert_{L^2}^2 := (\nabla \varphi, \kappa \nabla \varphi)_{L^2(\Omega)}.
\end{gather}
Based on this projector, the exact solution can be uniquely characterised as
\begin{gather}
    \phi := \bar{\phi} + \phi^{\prime}, \quad \text{with } \bar{\phi} \in \bar{V}, \: \phi^{\prime} \in V'. \label{eq:split_phi}
\end{gather}
Additionally, the resolved and unresolved scales are orthogonal in the energy norm
\begin{gather}
    (\nabla \bar{\varphi}, \kappa \nabla \varphi')_{L^2(\Omega)} = 0, \quad \forall \bar{\varphi} \in \bar{V}, \: \forall \varphi' \in V'. \label{eq:ortho_dir}
\end{gather}
If we test \eqref{eq:adv_diff} with $v^h \in \bar{V}$, we get
\begin{gather}
    (v^h, \underline{c} \nabla \phi)_{L^2(\Omega)} + (\nabla v^h, \kappa \nabla \phi)_{L^2(\Omega)} = (v^h, f)_{L^2(\Omega)} + \int_{\partial \Omega} v^h h \: \mathrm{d}\Gamma_n, \quad \forall v^h \in \bar{V}.
\end{gather}
Substituting the split form from \eqref{eq:split_phi} in the above equation gives
\begin{gather}
    (v^h, \underline{c} \nabla \bar{\phi})_{L^2(\Omega)} + (v^h, \underline{c} \nabla \phi^{\prime})_{L^2(\Omega)} + (\nabla v^h, \kappa \nabla \bar{\phi})_{L^2(\Omega)} = (v^h, f)_{L^2(\Omega)} + \int_{\partial \Omega} v^h h \: \mathrm{d}\Gamma_n, \quad \forall v^h \in \bar{V}. \label{eq:bar_phi_ad}
\end{gather}
The term $(\nabla v^h, \kappa \nabla \phi^{\prime})_{L^2(\Omega)}$ is exactly zero due to the orthogonal property in \eqref{eq:ortho_dir}. However, $(v^h, c \nabla \phi^{\prime})_{L^2(\Omega)}$ is certainly non-zero. If we ignore the term $(v^h, c \nabla \phi^{\prime})_{L^2(\Omega)}$, we get the well-known Galerkin weak formulation which will not yield the projected solution i.e. $\bar{\phi} \neq \mathcal{P}\phi$.
Including $(v^h, c \nabla \phi^{\prime})_{L^2(\Omega)}$ in the equation allows us to compute the solution $\bar{\phi}$ which is the optimal projection of the exact solution.

We thus need to formulate an additional equation for the unresolved scales to allow the $(v^h, c \nabla \phi^{\prime})_{L^2(\Omega)}$ term to be determined. To obtain this expression for $\phi^{\prime}$, we consider the fine-scale problem which we express as a Poisson problem. The fine-scale problem for the advection-diffusion problem expressed as a symmetric Poisson problem reads
\begin{gather}
    -\nabla \cdot (\kappa \nabla \phi^{\prime}) + \mathcal{P}^T \bar{\lambda} = f - \underline{c} \nabla \bar{\phi} - \underline{c} \nabla \phi^{\prime}, \quad \text{in } \Omega \label{eq:poisson_fine} \\
    \mathcal{P} \phi^{\prime} = 0,
\end{gather}
which we can simplify using the symmetric Fine-Scale Greens' function as
\begin{gather}
    \phi^{\prime} = \mathcal{G}^{\prime} \mathscr{R}\bar{\phi} - \mathcal{G}^{\prime}\underline{c} \nabla \phi^{\prime},
\end{gather}
where we have used the coarse scale residual $\mathscr{R}\bar{\phi} = f - \underline{c} \nabla \bar{\phi}$ with the Laplacian of $\bar{\phi}$ being eliminated given that $\mathcal{G}^{\prime} \mathcal{P}^T = 0$. We can then obtain the following expression for the gradient of $\phi^{\prime}$
\begin{gather}
    \nabla \phi^{\prime} = \nabla \mathcal{G}^{\prime}  \mathscr{R}\bar{\phi} - \nabla \mathcal{G}^{\prime}\underline{c} \nabla \phi^{\prime}. \label{eq:dphi_prime}
\end{gather}
We once again note that we only need the (Fine-Scale) Greens' function for the symmetric (linear) Poisson problem and not the full advection-diffusion problem. Since the advection term is linear, we can simplify \eqref{eq:dphi_prime} to obtain the Suyash-Greens' operator and we get
\begin{gather}
    \nabla  \phi^{\prime} = \sigma^{\nabla}_{SG} \mathscr{R}\bar{\phi}, \quad \quad \text{with } \sigma^{\nabla}_{SG} := (\mathbb{I} + \nabla \mathcal{G}^{\prime}\underline{c})^{-1} \nabla \mathcal{G}^{\prime}.
\end{gather}
This expression for $\nabla  \phi^{\prime}$ can be directly plugged into \eqref{eq:bar_phi_ad} which gives
\begin{gather}
    (v^h, \underline{c} \nabla \bar{\phi})_{L^2(\Omega)} + (v^h, \underline{c} \nabla \sigma^{\nabla}_{SG} \mathscr{R}\bar{\phi})_{L^2(\Omega)} + (\nabla v^h, \kappa \nabla \bar{\phi})_{L^2(\Omega)} = (v^h, f)_{L^2(\Omega)} + \int_{\partial \Omega} v^h h \: \mathrm{d}\Gamma_n \\
    \begin{split}
        (v^h, \underline{c} \nabla \bar{\phi})_{L^2(\Omega)} - (v^h, \underline{c} \sigma^{\nabla}_{SG} \underline{c} \nabla \bar{\phi})_{L^2(\Omega)} + &(\nabla v^h, \kappa \nabla \bar{\phi})_{L^2(\Omega)} = \\ &(v^h, f)_{L^2(\Omega)} - (v^h, \underline{c} \sigma^{\nabla}_{SG} f)_{L^2(\Omega)} + \int_{\partial \Omega} v^h h \: \mathrm{d}\Gamma_n.
    \end{split}
\end{gather}
Solving the above system directly yields the optimal solution $\bar{\phi}$ which is the $H^1_0$ projection of the exact solution onto the finite-dimensional polynomial subspace provided that we have the exact Greens' function for the Poisson equation.


\subsection{Steady advection-diffusion with a mixed formulation}
We may, alternatively, also consider the advection-diffusion problem in a mixed formulation for which the full problem statement and its underlying symmetric problem read 

\begin{minipage}{0.49\linewidth}
    \begin{gather}
        \kappa^{-1} \underline{q} - \nabla \phi = 0, \quad \text{in } \Omega \label{eq:adv_diff_mix_1}\\
        \underline{c} \cdot \kappa^{-1} \underline{q} - \nabla \cdot \underline{q} = f, \quad \text{in } \Omega \label{eq:adv_diff_mix_2} \\
        \phi = g, \quad \text{on } \Gamma_d \\
        \underline{q} \cdot \hat{n} = h, \quad \text{on } \Gamma_n \\
        \text{with } \Gamma_d \cup \Gamma_n = \partial \Omega.
    \end{gather}
\end{minipage}
\begin{minipage}{0.49\linewidth}
    \begin{gather}
        \kappa^{-1} \underline{q} - \nabla \phi = 0, \quad \text{in } \Omega \\
        -\nabla \cdot \underline{q} = f, \quad \text{in } \Omega \\
        \phi = g, \quad \text{on } \Gamma_d \\
        \underline{q} \cdot \hat{n} = h, \quad \text{on } \Gamma_n \\
        \text{with } \Gamma_d \cup \Gamma_n = \partial \Omega.
    \end{gather}
\end{minipage}
\vspace{0.5em}

\noindent Similar to the mixed formulation for the Poisson equation, here we seek $(\underline{q}, \phi) \in H(\mathrm{div}, \Omega) \times L^2(\Omega)$ given $f \in L^2(\Omega)$, and for our finite-dimensional setting we take $\bar{V}$ to be a closed (polynomial) subspace of $H(\mathrm{div}, \Omega)$ and $\bar{W}$ to be a closed (polynomial) subspace of $L^2(\Omega)$ provided they satisfy the inf-sup conditions.

We have already established that the pair $(\underline{q}, \phi)$ can be uniquely characterised as the saddle point of the energy functional. Subsequently, we may define a projection as follows
\begin{gather}
    {I}(\underline{\bar{q}} - \underline{q}, \bar{\phi} - \phi; 0, 0) := \frac{1}{2} \int_{\Omega} \left(\underline{\bar{q}} - \underline{q} \right)^T \kappa^{-1} \left(\underline{\bar{q}} - \underline{q} \right)  \: \mathrm{d}\Omega + \int_{\Omega} (\bar{\phi} - \phi) \nabla \cdot (\underline{\bar{q}} - \underline{q}) \: \mathrm{d}\Omega. \label{eq:proj_mixed}
\end{gather}
Using the projector, we can separate the resolved and unresolved scales and uniquely express the continuous quantities $\underline{q}$ and $\phi$ as
\begin{gather}
    \underline{q} := \bar{\underline{q}} + \underline{q}^{\prime}, \quad \text{with } \bar{\underline{q}} \in \bar{V}, \: \underline{q}^{\prime} \in V' \\
    \phi := \bar{\phi} + \phi^{\prime}, \quad \text{with } \bar{\phi} \in \bar{W}, \: \phi^{\prime} \in W'.
\end{gather}
Additionally, we have the following orthogonality conditions
\begin{gather}
    (\bar{\varphi}, \kappa^{-1} \varphi')_{L^2(\Omega)} + (\nabla \cdot \bar{\varphi}, \varsigma^{\prime})_{L^2(\Omega)} = 0, \quad \forall \bar{\varphi} \in \bar{V}, \: \forall \varphi' \in V', \: \forall \varsigma^{\prime} \in W^{\prime} \label{eq:mix_ortho_1}\\
    (\bar{\zeta}, \nabla \cdot \varphi^{\prime})_{L^2(\Omega)} = 0, \quad \forall \bar{\zeta} \in \bar{W}, \: \forall \varphi^{\prime} \in V^{\prime}. \label{eq:mix_ortho_2}
\end{gather}

We test \eqref{eq:adv_diff_mix_1} with $v^h \in \bar{V}$ and test \eqref{eq:adv_diff_mix_2} with $\eta^h \in \bar{W}$ to get
\begin{gather}
    \bilnear{\underline{v}^h}{\kappa^{-1} \underline{\bar{q}}} + \bilnear{\nabla \cdot \underline{v}^h}{\bar{\phi}} = \int_{\partial \Omega} \underline{v}^h g \cdot \hat{n} \: \mathrm{d}\Gamma_d, \quad \forall v^h \in \bar{V}  \\
    \bilnear{\eta^h}{\underline{c} \cdot \kappa^{-1} \underline{\bar{q}}} + \bilnear{\eta^h}{\underline{c} \cdot \kappa^{-1} \underline{{q}}^{\prime}} - \bilnear{\eta^h}{\nabla \cdot \underline{\bar{q}}} = \bilnear{\eta^h}{f}, \quad \forall \eta^h \in \bar{W},
\end{gather}
where we use the orthogonality conditions from \eqref{eq:mix_ortho_1} and \eqref{eq:mix_ortho_2} to eliminate some of the fine scale terms. Once again, the only fine-scale term which does not drop out is the bilinear form emerging from the convection term. The fine-scale problem for the mixed formulation is identical to that for the direct formulation highlighted in \eqref{eq:poisson_fine}. However, the projector $\mathcal{P}$ in question is now different, hence, it must be adjusted accordingly using \eqref{eq:proj_mixed}. In short, we can arrive at the following variational form
\begin{gather}
    \bilnear{\underline{v}^h}{\kappa^{-1} \underline{\bar{q}}} + \bilnear{\nabla \cdot \underline{v}^h}{\bar{\phi}} = \int_{\partial \Omega} \underline{v}^h g \cdot \hat{n} \: \mathrm{d}\Gamma_d, \quad \forall v^h \in \bar{V} \label{eq:weak_form_mix1} \\
    \begin{split}
        \bilnear{\eta^h}{\underline{c} \cdot \kappa^{-1} \underline{\bar{q}}} - \bilnear{\eta^h}{\underline{c} \cdot \kappa^{-1} \sigma_{SG}^{\nabla} \underline{c} \cdot \bar{\underline{q}}} - &\bilnear{\eta^h}{\nabla \cdot \underline{\bar{q}}} = \\ \bilnear{\eta^h}{f}& - \bilnear{\eta^h}{\underline{c} \cdot \kappa^{-1} \sigma_{SG}^{\nabla} f}, \quad \forall \eta^h \in \bar{W}. \label{eq:weak_form_mix2}
    \end{split}
\end{gather}
Just as the case for the direct formulation, solving the above system yields the optimal projection of the exact solution provided that we have the exact Greens' function for the Poisson problem to construct $\sigma_{SG}^{\nabla}$. 

\section{Greens' functions and a priori error estimates}
\label{sec:greens_func}

In this section, we briefly describe the construction of the Fine-Scale Greens' function established in \cite{Shrestha2024ConstructionScales} and present how we numerically obtain the classic Greens' function for the symmetric operator. Moreover, we also present an abstract perspective of the VMS approach using the numerical Greens' function. Lastly, we present an analysis of the error estimates for the proposed methodology and conclude the section with a rationale behind our definition of the optimal projector.

\subsection{Fine-Scale Greens' function}
In \cite{Shrestha2024ConstructionScales} we have already established how one can explicitly compute the Fine-Scale Greens' function when given the classic Greens' function and the dual basis functions associated with the projector. We shall briefly recap the main concept of \cite{Shrestha2024ConstructionScales} in this section. We first recall that the expression for the Fine-Scale Greens' function is given by
\begin{equation}
    \mathcal{G}' = \mathcal{G} - \mathcal{G}{\bm{\mu}}^T \left[{\bm{\mu}} \mathcal{G} {\bm{\mu}}^T \right]^{-1} {\bm{\mu}} \mathcal{G},
    \label{eq:g_prime_mu}
\end{equation}
where $\mathcal{G}$ is the classic Greens' function and the $\bm{\mu}$ is a set of linear functionals associated with the projector. In \cite{Shrestha2024ConstructionScales} we have shown that the $\bm{\mu}$ can be uniquely computed for any projector using dual basis functions. In short, if we have a projector $\mathcal{P} : V \rightarrow \bar{V}$ defined with a norm $\lVert \cdot \rVert_X$, then the dual basis functions ($\mu$) must satisfy $(\mu_i, \psi_{j})_{X} = \delta_{i,j}$ where $\psi_{j}$ are the basis functions spanning the $\bar{V}$ subspace. Thus, if we have the classic Greens' function, we are able to compute the Fine-Scale Greens' function for any given projector. By now, we also know that we do not need the Greens' function for the full differential operator, but only that for the symmetric operator.

\subsection{Approximate classic Greens' functions}
\label{subsec:greens_approx}
Despite the fact that the Greens' function for the symmetric operator is linear in many cases, there can still be numerous difficulties in acquiring and computing its integral. We now address the major question of how to numerically obtain the classic Greens' function $\mathcal{G}$ for the symmetric operator. We have previously established that \eqref{eq:proj_weak_fd} can be used to solve for the optimal solution to any given symmetric operator. Numerically, we would do this by selecting a set of basis functions that span $\bar{V}$ and constructing an algebraic system by evaluating the bilinear form in \eqref{eq:proj_weak_fd}. Suppose that we have chosen a set of basis functions ${\psi}_i(\bm{x})$ that span $\bar{V}$ and we organise them in a row vector and store the expansion coefficients of the basis functions i.e. the degrees of freedom $\mathcal{N}(\bar{\phi})$, in a column vector as follows
\begin{gather}
    \bm{\psi}(\bm{x}) := \left[{\psi}_0(\bm{x}), {\psi}_1(\bm{x}), \hdots {\psi}_{\mathrm{N}}(\bm{x}) \right], \quad \quad \quad \quad \quad \mathcal{N}(\bar{\phi}) := \left[\begin{array}{c}
        \mathcal{N}(\bar{\phi}_0)  \\
        \mathcal{N}(\bar{\phi}_1) \\
        \vdots  \\
        \mathcal{N}(\bar{\phi}_{\mathrm{N}})
    \end{array} \right].
\end{gather}
The final algebraic system constructed through \eqref{eq:proj_weak_fd} may then be expressed as follows
\begin{gather}
    \mathbb{A} \mathcal{N}(\bar{\phi}) = \hat{f} \quad \Longrightarrow \quad \bar{\phi} = \bm{\psi}(\bm{x}) \mathcal{N}(\bar{\phi}),
\end{gather}
where $\mathbb{A}$ is the discrete matrix operator and $\hat{f}$ is a vector with the bilinear forms of the right-hand side of \eqref{eq:proj_weak_fd} as its entries. The matrix operator $\mathbb{A}$ acts as a discrete surrogate for the continuous linear operator $\mathcal{L}$. Subsequently, we have that $\mathbb{A}^{-1}$ is a surrogate for $\mathcal{L}^{-1} = \mathcal{G}$ i.e. the classic Greens' function of the symmetric problem. We thus have
\begin{align}
    \mathcal{G} f := \int_{\Omega} g(\bm{x}, \bm{s}) f(\bm{s}) \: \mathrm{d} \Omega_s &\Longleftrightarrow \bm{\psi}(\bm{x}) \mathbb{A}^{-1} \int_{\Omega} (\bm{\psi}(\bm{s}))^T f(s) \: \mathrm{d}\Omega_s \\
    &\Longleftrightarrow \int_{\Omega} \bm{\psi}(\bm{x}) \mathbb{A}^{-1} 
 (\bm{\psi}(\bm{s}))^T f(\bm{s}) \: \mathrm{d}\Omega_s.
\end{align}
Consequently, an approximation of the Greens' function is given by
\begin{gather}
    g(\bm{x}, \bm{s}) \approx g_h(\bm{x}, \bm{s}) = \bm{\psi}(\bm{x}) \mathbb{A}^{-1} 
 (\bm{\psi}(\bm{s}))^T.
    \label{eq:greens_approx}
\end{gather}
The above expression provides a consistent approximation for the classic Greens' function provided that a suitable basis $\bm{\psi}(\bm{x})$ is chosen. This concept has also been employed in literature, see for example \cite{Hughes2007VariationalMethods, Melchers2024NeuralEquations, Praveen2023PrincipledData}. In \ref{app:A} we demonstrate this methodology by presenting plots of the Greens' function for the Poisson problem.




\subsection{The abstract VMS framework}
In an abstract sense, the VMS approach may be explained as follows: Given a Hilbert space $V$ with $\bar{V}$ as a closed subspace of $V$, we seek the optimal projection of the continuous solution onto $\bar{V}$. The infinite-dimensional unresolved scales live in $V^{\prime}$, the orthogonal complement of $\bar{V}$. If we have the analytical (Fine-Scale) Greens' functions, we have the situation depicted in \Cref{fig:VMS_abstract_1} where we would be able to recover all the fine scales provided we are able to accurately integrate the Fine-Scale Greens' function. 

If we instead opt to use the numerically computed Greens' function using \eqref{eq:greens_approx}, we must work with two meshes, a coarse mesh and a fine mesh. 
We formulate the VMS approach on the coarse mesh and we use the fine mesh to compute the Greens' function using \eqref{eq:greens_approx}. In this particular implementation, we modify the refinement level between the two meshes by changing the polynomial degree while maintaining the same number of elements. The coarse mesh has a polynomial degree of $p$ and the fine mesh has a polynomial degree $p + k$ with $k > 0$. This situation is depicted in \Cref{fig:VMS_abstract_2} where we have the orthogonal complement of $\bar{V}$, now denoted as $V^{\prime}_k$ which is a closed finite-dimensional subspace wherein the fine scales are approximated. In this scenario, we first generate a basis spanning $\bar{V} \oplus V^{\prime}_k$ and use them in \eqref{eq:greens_approx} to compute the approximate classic Greens' function. Once the classic Greens' function is computed we subtract away the components that live in $\bar{V}$ which yields the Fine-Scale Greens' function defined solely in $V^{\prime}_k$. The proposed framework shares similarities with \cite{Giraldo2023ANorms}, but differs in two key aspects: we employ a Greens' function approach using only the symmetric part of the operator and we additionally incorporate the mixed formulation of the problem.
\begin{figure}[ht]
\begin{minipage}{0.49\linewidth}
    \centering
    \begin{tikzpicture}
    \draw[dotted, thick] (0,0) circle (3cm);
    
    \filldraw[dotted, fill=gray!10, pattern={Lines[
                  distance=2mm,
                  angle=-45,
                  line width=0.1mm
                 ]}] (0,0) circle (3cm);
    
    \filldraw[fill=gray!70] (0,0) circle (1cm);

    \node at (2.5, -2.5) {$V$};

    \node at (0, 0) {$\bar{V}$};

    \node at (1.5, 1.5) {${V}^{\prime}$};
\end{tikzpicture}
    \caption{Abstract framework at the continuous level}
    \label{fig:VMS_abstract_1}
\end{minipage}
\begin{minipage}{0.49\linewidth}
    \centering
    \begin{tikzpicture}
    \draw[dotted, thick] (0,0) circle (3cm);
    
    \filldraw[fill=gray!1, pattern={Lines[
                  distance=2mm,
                  angle=-45,
                  line width=0.05mm
                 ]}] (0,0) circle (2cm);
    
    \filldraw[fill=gray!70] (0,0) circle (1cm);

    \node at (2.5, -2.5) {$V$};

    \node at (0, 0) {$\bar{V}$};

    \node at (1, 1) {${V}^{\prime}_k$};
\end{tikzpicture}
    \caption{Abstract framework at the discrete level}
    \label{fig:VMS_abstract_2}
\end{minipage}
\end{figure}

\subsection{Error analysis}
In this section, we present an error analysis of the proposed approach. We start the error analysis by considering the well-established analysis for symmetric problems \cite{Brenner2008TheMethods, Ern2021FiniteI}. We revisit the generic weak form of a symmetric (elliptic) problem described in \Cref{sec:optmal_sym}
\begin{gather}
    a(v, \phi) = \bilnear{v}{f} + \bilnear{v}{\{g, h\}}, \quad \forall v \in V. \label{eq:bilinear_sym_gen}
\end{gather}
We assume that $a(\cdot, \cdot)$ is bounded and coercive in $V$ and $\bilnear{\cdot}{f} + \bilnear{\cdot}{\{g, h\}}$ is continuous in $V$, hence the above problem has a unique solution. Using the Galerkin approach, we wish to assess the error of solving the above equation in a finite-dimensional subspace of $V$. To do so, we take the difference between \eqref{eq:bilinear_sym_gen} with both test and trial spaces restricted to the finite-dimensional subspace and \eqref{eq:bilinear_sym_gen} with only the test space being restricted
\begin{align}
    a(v^h, \phi^h) &= \bilnear{v^h}{f} + \bilnear{v^h}{\{g, h\}}, \quad \forall v^h \in \bar{V} \nonumber\\
    -\{ a(v^h, \phi) &= \bilnear{v^h}{f} + \bilnear{v^h}{\{g, h\}}\}, \quad \forall v^h \in \bar{V} \nonumber\\
    \cline{1-2}
    a(v^h, & \phi^h - \phi) = 0, \quad \forall v^h \in \bar{V}. \label{eq:glk_ortho_sym}
\end{align}
The expression in \eqref{eq:glk_ortho_sym} tells us that the error or unresolved scales are orthogonal to $\bar{V}$ in the norm induced by the bilinear form $a(\cdot, \cdot)$. To get an error estimate, we must consider \emph{a} norm to measure the error. Since the bilinear form $a(\cdot, \cdot)$ emerges from a symmetric problem, it induces a well-defined norm in $V$, and we choose this norm to carry out the analysis
\begin{gather}
    \lVert \cdot \rVert_V^2 := a(\cdot, \cdot).
\end{gather}
We thus have 
\begin{align}
    \lVert \phi^h - \phi \rVert_V^2 &= a(\phi^h - \phi, \phi^h - \phi) \\
    &= \cancel{a(\phi^h - \phi, \underbrace{\phi^h - v^h}_{\in \bar{V}})} + a(\phi^h - \phi, v^h - \phi) \\
    &\leq \lVert \phi^h - \phi \rVert_V \lVert v^h - \phi \rVert_V \quad (\text{Cauchy-Schwarz})\\
    \lVert \phi^h - \phi \rVert_V &\leq \lVert v^h - \phi \rVert_V, \quad \forall v^h \in \bar{V} \quad \Leftrightarrow \lVert \phi^h - \phi \rVert_V = \inf_{v^h \in \bar{V}} \lVert v^h - \phi \rVert_V,
\end{align}
which says that the Galerkin solution is in fact \emph{the} optimal projection of the exact solution onto $\bar{V}$ in the energy norm. \\ \\
Next, we consider the generic form of a skew-symmetric problem described in \Cref{sec:optmal_skew_sym}
\begin{gather}
    a(v, \phi) + c(v, \phi) = (v, f)_{L^2(\Omega)} + b(v, \{g, h\}), \quad \forall v \in V.
    \label{eq:bilinear_non_sym_gen}
\end{gather}
Using the same assumptions on the boundedness and coercivity of $a(\cdot, \cdot) + c(\cdot, \cdot)$ and continuity of $\bilnear{\cdot}{f} + \bilnear{\cdot}{\{g, h\}}$ as for the symmetric case, we perform the same error analysis. 
\begin{align}
    a(v^h, \phi^h) + c(v^h, \phi^h)&= \bilnear{v^h}{f} + \bilnear{v^h}{\{g, h\}}, \quad \forall v^h \in \bar{V} \nonumber\\
    -\{ a(v^h, \phi) + c(v^h, \phi)&= \bilnear{v^h}{f} + \bilnear{v^h}{\{g, h\}}\}, \quad \forall v^h \in \bar{V} \nonumber\\
    \cline{1-2}
    a(v^h, \phi^h - \phi) + &c(v^h, \phi^h - \phi) = 0, \quad \forall v^h \in \bar{V}. \label{eq:glk_ortho_non_sym}
\end{align}
Again, we choose a norm to estimate the error. We argue that the bilinear form $a(\cdot, \cdot)$ (associated with the elliptic problem) induces a well-defined norm on $V$ whereas the bilinear form $c(\cdot, \cdot)$ does not induce a norm (it violates positive definite property). We thus choose to use the bilinear form $a(\cdot, \cdot)$ as done for the symmetric problem
\begin{align}
    \lVert \phi^h - \phi \rVert_V^2 &= a(\phi^h - \phi, \phi^h - \phi) \\
    &= a(\phi^h - \phi, \underbrace{\phi^h - v^h}_{\in \bar{V}}) + a(\phi^h - \phi, v^h - \phi) \\
    &= -c(\phi^h - v^h, \phi^h - \phi) + a(\phi^h - \phi, v^h - \phi) \\
    &\leq -c(\phi^h - v^h, \phi^h - \phi) + \lVert \phi^h - \phi \rVert_{V} \lVert v^h - \phi \rVert_{V} \\
    \lVert \phi^h - \phi \rVert_V &\leq -\frac{c(\phi^h - v^h, \phi^h - \phi)}{\lVert \phi^h - \phi \rVert_V} + \lVert v^h - \phi \rVert_{V}, \quad \forall v^h \in \bar{V}.
\end{align}
Clearly, the Galerkin approximation for the skew-symmetric problem is sub-optimal in the $\lVert \cdot \rVert_V$ norm given that $c(\phi^h - v^h, \phi^h - \phi) \neq 0$. Moreover, the Galerkin solution is \emph{not} the optimal projection of the exact solution. \\ \\
We now switch to the VMS formulation where we apply a direct sum decomposition of $V$ based on a projector $\mathcal{P}$. We select our optimal projector which is defined using the norm induced by the symmetric bilinear form $a(\cdot, \cdot)$. The direct sum implies the following
\begin{gather}
    V := \bar{V} \oplus V^{\prime} \quad \xrightarrow{ } v = v^h + v^{\prime}, \phi = \bar{\phi} + \phi^{\prime}, \quad \text{where } a(\bar{\varphi}, \varphi^{\prime}) = a(\varphi^{\prime}, \bar{\varphi}) = 0, \quad \forall \bar{\varphi}\in \bar{V}, \forall {\varphi}^{\prime} \in V^{\prime}.
\end{gather}
Filling in this decomposition in \eqref{eq:bilinear_non_sym_gen} and eliminating the terms conforming to the orthogonality condition gives
\begin{gather}
    a(v^h, \bar{\phi}) + a(v^{\prime}, \phi^{\prime}) + c(v^h + v^{\prime}, \bar{\phi}) + c(v^h + v^{\prime}, \phi^{\prime}) = \bilnear{v^h + v^{\prime}}{f} + \bilnear{v^h + v^{\prime}}{\{g, h\}}.
\end{gather}
Since the equation must hold for $\forall v^h \in \bar{V}, \forall v^{\prime} \in V^{\prime}$, we may separate the two as follows
\begin{align}
    a(v^h, \bar{\phi}) + c(v^h, \bar{\phi}) + c(v^h, \phi^{\prime}) &= \bilnear{v^h}{f} + \bilnear{v^h}{\{g, h\}}, \quad \forall v^h \in \bar{V} \label{eq:resolved_eq}\\
    a(v^{\prime}, \phi^{\prime}) + c(v^{\prime}, \bar{\phi}) + c(v^{\prime}, \phi^{\prime}) &= \bilnear{v^{\prime}}{f} + \bilnear{v^{\prime}}{\{g, h\}}, \quad \forall v^{\prime} \in V^{\prime}. \label{eq:unresolved_eq}
\end{align}
We see that we get two sets of coupled equations where the coupling is introduced through the non-symmetric bilinear form. In \cite{Shrestha2024ConstructionScales} we proposed to solve \eqref{eq:resolved_eq} and \eqref{eq:unresolved_eq} in an iterative fashion using the (Fine-Scale) Greens' function associated with the symmetric term. In the present work, we choose to compute an approximate classic Greens' function by essentially computing the inverse of $a(v, \phi)$ and constructing the Suyash-Greens' operator to monolithically solve \eqref{eq:unresolved_eq}. This latter case may be summarised as follows
\begin{align}
    a(v^h, \bar{\phi}) + c(v^h, \bar{\phi}) + c(v^h, \phi^{\prime}_k) &= \bilnear{v^h}{f} + \bilnear{v^h}{\{g, h\}}, \quad \forall v^h \in \bar{V}, \label{eq:resolved_eq_aprox}
\end{align}
where $\phi^{\prime}_k \in V^{\prime}_k$ is an approximation of the unresolved scales computed using the Suyash-Greens' operator with the approximate classic Greens' function for the symmetric operator. We now perform the error analysis in a similar fashion as before where we take the difference between \eqref{eq:resolved_eq_aprox} and \eqref{eq:bilinear_non_sym_gen} with the restricted test space
\begin{align}
    &a(v^h, \bar{\phi}) + c(v^h, \bar{\phi}) + c(v^h, \phi^{\prime}_k) = \bilnear{v^h}{f} + \bilnear{v^h}{\{g, h\}}, \quad \forall v^h \in \bar{V} \nonumber\\
    -\{&a(v^h, \phi) + c(v^h, \phi) = \bilnear{v^h}{f} + \bilnear{v^h}{\{g, h\}}, \quad \forall v^h \in \bar{V} \nonumber\\
    \cline{1-2}
    &a(v^h, \bar{\phi} - \phi ) + c(v^h, \phi^{\prime}_k - \phi^{\prime}) = 0, \quad \forall v^h \in \bar{V},
\end{align}
where we use the fact that $\phi^{\prime} = \phi - \bar{\phi}$ on the second bilinear form. Using the previously defined norm, we get the following expression for the error
\begin{align}
    \lVert \bar{\phi} - \phi \rVert_V^2 &= a(\bar{\phi} - \phi, \bar{\phi} - \phi) \\
    &= a(\bar{\phi} - \phi, \bar{\phi} - v^h) + a(\bar{\phi} - \phi, v^h - \phi) \\
    &= -c(\bar{\phi} - v^h, \phi^{\prime}_k - \phi^{\prime}) + a(\bar{\phi} - \phi, v^h - \phi) \\
    &\leq -c(\bar{\phi} - v^h, \phi^{\prime}_k - \phi^{\prime}) + \lVert \bar{\phi} - \phi \rVert_V \lVert v^h - \phi \rVert_V \\
    \lVert \bar{\phi} - \phi \rVert_V &\leq -\frac{c(\bar{\phi} - v^h, \phi^{\prime}_k  - \phi^{\prime})}{\lVert \bar{\phi} - \phi \rVert_V}  + \lVert v^h - \phi \rVert_V, \quad \forall v^h \in \bar{V}.
\end{align}
If we have the analytical classic Greens' function of the symmetric operator, we are able to recover the exact unresolved scales $\phi^{\prime}_k - \phi^{\prime} = 0$ as presented in \cite{Shrestha2024ConstructionScales}. This would be the case depicted in \Cref{fig:VMS_abstract_1}, and we would thus have
\begin{gather}
    \lVert \bar{\phi} - \phi \rVert_V \leq \cancel{-\frac{c(\bar{\phi} - v^h, \phi^{\prime}_k  - \phi^{\prime})}{\lVert \bar{\phi} - \phi \rVert_V}} + \lVert v^h - \phi \rVert_V, \quad \forall v^h \in \bar{V} \\
    \lVert \bar{\phi} - \phi \rVert_V = \inf_{v^h \in \bar{V}} \lVert v^h - \phi \rVert_V,
\end{gather}
i.e. the numerical solution $\bar{\phi}$ is the optimal projection of the exact solution. This approach inherently requires access to the exact classical Greens' function. However, the newly proposed method (\Cref{fig:VMS_abstract_2}), which employs an approximate Greens' function, circumvents this requirement. The trade-off is that $\phi^{\prime}_k - \phi^{\prime}$ is non-zero but only converges to zero as we improve the approximation of the classic Greens' function by increasing the mesh increment parameter $k$. Hence as the approximate Greens' function converges to the exact one, $\bar{\phi}$ converges to the optimal projection of the exact solution. The analysis relies solely on the well-posedness of the weak form, specifically its boundedness and coercivity (see \cite{Canuto2006EnhancedProblems} for details), ensuring that the proposed method is applicable as long as the problem remains well-posed. We can perform a similar analysis for the mixed formulation case and arrive at the same conclusion.
\subsection{Rationale behind the optimal projector}
We may now take the time to reflect upon the choice of norm used to define the optimal projector. Here, the norm induced by the symmetric bilinear form was chosen. This can certainly be seen as an unconventional or incomplete choice as one can easily see that the projector does not 'see' the non-symmetric components of the solution. However, we may ask ourselves, if the projector even needs to 'see' the symmetric or non-symmetric terms of the solution at all. The true physics of the system we are solving is embedded in the infinite-dimensional exact solution and the projector is simply a means of obtaining a finite-dimensional representation of it. That said, it is always beneficial to incorporate the system's physics into the projector. Hence, from a purely mathematical perspective, we may choose any canonical/convectional projector (eg the $L^2$, $H^1$, or other equivalent projectors) defined in $V$ and we can follow the procedure described in \cite{Shrestha2024ConstructionScales} to find the corresponding Fine-Scale Greens' function and employ a VMS approach to retrieve the projected solution. Whereas, from a purely physical perspective, one would opt to choose a projector that encompasses all the physics of the system which may well include terms that do not strictly conform to mathematical formalisms, such as the properties of a norm. Noting this, we have made a particular choice of employing the symmetric component to define the optimal projector which we view as a balanced compromise of the mathematical and physical perspectives. The defined optimal projector includes the symmetric components of the system which has the property of inducing a well-defined norm in $V$ whereas the non-symmetric terms are excluded as they formally do not induce a norm in $V$. The work presented in \cite{Giraldo2023ANorms} describes a method where the non-symmetric bilinear form is also included in the norm. This approach with the inclusion of the non-symmetric form is proven to be effective for the mesh adaption application considered in \cite{Giraldo2023ANorms}. In light of this, we are prompted to explore deeper questions pertaining to whether our chosen projection is in fact the ``best" choice among the infinite sea of possible projectors or whether an ultimate optimal projection even exists. However, this deviates from the goal of the present study, which seeks to recover the chosen projection of the exact solution considering the canonical projectors detailed in \Cref{sec:optmal_sym}.
\section{Numerical tests}
\label{sec:numerical_ex}
In this section, we present the numerical tests performed to assess the proposed methodology. We performed the tests by considering steady advection-diffusion problems in both direct and mixed formulations described as follows

\begin{minipage}{0.49\linewidth}
    \begin{gather}
        \underline{c} \nabla \phi - \nabla \cdot (\kappa \nabla \phi) = f, \quad \text{in } \Omega = [0, 1]^{d}\label{eq:adv_diff3} \\
        \phi = 0, \quad \text{on } \partial \Omega.
    \end{gather}
\end{minipage}
\begin{minipage}{0.49\linewidth}
    \begin{gather}
        \kappa^{-1} \underline{q} - \nabla \phi = 0, \quad \text{in } \Omega = [0, 1]^{d}\label{eq:adv_diff_mix_3}\\
        \underline{c} \cdot \kappa^{-1} \underline{q} - \nabla \cdot \underline{q} = f, \quad \text{in } \Omega = [0, 1]^{d}\label{eq:adv_diff_mix_4} \\
        \phi = 0, \quad \text{on } \partial \Omega
    \end{gather}
\end{minipage}
\vspace{0.5em}

\noindent In the direct formulation, we seek $\bar{\phi}$ in a finite-dimensional subspace of $H^1_0(\Omega)$ where the solution values are forced to be continuous across the element boundaries and the Dirichlet boundary condition is strongly imposed. In the mixed formulation, we seek the pair $(\underline{\bar{q}}, \bar{\phi})$ in finite-dimensional subspaces $H(\mathrm{div}, \Omega) \times L^2(\Omega)$ where $\bar{\phi}$ is allowed to be discontinuous across the elements and the fluxes $\underline{\bar{q}}$ are forced to be continuous with the Dirichlet boundary condition being weakly imposed. For simplicity, we consider the 1D problem with constant coefficients where we take $c = 1$ and $\nu$ being a positive constant. Furthermore, we consider the 2D case with a constant velocity vector $\underline{c} = \bm{1}$ and a spatially varying anisotropic diffusion tensor $\kappa$ which we construct to be positive definite. We study the 1D problem in both direct and mixed formulation and the 2D problem in mixed formulation. 

\subsection{1D steady advection-diffusion: direct formulation}
We start with the 1D direct formulation for which the variational form of \eqref{eq:adv_diff3} employing the VMS approach with $d = 1$ reads
\begin{gather}
    \bilnear{v^h}{\parddx{\bar{\phi}}{x}} - \bilnear{v^h}{\sigma^{\partial_x}_{SG} \parddx{\bar{\phi}}{x}} + \nu \bilnear{\parddx{v^h}{x}}{\parddx{\bar{\phi}}{x}} = \bilnear{v^h}{f} - \bilnear{v^h}{\sigma^{\partial_x}_{SG} f}.
\end{gather}
We have used the linearity of the inner product to simplify the above equation noting that $c = 1$ and $\kappa = \nu$ are constants. For the numerical tests that follow, we take the diffusion coefficient to be $\nu = 0.01$ (Peclet number of $\alpha = 100$) and constant source term $f = 1$. We can find the exact solution for this case to be given by
\begin{gather}
    \phi_{exact} = \frac{1}{c} \left(x - \frac{e^{\alpha(x - 1)} - e^{-\alpha}}{1 - e^{-\alpha}} \right).
\end{gather}

\begin{figure}[htp]
\begin{subfigure}{0.49\linewidth}
    \centering
    \includegraphics[width = \linewidth]{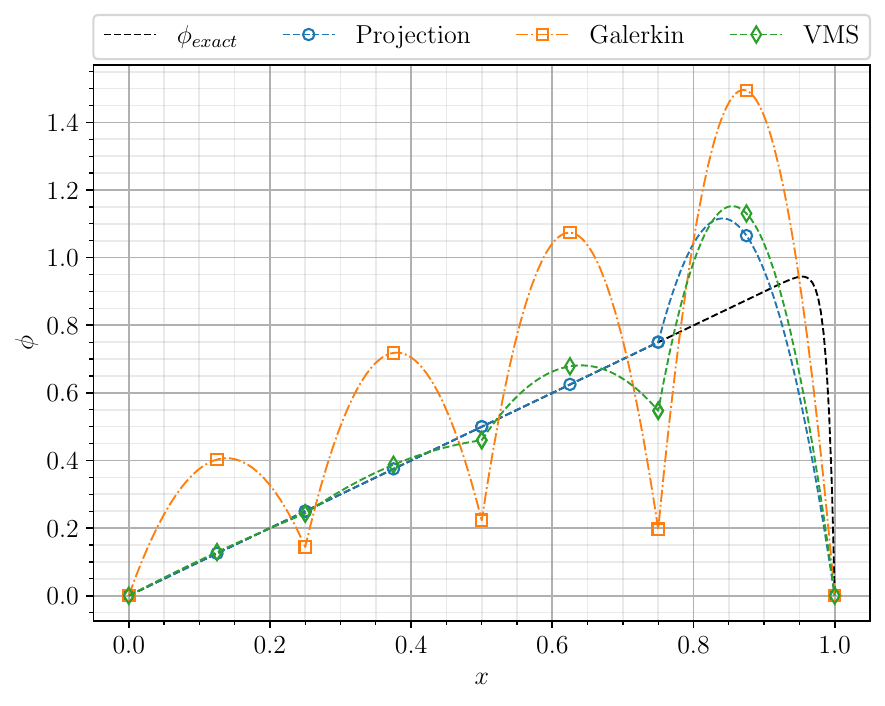}
    \caption{$k = 2$}
    \label{fig:advDiff1D_u_bar_nodal_p2_k2}
\end{subfigure}
\begin{subfigure}{0.49\linewidth}
    \centering
    \includegraphics[width = \linewidth]{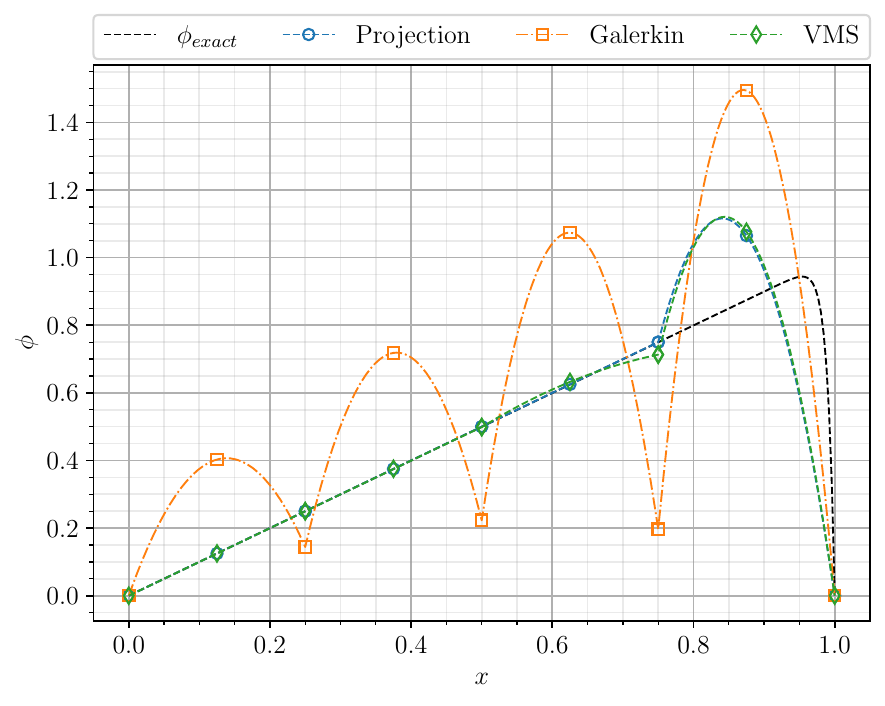}
    \caption{$k = 4$}
    \label{fig:advDiff1D_u_bar_nodal_p2_k4}
\end{subfigure}
\caption{Solutions for the direct formulation of the steady advection-diffusion equation on a coarse mesh with polynomial degree $p = 2$ and $N = 4$ elements for varying the fine mesh refinement parameter $k$}
\label{fig:advDiff1D_u_bar_nodal_p2}
\end{figure}

We start by looking at the solutions on a coarse mesh for varying refinements of the fine mesh in \Cref{fig:advDiff1D_u_bar_nodal_p2}. 
\noindent The two plots in \Cref{fig:advDiff1D_u_bar_nodal_p2} illustrate the solutions to the 1D steady advection-diffusion equation, comparing the traditional Galerkin approach to the proposed VMS approach. Each plot displays the exact solution, the optimal projection onto the coarse mesh, and the Galerkin solution. The traditional Galerkin scheme suffers from being highly oscillatory primarily due to the scheme being an effectively centred scheme. Despite being oscillatory, the Galerkin solution is stable, which is an important point to note. The plot in \Cref{fig:advDiff1D_u_bar_nodal_p2_k2} depicts the solution of the VMS approach with a refinement parameter $k=2$, which demonstrates a significant improvement over the Galerkin solution, yielding a solution that closely aligns with the optimal projection. The second plot in \Cref{fig:advDiff1D_u_bar_nodal_p2_k4} further explores this VMS approach with a higher refinement parameter $k=4$, resulting in a solution that even more accurately approximates the optimal projection.

\begin{figure}[htp]
\begin{subfigure}{0.49\linewidth}
    \centering
    \includegraphics[width = \linewidth]{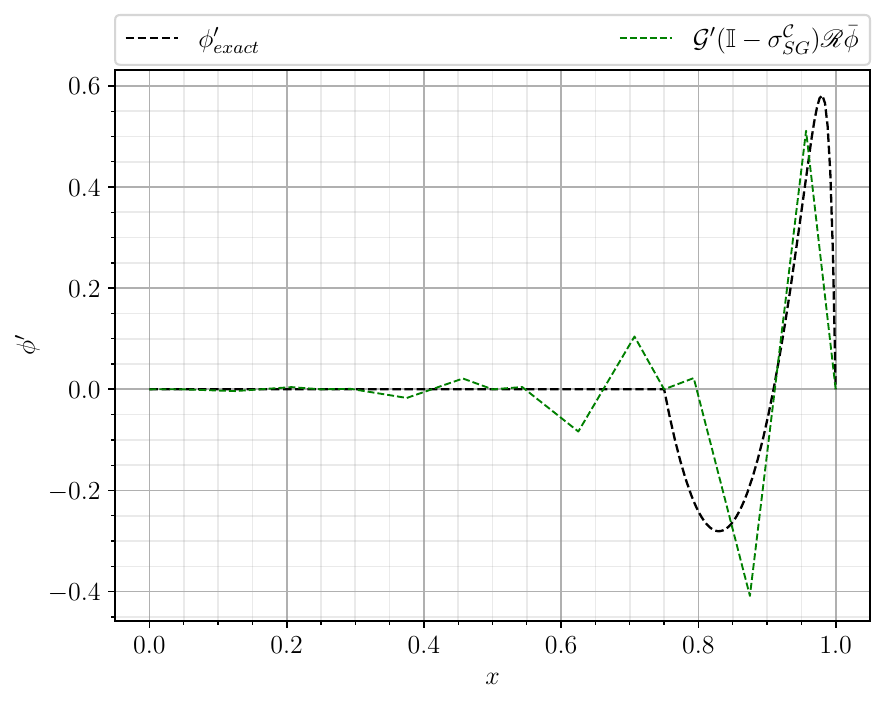}
    \caption{$k = 2$}
    \label{fig:advDiff1D_u_prime_nodal_p2_k2}
\end{subfigure}
\begin{subfigure}{0.49\linewidth}
    \centering
    \includegraphics[width = \linewidth]{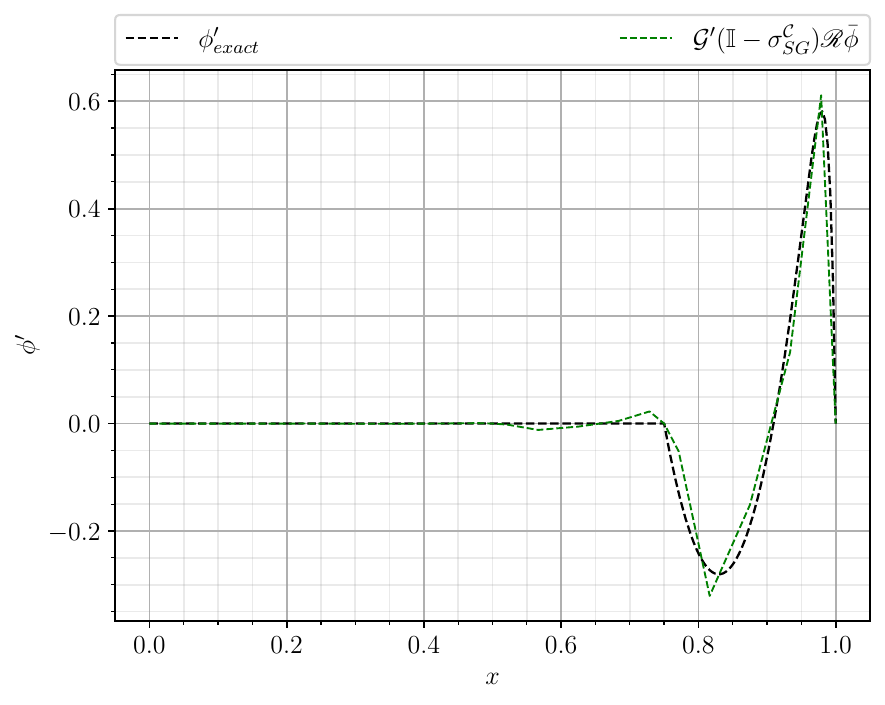}
    \caption{$k = 4$}
    \label{fig:advDiff1D_u_prime_nodal_p2_k4}
\end{subfigure}
\caption{Unresolved scales for the direct formulation of the steady advection-diffusion equation on a coarse mesh with polynomial degree $p = 2$ and $N = 4$ elements for varying the fine mesh refinement parameter $k$}
\label{fig:advDiff1D_u_prime_nodal_p2}
\end{figure}

Next, we consider the plots for the unresolved scales corresponding to the discussed cases computed using the Fine-Scale Greens' function. The plots in \Cref{fig:advDiff1D_u_prime_nodal_p2} show the exact unresolved scales corresponding to the optimal projection onto a coarse mesh with polynomial degree $p = 2$ with $N = 4$ elements along with the computed fine scales for varying $k$. \Cref{fig:advDiff1D_u_prime_nodal_p2_k2} shows the computed fine scales for $k = 2$. Note that the computed fine scales only roughly align with the exact one and the general prediction is quite poor. In the refined case with $k = 4$ shown in \Cref{fig:advDiff1D_u_prime_nodal_p2_k4}, we see that the predictions are improved and the computed fine scales more closely match the exact one, although not exactly. The mismatch in the unresolved scales is reflected in the small mismatch between the VMS solution and the exact projection shown in \Cref{fig:advDiff1D_u_bar_nodal_p2}. 

We now move on to assess the $h$-$p$ convergence of the schemes. As the optimal projector in the direct formulation is based on the norm equivalent to $H^1$, we measure the error using the $H^1$ norm as follows
\begin{align}
    error \: wrt \: the  \: exact \: solution : \: &\left(\lVert \bar{\phi} - \phi \rVert_{L^2}^2 + \lVert \nabla \bar{\phi} - \nabla \phi \rVert_{L^2}^2 \right)^{\frac{1}{2}} \\
    error \: wrt \: the   \: exact \: projection : \: &\left(\lVert \bar{\phi} - \mathcal{P}\phi \rVert_{L^2}^2 + \lVert \nabla \bar{\phi} - \nabla \mathcal{P} \phi \rVert_{L^2}^2\right)^{\frac{1}{2}} \\
    error \: wrt \: the  \: exact \: fine \: scales : \: &\left(\lVert \phi^{\prime}_k - \phi^{\prime} \rVert_{L^2}^2 + \lVert \nabla \phi^{\prime}_k - \nabla \phi^{\prime} \rVert_{L^2}^2\right)^{\frac{1}{2}},
\end{align}
where $\phi$ represents the exact solution, $\mathcal{P}\phi$ denotes the exact projection, $\phi^{\prime}$ signifies the exact unresolved scales, and $\bar{\phi}$ and $\phi^{\prime}_k$ are the computed resolved and unresolved scales, respectively. The computation of all the error norms is done using a Gauss-Lobatto quadrature rule with a degree of precision of 25 to ensure the exact solutions are properly captured.

\begin{figure}[htp]
\begin{subfigure}{0.49\linewidth}
    \centering
    \includegraphics[width = \linewidth]{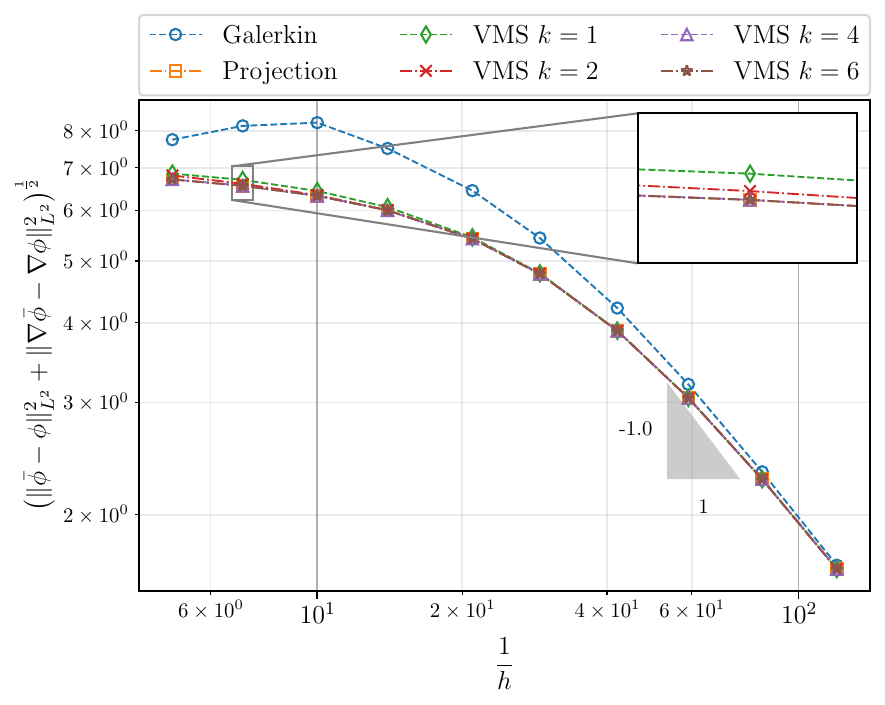}
    \caption{$p = 1$}
    \label{fig:advDiff1D_convergence_nodal_p1}
\end{subfigure}
\begin{subfigure}{0.49\linewidth}
    \centering
    \includegraphics[width = \linewidth]{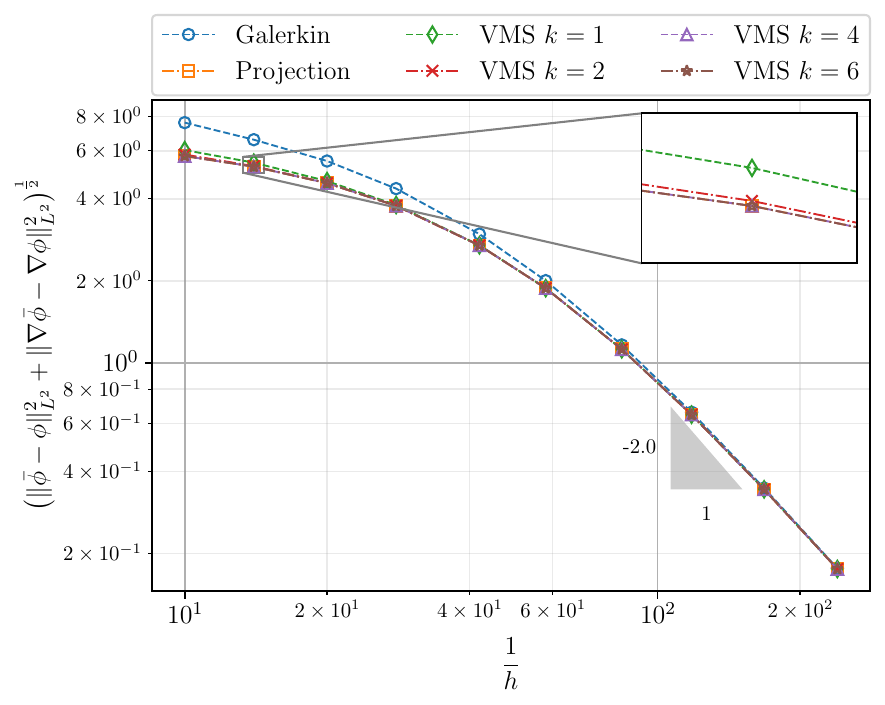}
    \caption{$p = 2$}
    \label{fig:advDiff1D_convergence_nodal_p2}
\end{subfigure} \\
\begin{subfigure}{0.49\linewidth}
    \centering
    \includegraphics[width = \linewidth]{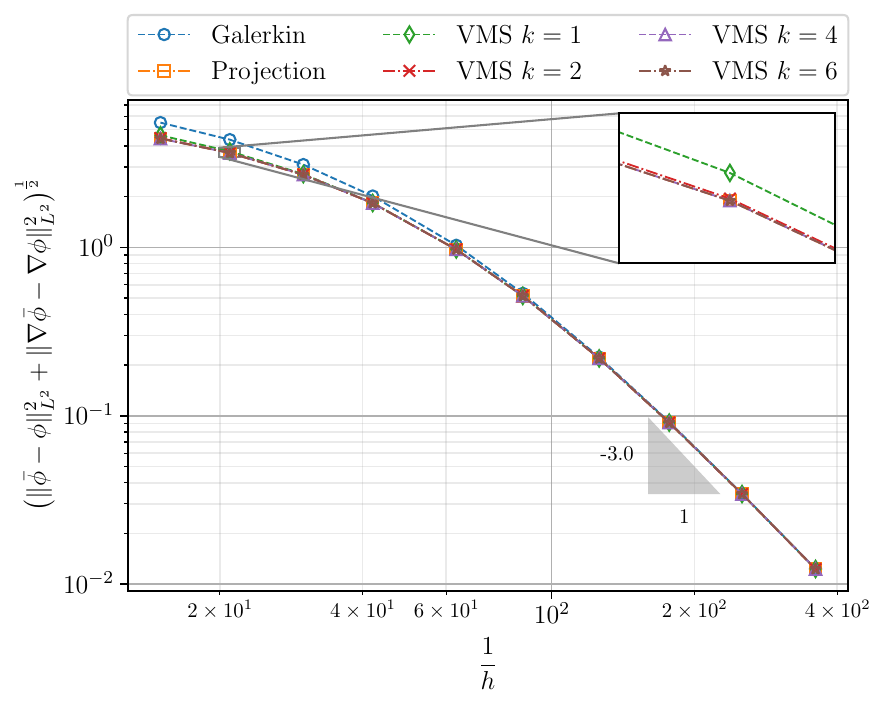}
    \caption{$p = 3$}
    \label{fig:advDiff1D_convergence_nodal_p3}
\end{subfigure}
\begin{subfigure}{0.49\linewidth}
    \centering
    \includegraphics[width = \linewidth]{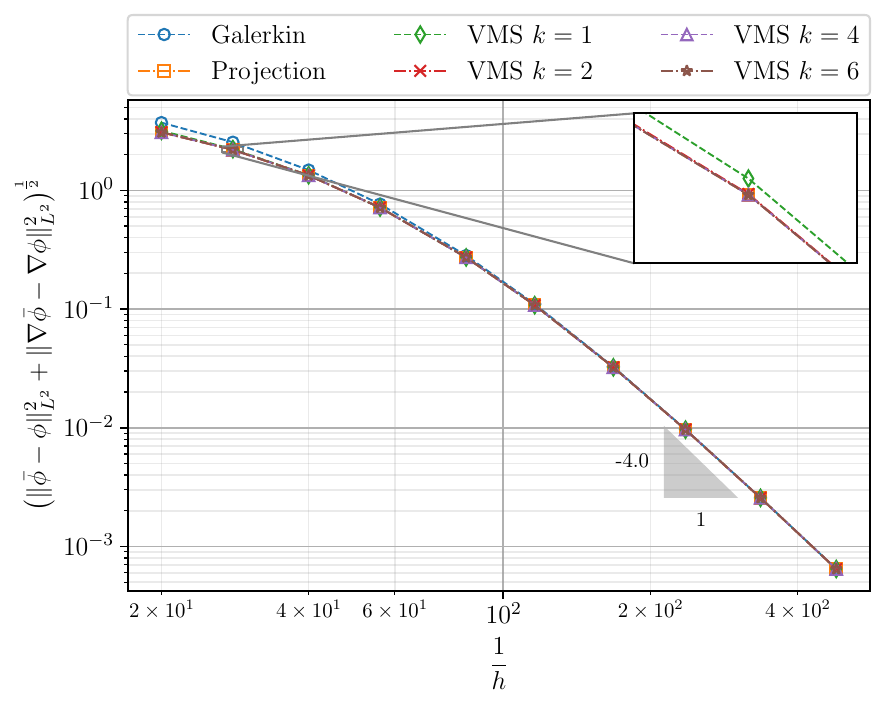}
    \caption{$p = 4$}
    \label{fig:advDiff1D_convergence_nodal_p4}
\end{subfigure}
\caption{$h$-$p$ convergence for the direct formulation of the steady advection-diffusion equation}
\label{fig:advDiff1D_convergence_nodal}
\end{figure}

We begin by examining the $h$-$p$ convergence of the resolved scales, as illustrated in \Cref{fig:advDiff1D_convergence_nodal}. This figure depicts the convergence behaviour of the Galerkin solution, the exact projection, and the VMS solution for different values of $k$. We observe that on very coarse meshes, the Galerkin solution possesses the largest error, while the projection, by design, has the lowest error. The error of the VMS solution falls between that of the Galerkin solution and the projection, but it rapidly converges towards the projection as the fine mesh refinement parameter $k$ is increased. Overall, as the coarse mesh is refined, all the schemes converge at the expected rate of $p$ in the $H^1$ error norm.
\begin{figure}[htp]
\begin{subfigure}{0.49\linewidth}
    \centering
    \includegraphics[width = \linewidth]{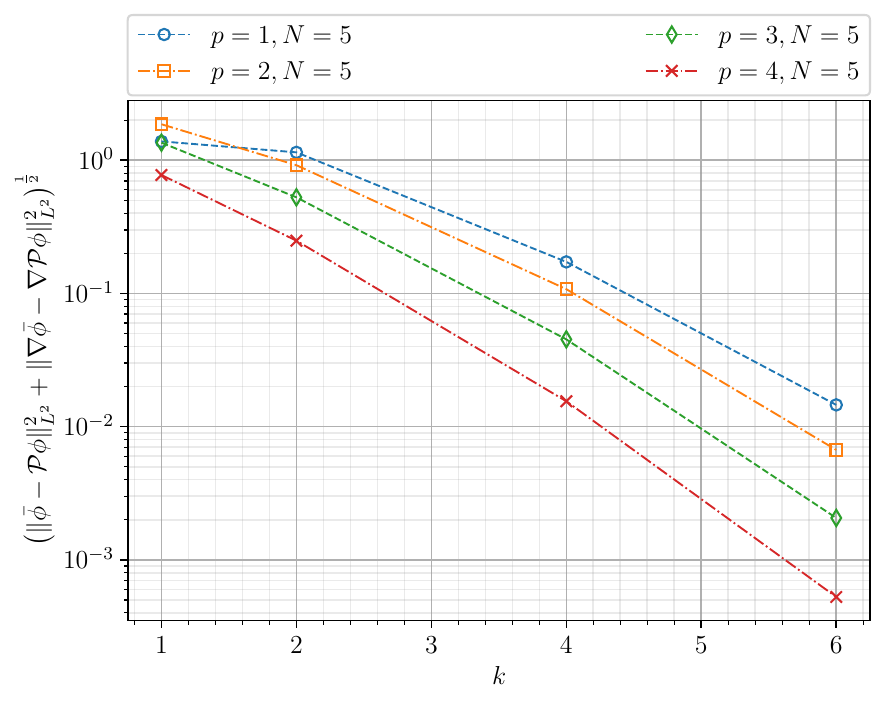}
    \caption{$k$-convergence of the VMS solution}
    \label{fig:advDiff1D_convergence_nodal_u_bar_k}
\end{subfigure}
\begin{subfigure}{0.49\linewidth}
    \centering
    \includegraphics[width = \linewidth]{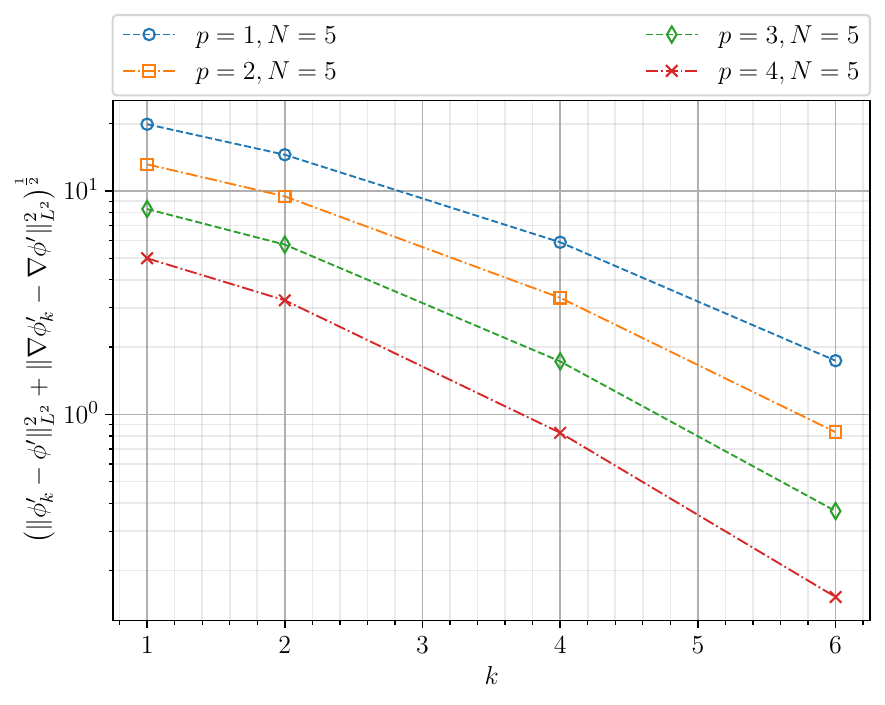}
    \caption{$k$-convergence of the unresolved scales}
    \label{fig:advDiff1D_convergence_nodal_u_prime_k}
\end{subfigure}
\caption{$k$-convergence VMS solution and unresolved scales for the direct formulation of the 1D steady advection-diffusion equation on coarse meshes with varying polynomial degree}
\label{fig:advDiff1D_convergence_nodal_u_bar_u_prime}
\end{figure}
We observe a similar convergence behaviour for the convergence of the resolved scales to the exact projection and the convergence of the fine scales as depicted in \Cref{fig:advDiff1D_convergence_nodal_u_bar_u_prime}. The first significant observation from \Cref{fig:advDiff1D_convergence_nodal_u_bar_u_prime} is the exponential convergence achieved with $k$ refinement. An exception is noted on the coarse mesh with $p=1$, $N=5$, and $k=1$, where a local super-convergence is observed, as seen through \Cref{fig:advDiff1D_convergence_nodal_u_bar_k}. This super-convergence is specific to the test case considered and does not appear in more general cases. Furthermore, when comparing the error values of the resolved scales with respect to the projection in \Cref{fig:advDiff1D_convergence_nodal_u_bar_k} with the error values for the unresolved scales in \Cref{fig:advDiff1D_convergence_nodal_u_prime_k}, it is evident that for a given coarse mesh, the error in the unresolved scales is consistently larger than the error in the resolved scales relative to the projection. This indicates that a highly accurate estimation of the fine scales is not essential for obtaining a good approximation of the projection. This can be related back to the error analysis from \Cref{sec:greens_func} where we noted that the bilinear form $c(\bar{\phi} - v^h, \phi^{\prime}_k  - \phi^{\prime})$ dictates whether we achieve optimality. As observed through \Cref{fig:advDiff1D_convergence_nodal_u_prime_k}, we see that the computed unresolved scales may not closely match the exact ones, however, their contribution in the bilinear form is evidently weak which thus yields a numerical solution that closely matches the exact projection. This fact can also attributed to the fact that we work in Sobolev spaces, where pointwise values are less relevant and integrals are of primary importance.
\subsection{1D steady advection-diffusion: mixed formulation}
\begin{figure}[htp]
\begin{subfigure}{0.49\linewidth}
    \centering
    \includegraphics[width = \linewidth]{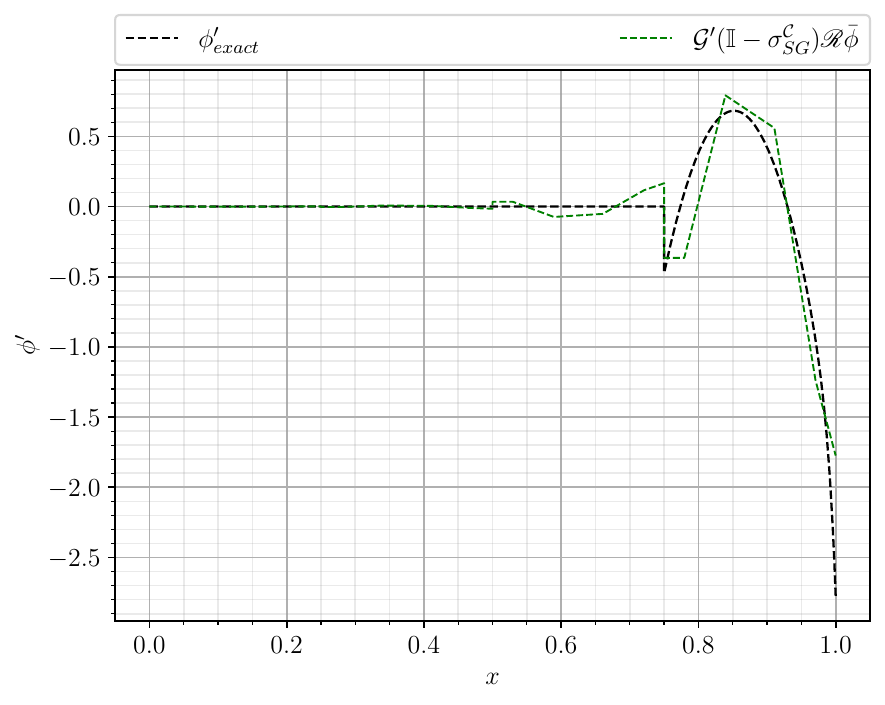}
    \caption{$k = 2$}
    \label{fig:advDiff1D_u_prime_edge_p3_k2}
\end{subfigure}
\begin{subfigure}{0.49\linewidth}
    \centering
    \includegraphics[width = \linewidth]{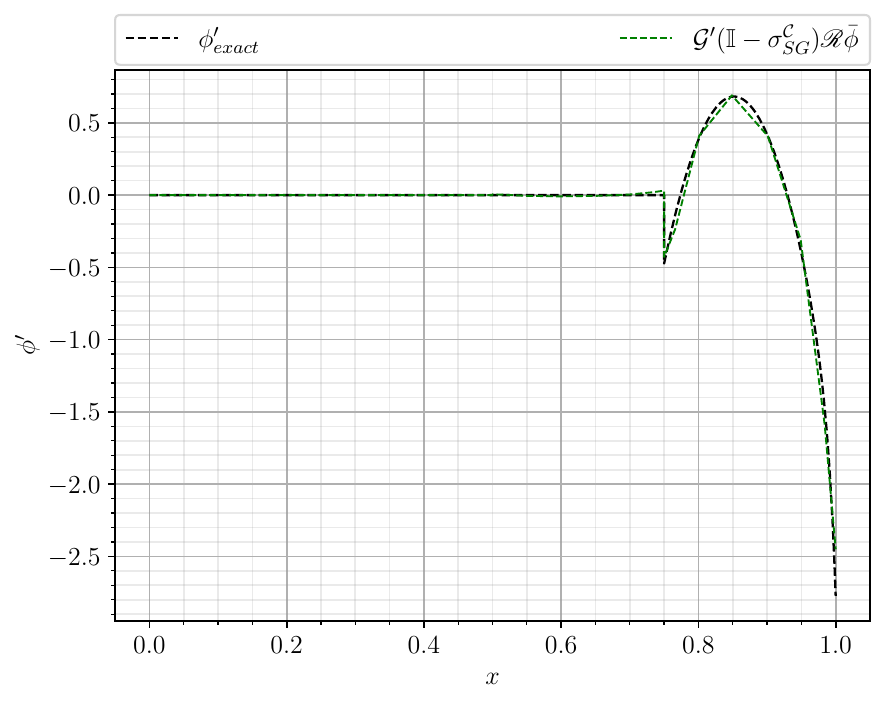}
    \caption{$k = 4$}
    \label{fig:advDiff1D_u_prime_edge_p3_k4}
\end{subfigure}
\caption{Unresolved scales for the mixed formulation of the steady advection-diffusion equation on a coarse mesh with polynomial degree $p = 3$ and $N = 4$ elements for varying the fine mesh refinement parameter $k$}
\label{fig:advDiff1D_u_prime_edge_p3}
\end{figure}

Moving on from the direct formulation, we now consider the 1D advection-diffusion equation in the mixed formulation. Applying the simplification with the constant advection speed and the constant diffusion coefficient, the variational form of the problem reads
\begin{gather}
    \bilnear{v^h}{\bar{q}} + \bilnear{\parddx{v^h}{x}}{\bar{\phi}} = 0, \quad \forall v^h \in \bar{V}  \\
    \begin{split}
        \bilnear{\eta^h}{\bar{q}} - \bilnear{\eta^h}{\sigma_{SG}^{\partial_x} \bar{q}} - &\nu \bilnear{\eta^h}{\parddx{\bar{q}}{x}} = \\ &\bilnear{\eta^h}{f} - \bilnear{\eta^h}{\sigma_{SG}^{\partial_x} f}, \quad \forall \eta^h \in \bar{W}.
    \end{split}
\end{gather}
We consider the same test case as for the direct formulation with $f = 1$ and $\nu = 0.01$. 
\begin{figure}[htp]
\begin{subfigure}{0.49\linewidth}
    \centering
    \includegraphics[width = \linewidth]{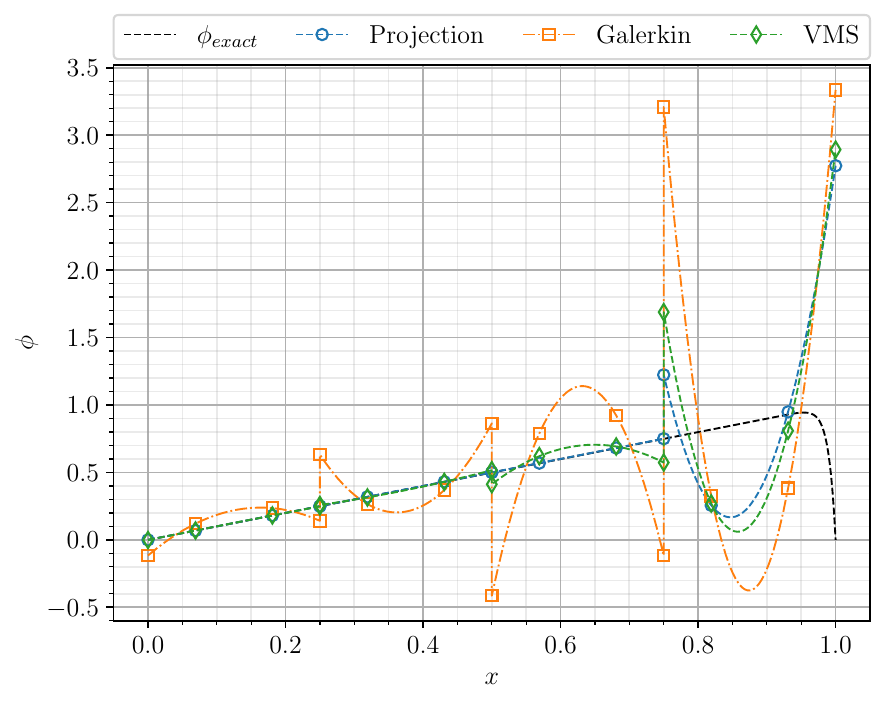}
    \caption{$k = 2$}
    \label{fig:advDiff1D_u_bar_edge_p3_k2}
\end{subfigure}
\begin{subfigure}{0.49\linewidth}
    \centering
    \includegraphics[width = \linewidth]{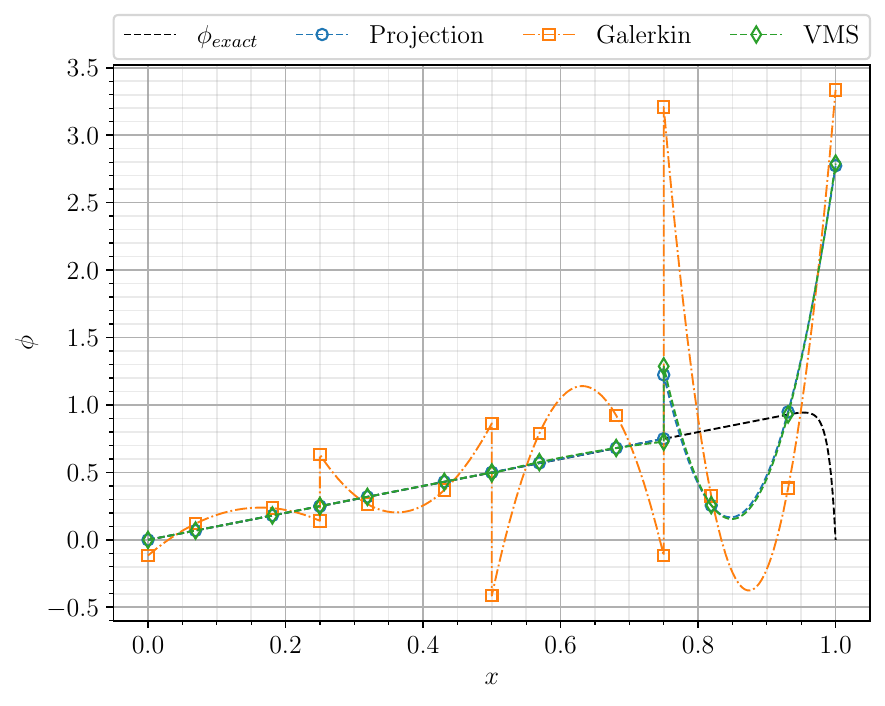}
    \caption{$k = 4$}
    \label{fig:advDiff1D_u_bar_edge_p3_k4}
\end{subfigure}
\caption{Solutions for the mixed formulation of the steady advection-diffusion equation on a coarse mesh with polynomial degree $p = 3$ and $N = 4$ elements for varying the fine mesh refinement parameter $k$}
\label{fig:advDiff1D_u_bar_edge_p3}
\end{figure}
In line with expectations, for the mixed formulation, we observe very similar behaviour to that of the direct formulation. The plots in \Cref{fig:advDiff1D_u_bar_edge_p3} depict the solutions to the 1D steady advection-diffusion equation now for the mixed formulation. Once again, the Galerkin scheme produces a solution that is stable but polluted with oscillations. The VMS solution for $k = 2$ shown in \Cref{fig:advDiff1D_u_bar_edge_p3_k2} shows improvement over the Galerkin solution where the solution closely matches the optimal projection. Moreover, the VMS solution for $k = 4$ comes even closer to the optimal projection with the two curves nearly lying on top of each other as seen through \Cref{fig:advDiff1D_u_bar_edge_p3_k4}.
The same observations can be made for \Cref{fig:advDiff1D_u_prime_edge_p3} where the unresolved scales are only roughly estimated for $k = 2$ and it ends up closely matching the exact unresolved scales for $k = 4$.

\begin{figure}[H]
\begin{subfigure}{0.49\linewidth}
    \centering
    \includegraphics[width = \linewidth]{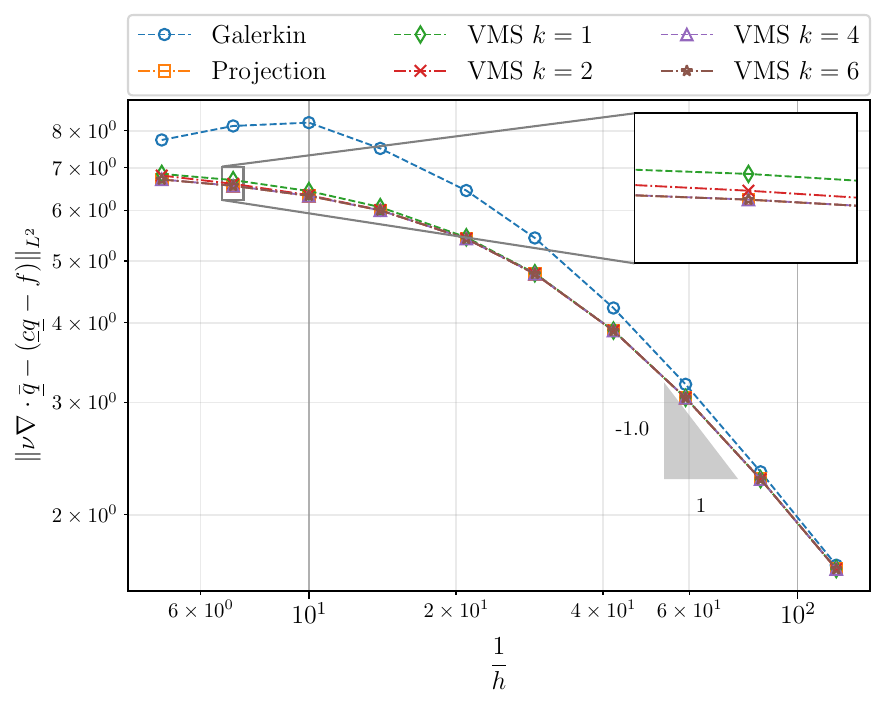}
    \caption{$p = 1$}
    \label{fig:advDiff1D_convergence_edge_p1}
\end{subfigure}
\begin{subfigure}{0.49\linewidth}
    \centering
    \includegraphics[width = \linewidth]{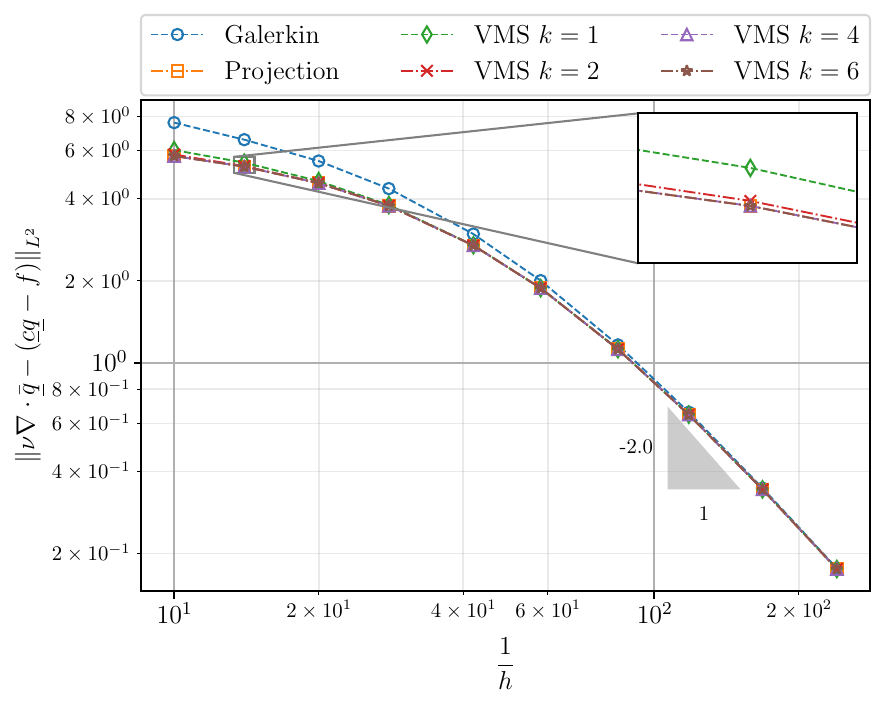}
    \caption{$p = 2$}
    \label{fig:advDiff1D_convergence_edge_p2}
\end{subfigure} \\
\begin{subfigure}{0.49\linewidth}
    \centering
    \includegraphics[width = \linewidth]{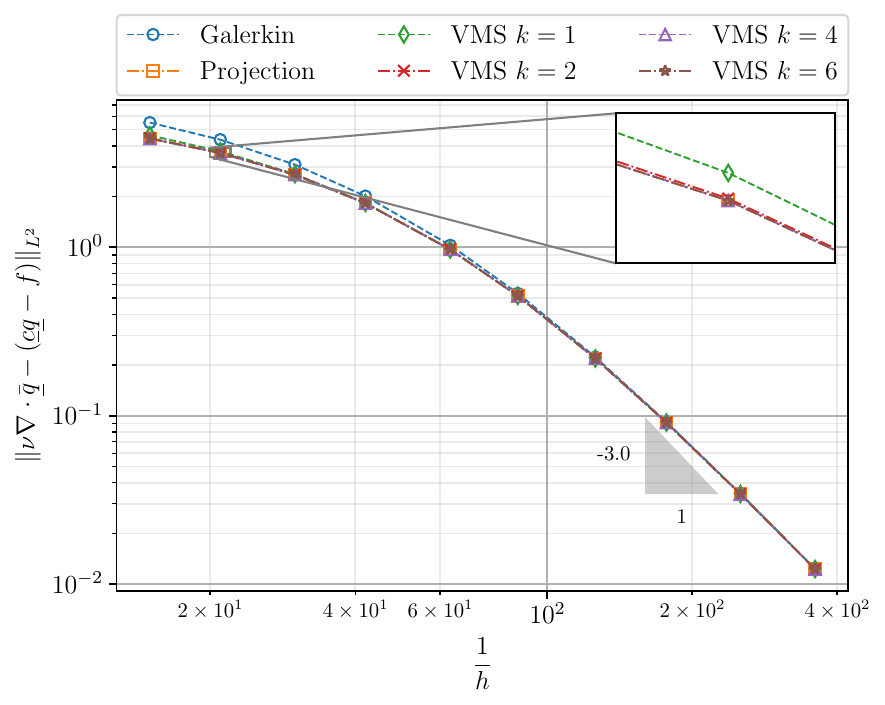}
    \caption{$p = 3$}
    \label{fig:advDiff1D_convergence_edge_p3}
\end{subfigure}
\begin{subfigure}{0.49\linewidth}
    \centering
    \includegraphics[width = \linewidth]{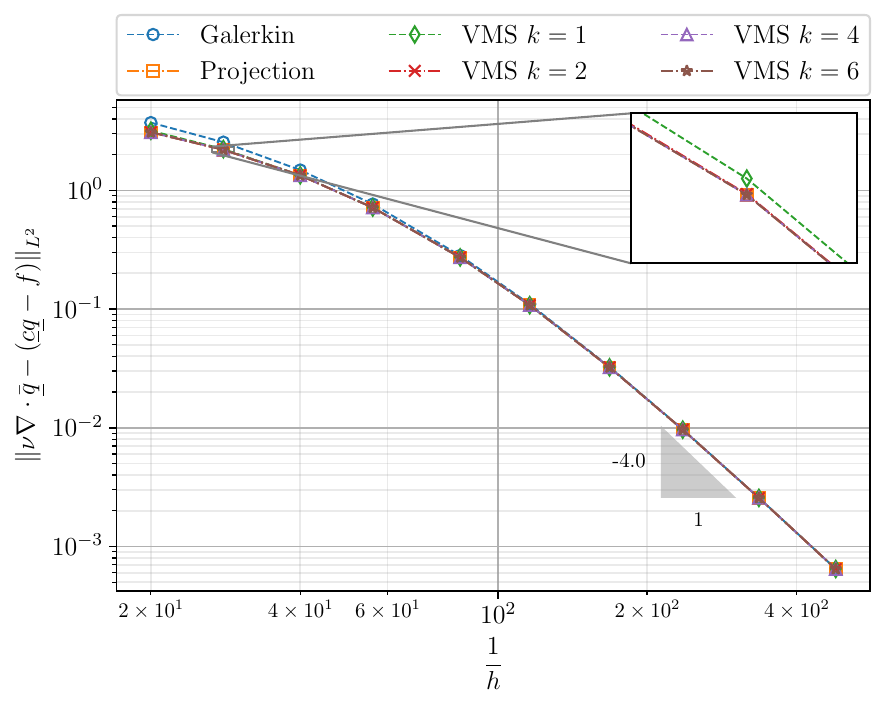}
    \caption{$p = 4$}
    \label{fig:advDiff1D_convergence_edge_p4}
\end{subfigure}
\caption{$h$-$p$ convergence for the mixed formulation of the steady advection-diffusion equation}
\label{fig:advDiff1D_convergence_edge}
\end{figure}

We now assess the convergence of the different schemes as was done for the direct formulation. We stick with the $H^1$ error norm for computing the error between the VMS solution and the exact projection, and the error between computed and exact fine scales. For the error with respect to the exact solution, on the other hand, we use a different norm. Normally, we would construct this error norm based on the minimisation problem associated with the optimal projector to compute the error respect to the exact solution. However, we have a saddle point problem for the mixed formulation instead of a minimisation problem as highlighted in \eqref{eq:mix_proj}. Hence, the optimal projector does not minimise a functional but rather compromises the error in the solution gradients while ensuring that the divergence of the numerical solution matches the exact divergence. While there is a possibility to encapsulate this notion in a norm equivalent to the $H(\mathrm{div})$ norm by considering the operator norm associated with the projector, we exclude this analysis in the present paper. Instead, we opt to measure the error with respect to the exact solution using the $H(\mathrm{div})$ semi-norm, with the argument that the optimal projector prioritises the divergence error yielding the optimal solution whose divergence matches the exact divergence in $L^2$. We thus have the following error norms
\begin{align}
    error \: wrt \: the  \: exact \: solution : \: &\lVert \nu \nabla \cdot \bar{\underline{q}} - (\underline{c} \underline{q} - f) \rVert_{L^2} \label{eq:err_norm_mix_1} \\
    error \: wrt \: the  \: exact \: projection : \: &\left(\lVert \bar{\phi} - \mathcal{P}\phi \rVert_{L^2}^2 + \lVert {\underline{\bar{q}}} - \mathcal{P} \underline{q} \rVert_{L^2}^2\right)^{\frac{1}{2}} \label{eq:err_norm_mix_2} \\
    error \: wrt \: the  \: exact \: fine \: scales : \: &\left(\lVert \phi^{\prime}_k - \phi^{\prime} \rVert_{L^2}^2 + \lVert \underline{q}^{\prime}_k - \underline{q}^{\prime} \rVert_{L^2}^2\right)^{\frac{1}{2}}, \label{eq:err_norm_mix_3}
\end{align}
where $\underline{q}$, $\mathcal{P}\phi$, $\mathcal{P}\underline{q}$, and $\phi^{\prime}$ denote the exact solution gradient, the exact projection of the solution and its gradient, and the exact unresolved scales, respectively. Additionally, $\bar{\phi}$ and $\phi^{\prime}_k$ denote the computed resolved and unresolved scales, respectively. Once again we use over-integration to evaluate the integrals for the error norm using a Gauss-Lobatto quadrature with a degree of precision of 25. Assessing the convergence of the schemes using the specified error norms reveals the anticipated behaviour, namely, the Galerkin solution possesses the highest errors, the optimal projection achieves the lowest error, and the VMS solutions fall in between the two as seen through \Cref{fig:advDiff1D_convergence_edge}. The VMS solution rapidly collapses to the optimal projection with $k$ refinement, and all the schemes converge at the expected rate of $p$ once the mesh is sufficiently refined. 

\begin{figure}[H]
\begin{subfigure}{0.49\linewidth}
    \centering
    \includegraphics[width = \linewidth]{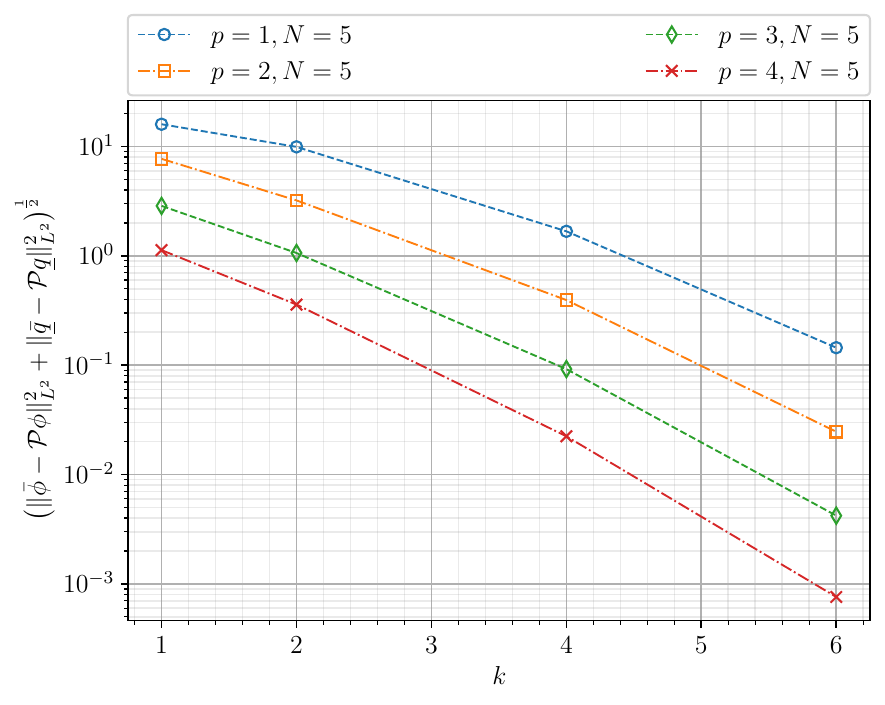}
    \caption{$k$-convergence of the VMS solution}
    \label{fig:advDiff1D_convergence_edge_u_bar_k}
\end{subfigure}
\begin{subfigure}{0.49\linewidth}
    \centering
    \includegraphics[width = \linewidth]{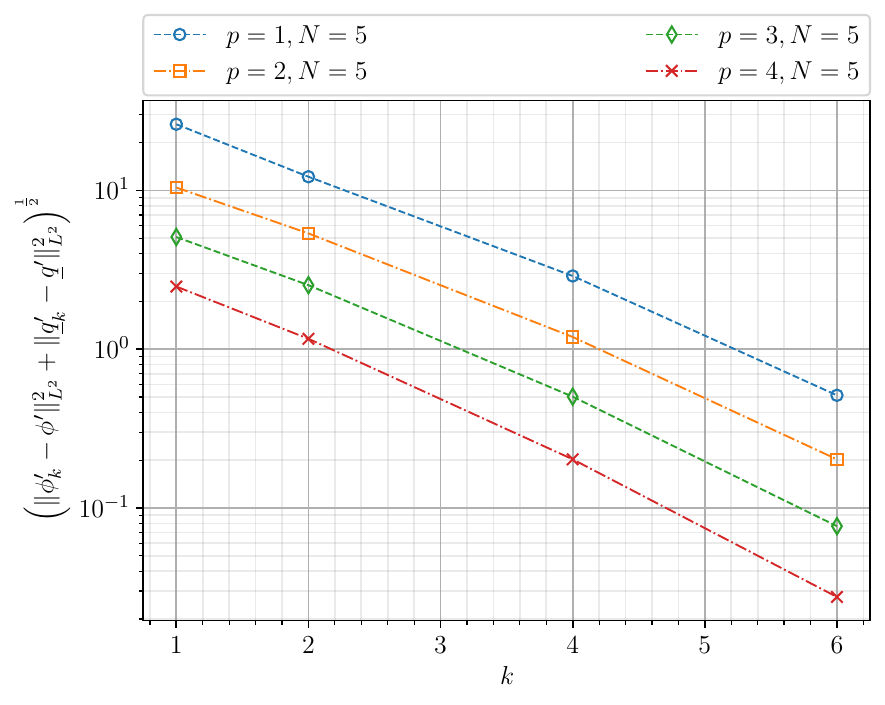}
    \caption{$k$-convergence of the unresolved scales}
    \label{fig:advDiff1D_convergence_edge_u_prime_k}
\end{subfigure}
\caption{$k$-convergence VMS solution and unresolved scales for the mixed formulation of the 1D steady advection-diffusion equation on coarse meshes with varying polynomial degree}
\label{fig:advDiff1D_convergence_edge_u_bar_u_prime}
\end{figure}

Concerning the convergence of the VMS solution and the unresolved scales to the exact projection and the exact unresolved scales, we observe exponential convergence with $k$ as depicted in \Cref{fig:advDiff1D_convergence_edge_u_bar_u_prime}. Like the case for the direct formulation, the error for the unresolved scales is significantly larger than that for the resolved scales as seen when comparing the error values in \Cref{fig:advDiff1D_convergence_edge_u_bar_k} and \Cref{fig:advDiff1D_convergence_edge_u_prime_k}. Hence the statement that a highly accurate estimation of the fine scales is not essential for obtaining a good approximation of the projection still holds.

\subsection{2D steady advection-diffusion: mixed formulation}
\begin{figure}[htp]
    \begin{minipage}{0.44\linewidth}
        \includegraphics[width=\linewidth]{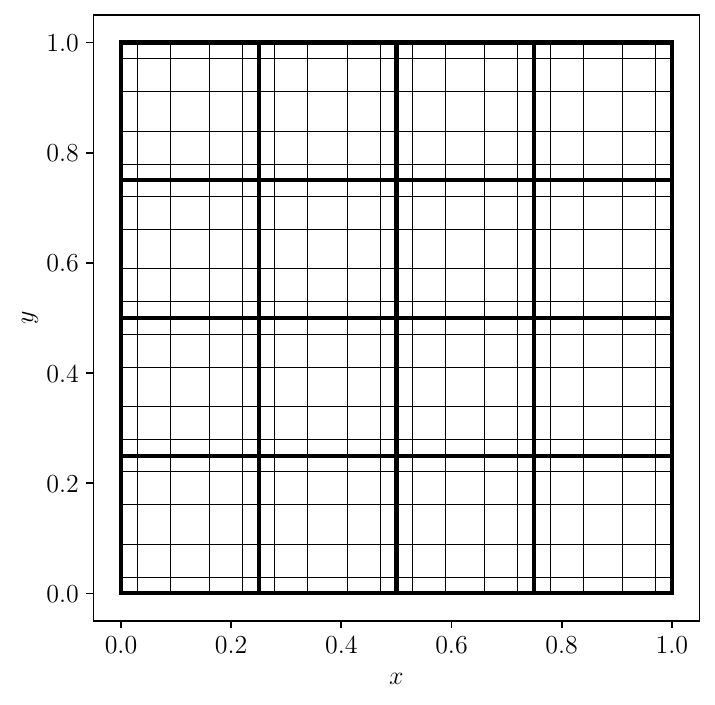}
        \caption{Orthogonal mesh with $4 \times 4$ elements}
        \label{fig:mesh_ortho}
    \end{minipage}
    \hfill
    \begin{minipage}{0.44\linewidth}
        \includegraphics[width=\linewidth]{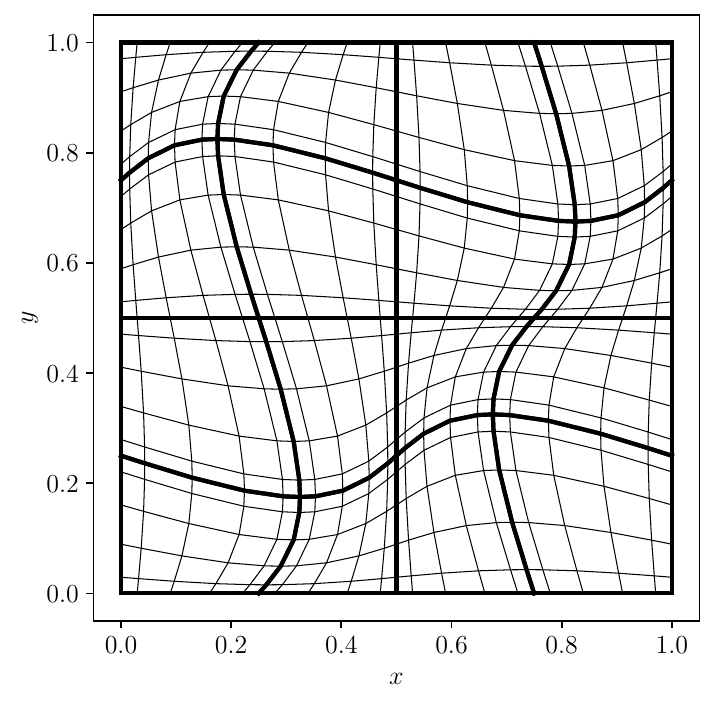}
        \caption{Curvilinear mesh with $4 \times 4$ elements}
        \label{fig:mesh_skew}
    \end{minipage}
\end{figure}

In this section, we demonstrate the extension of the proposed formulation to the multidimensional setting by considering the steady advection-diffusion problem in 2D. We consider the variational form shown in \eqref{eq:weak_form_mix1} and \eqref{eq:weak_form_mix2} and we perform the initial set of tests using the following advection speed and the diffusion tensor
\begin{gather}
    \underline{c} = \boldsymbol{1} = [1, 1]^T, \quad \quad \kappa = \left[\begin{array}{cc}
        \mu_x(x, y) & \frac{1}{2} \sqrt{\mu_x(x, y) \mu_y(x, y)} \\
       \frac{1}{2} \sqrt{\mu_x(x, y) \mu_y(x, y)}  & \mu_y(x, y)
    \end{array} \right] \\
    \text{where } \mu_x(x, y) = \nu (1 + (1 - \varepsilon) \sin{(2 \pi x) \cos{(2 \pi y}})),\quad  \mu_y(x, y) = \nu (1 + (1 - \varepsilon) \cos{(2 \pi x) \sin{(2 \pi y}})). \label{eq:diff_tensor}
\end{gather}
For these 2D test cases, we employ a manufactured solution procedure with the manufactured solution given by
\begin{gather}
    \phi_{exact} = \left(x - \frac{e^{\alpha(x - 1)} - e^{-\alpha}}{1 - e^{-\alpha}} \right) \left(y - \frac{e^{\alpha(y - 1)} - e^{-\alpha}}{1 - e^{-\alpha}} \right), \label{eq:manufactured_soln}
\end{gather}
which generates sharp layers in the domain for large values of the Peclet number $\alpha$. Moreover, we perform the tests on an orthogonal mesh and a curvilinear mesh on the domain $[0, 1]^2$. The mapping $\boldsymbol{\Phi} \: : \: (\xi, \eta) \rightarrow (x, y)$ from the reference domain $(\xi, \eta) \in [-1, 1]^2$ to the physical domain $(x, y) \in [0, 1]^2$ for the curvilinear mesh is given by
\begin{equation}
\left\{\begin{array} { l } 
{ \Phi_x = \hat { x } + 0.075 \sin{( 2 \pi \hat { x } )} \sin {( 2 \pi \hat { y } )}} \\
{ \Phi_y = \hat { y } - 0.075 \sin{( 2 \pi \hat { x } )} \sin {( 2 \pi \hat { y } )}}
\end{array} \text {, where } \left\{\begin{array}{l}
\hat{x}=0.5(1+\xi) \\
\hat{y}=0.5(1+\eta)
\end{array}\right.\right. \text {. }
\end{equation}
Examples of these meshes are depicted \Cref{fig:mesh_ortho} and \Cref{fig:mesh_skew}.

\begin{figure}[htp]
\begin{subfigure}{0.33\linewidth}
    \centering
    \includegraphics[width = \linewidth]{Images/advDiff2D_anisoV2/u_bar_ex_surf_Glk_nu_0.010_skew_0.0_N_4_p_3.pdf}
    \caption{Projection}
\end{subfigure}
\begin{subfigure}{0.33\linewidth}
    \centering
    \includegraphics[width = \linewidth]{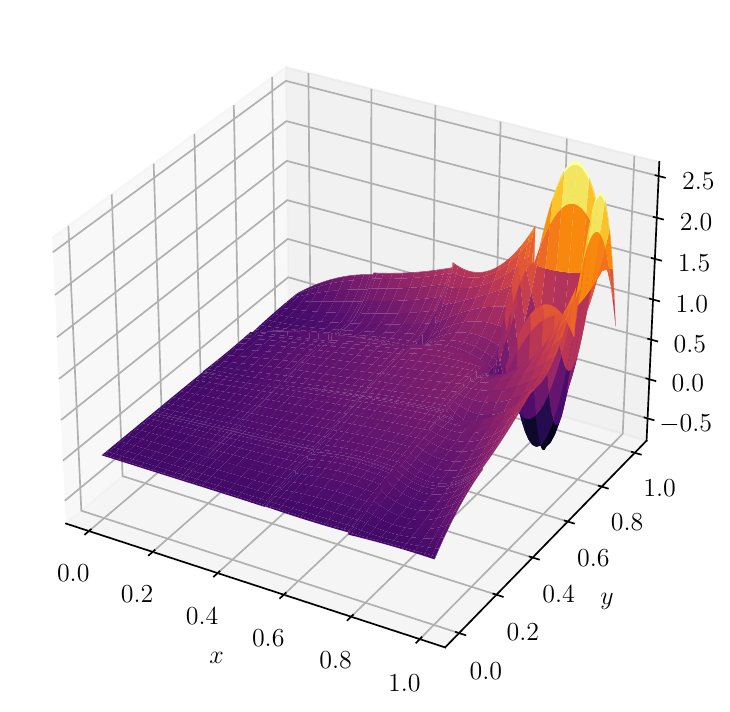}
    \caption{VMS solution with $k = 2$}
    \label{fig:advDiff2D_u_bar_VMS_p2_k1}
\end{subfigure}
\begin{subfigure}{0.33\linewidth}
    \centering
    \includegraphics[width = \linewidth]{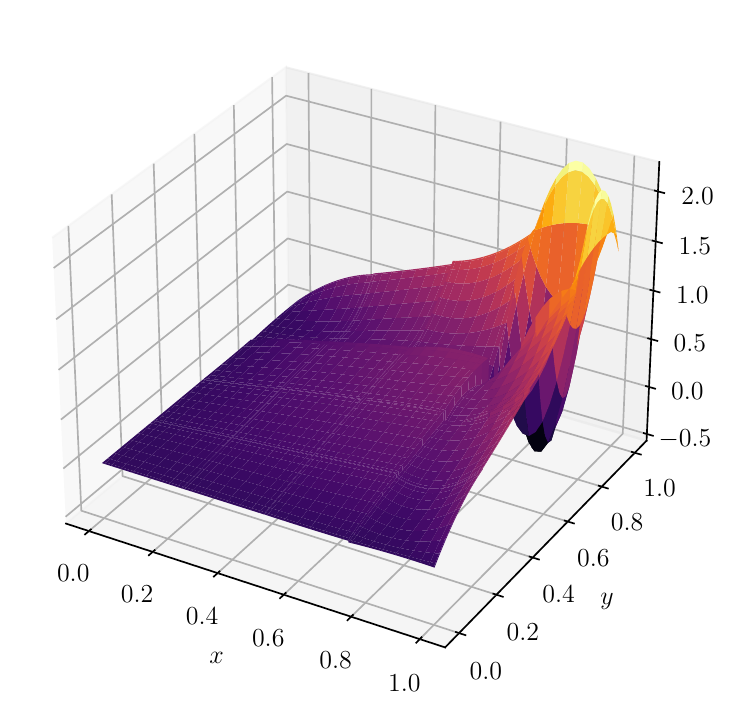}
    \caption{VMS solution with $k = 4$}
    \label{fig:advDiff2D_u_bar_VMS_p2_k4}
\end{subfigure}
\caption{Optimal projection of the exact solution and VMS solutions with varying $k$ on a coarse orthogonal mesh with polynomial degree {$p = 3$} and $4 \times 4$ elements}
\label{fig:advDiff2D_u_bar_VMS_p2}
\end{figure}

\begin{figure}[H]
\begin{subfigure}{0.33\linewidth}
    \centering
    \includegraphics[width = \linewidth]{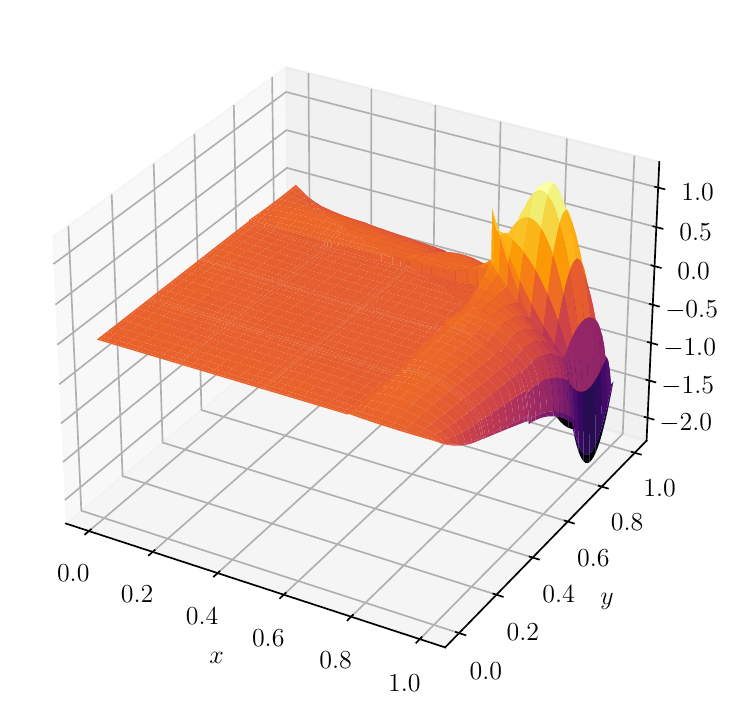}
    \caption{Exact unresolved scales}
\end{subfigure}
\begin{subfigure}{0.33\linewidth}
    \centering
    \includegraphics[width = \linewidth]{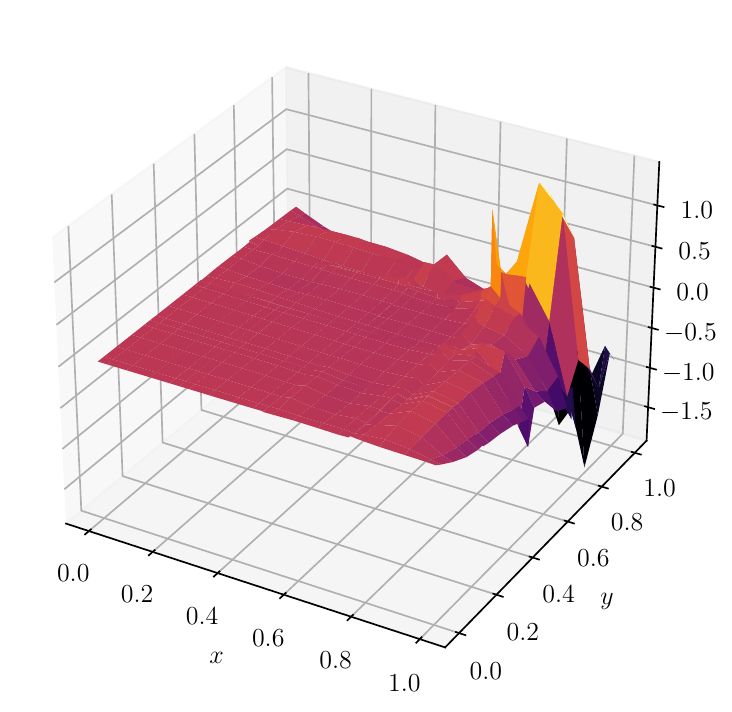}
    \caption{$\mathcal{G}^{\prime}(\mathbb{I} - \sigma_{SG}^{\mathcal{C}})\mathscr{R}\bar{\phi}$ with $k = 2$}
    \label{fig:advDiff2D_u_prime_VMS_p2_k1}
\end{subfigure}
\begin{subfigure}{0.33\linewidth}
    \centering
    \includegraphics[width = \linewidth]{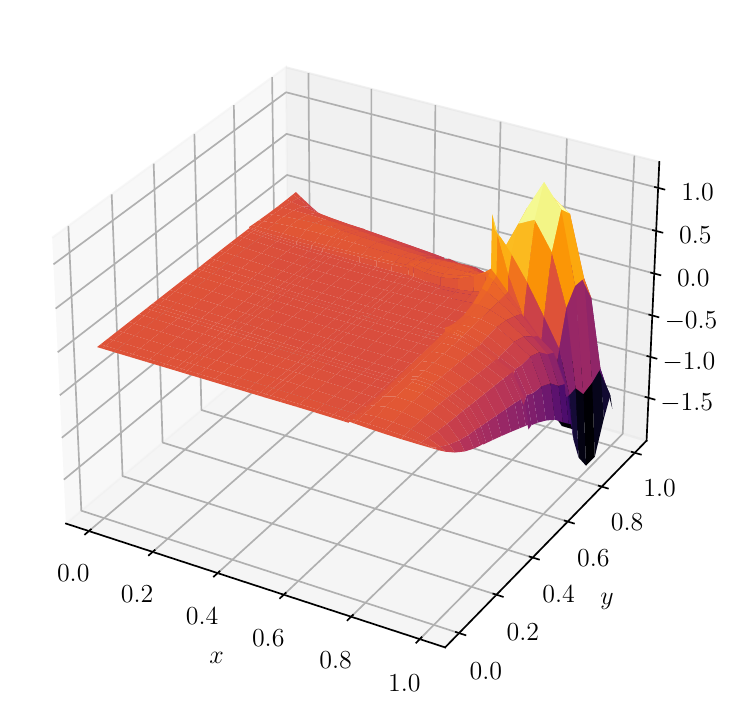}
    \caption{$\mathcal{G}^{\prime}(\mathbb{I} - \sigma_{SG}^{\mathcal{C}})\mathscr{R}\bar{\phi}$ with $k = 4$}
    \label{fig:advDiff2D_u_prime_VMS_p2_k4}
\end{subfigure}
\caption{Exact and computed unresolved scales with varying $k$ on a coarse orthogonal mesh with polynomial degree {$p = 3$} and $4 \times 4$ elements}
\label{fig:advDiff2D_u_prime_VMS_p2}
\end{figure}

We start the initial tests with $\nu = 0.01$ and $\varepsilon = 0.1$ which gives a minimum Peclet number of $\approx52$ and a local maximum value of a $1000$ where we use $\alpha = 100$ for the manufactured solution in \eqref{eq:manufactured_soln}. 
\begin{figure}[htp]
\begin{subfigure}{0.33\linewidth}
    \centering
    \includegraphics[width = \linewidth]{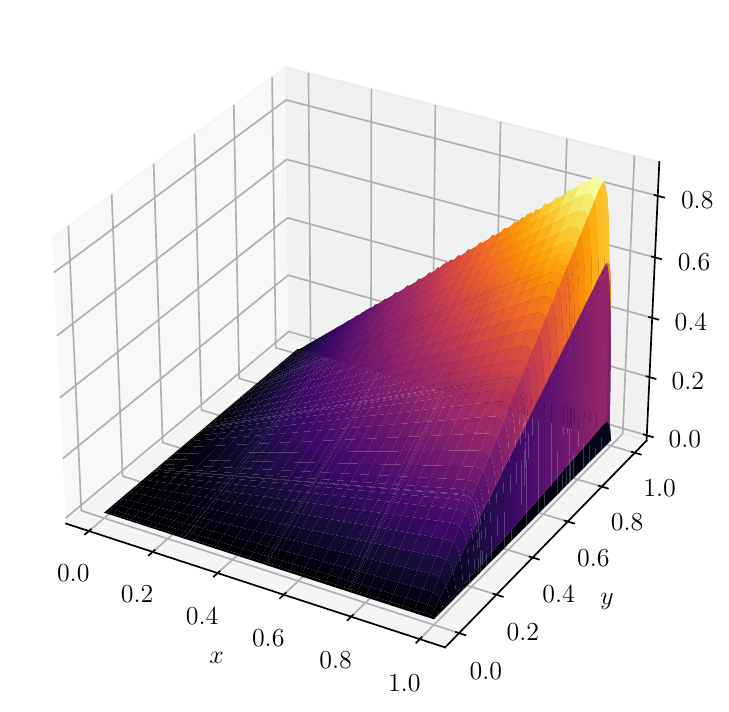}
    \caption{Exact solution}
    \label{fig:advDiff2D_u_ex}
\end{subfigure}
\begin{subfigure}{0.33\linewidth}
    \centering
    \includegraphics[width = \linewidth]{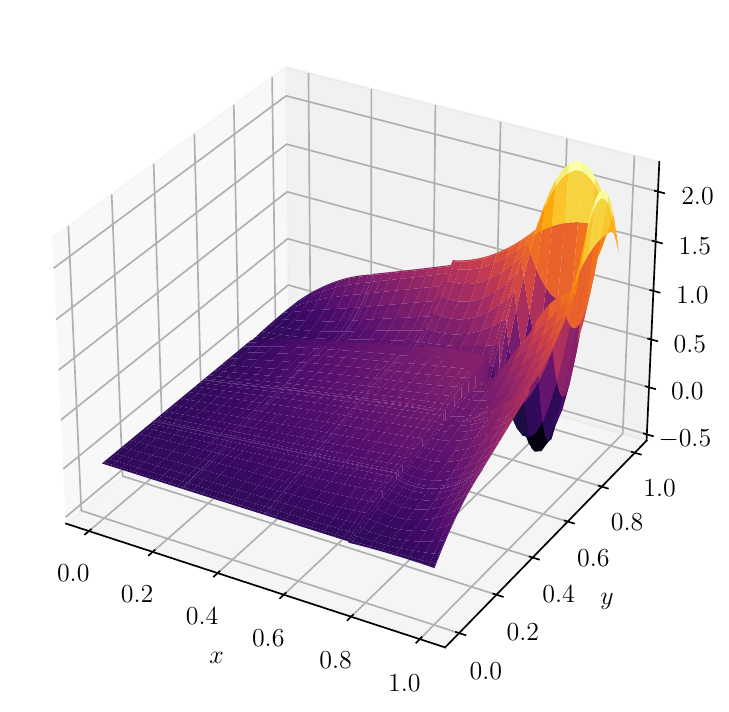}
    \caption{Projection}
    \label{fig:advDiff2D_u_bar_ex_p2}
\end{subfigure}
\begin{subfigure}{0.33\linewidth}
    \centering
    \includegraphics[width = \linewidth]{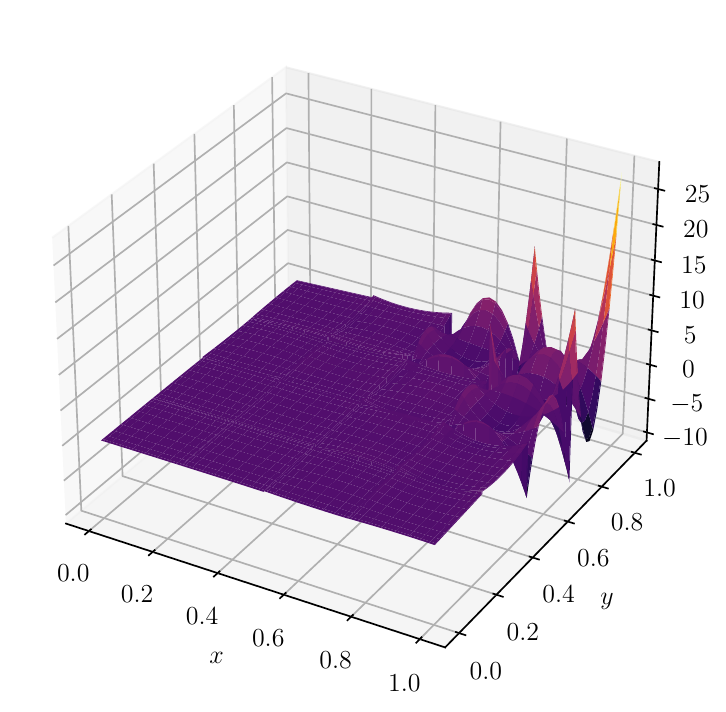}
    \caption{Galerkin solution }
    \label{fig:advDiff2D_u_Glk}
\end{subfigure}
\caption{Exact solution to 2D steady advection-diffusion problem along with its optimal projection onto an orthogonal mesh with polynomial degree $p = 3$ and $4 \times 4$ elements and the Galerkin solution on the same mesh}
\end{figure}
\Cref{fig:advDiff2D_u_ex} and \Cref{fig:advDiff2D_u_bar_ex_p2} show the plots of the exact (manufactured) solution and its projection onto a coarse orthogonal mesh, respectively. Subsequently, \Cref{fig:advDiff2D_u_Glk} shows the Galerkin solution on the same coarse mesh. We once again note that the Galerkin scheme is stable, however, unlike the optimal projection, yields a highly oscillatory solution. Moving on to \Cref{fig:advDiff2D_u_bar_VMS_p2}, we show the optimal projection and the VMS solutions computed on a fixed coarse orthogonal mesh with varying refinement parameter $k$. We can observe the same behaviour as highlighted for the 1D case, where we see that the VMS solution rapidly gets closer to the optimal solution as the fine mesh refinement parameter $k$ is increased. The same can be said for the plots of the unresolved scales plotted in \Cref{fig:advDiff2D_u_prime_VMS_p2}. 

\begin{figure}[htp]
\begin{subfigure}{0.33\linewidth}
    \centering
    \includegraphics[width = \linewidth]{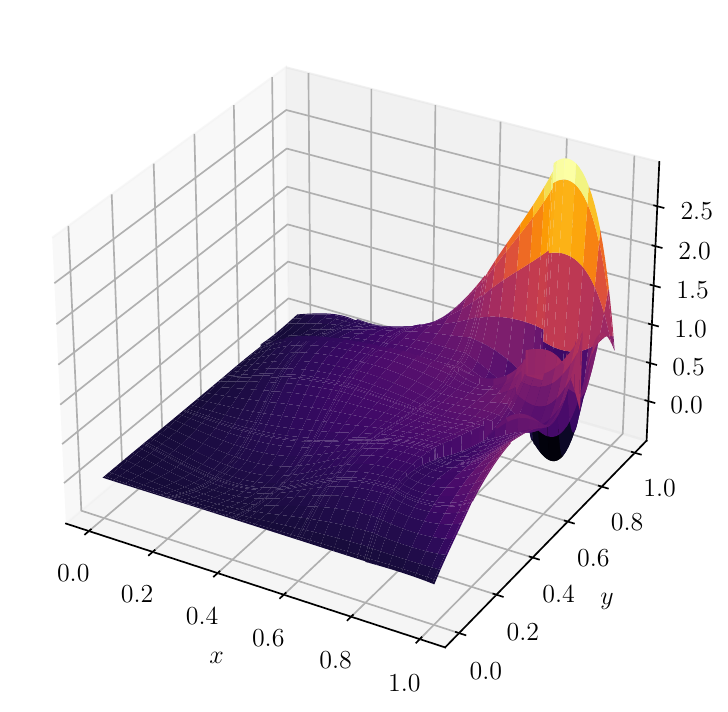}
    \caption{Projection}
    \label{fig:advDiff2D_u_bar_ex_p3_skew}
\end{subfigure}
\begin{subfigure}{0.33\linewidth}
    \centering
    \includegraphics[width = \linewidth]{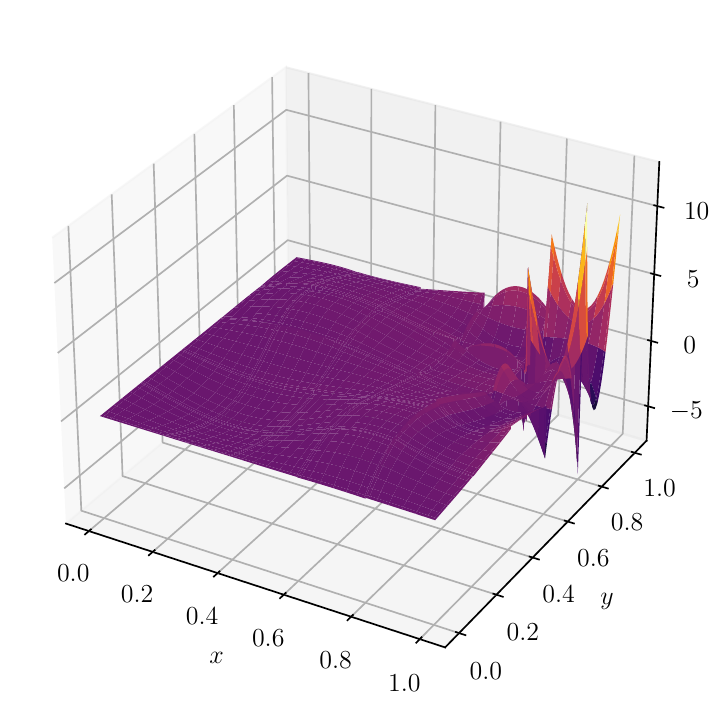}
    \caption{Galerkin solution}
    \label{fig:advDiff2D_u_Glk_skew}
\end{subfigure}
\begin{subfigure}{0.33\linewidth}
    \centering
    \includegraphics[width = \linewidth]{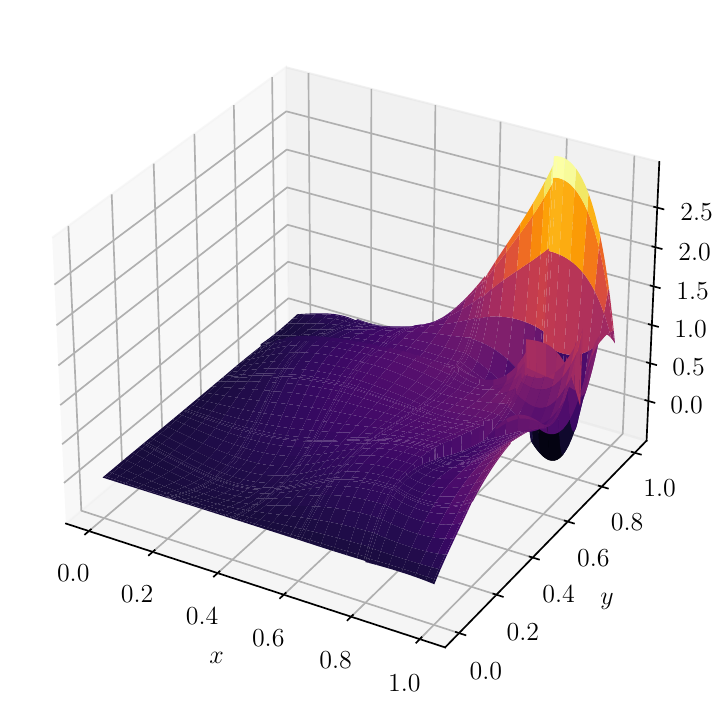}
    \caption{VMS solution with $k = 4$}
    \label{fig:advDiff2D_u_bar_VMS_p2_k4_skew}
\end{subfigure}
\caption{Optimal projection, Galerkin solution, and VMS solution on a coarse curvilinear mesh with polynomial degree $p = 3$ and $5 \times 5$ elements}
\end{figure}

\begin{figure}[htp]
\centering
\begin{subfigure}{0.33\linewidth}
    \centering
    \includegraphics[width = \linewidth]{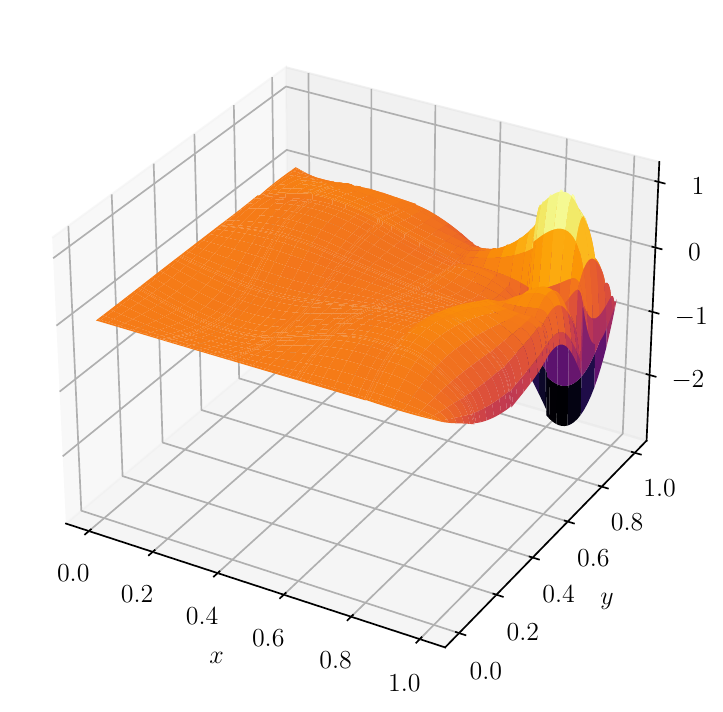}
    \caption{Exact unresolved scales}
\end{subfigure}
\begin{subfigure}{0.33\linewidth}
    \centering
    \includegraphics[width = \linewidth]{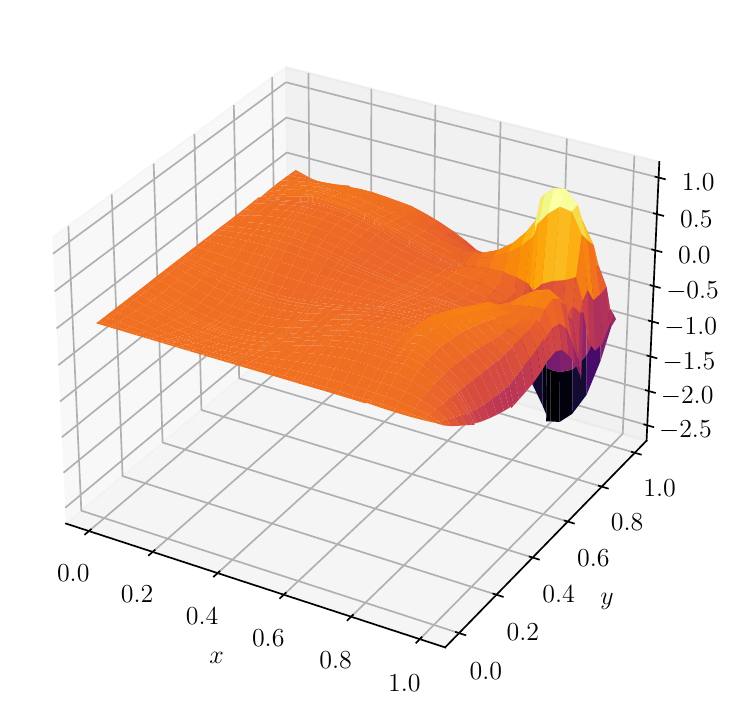}
    \caption{$\mathcal{G}^{\prime}(\mathbb{I} - \sigma_{SG}^{\mathcal{C}})\mathscr{R}\bar{\phi}$ with $k = 4$}
    \label{fig:advDiff2D_u_prime_VMS_p3_k4_skew}
\end{subfigure}
\caption{Exact and computed unresolved scales with varying $k$ on a coarse curvilinear mesh with polynomial degree {$p = 3$} and $5 \times 5$ elements}
\label{fig:advDiff2D_u_prime_VMS_p3_skew}
\end{figure}

Repeating the same test on a curvilinear mesh leads to the same observations as for the orthogonal mesh. Firstly, in \Cref{fig:advDiff2D_u_bar_ex_p3_skew} we have the optimal projection of the exact solution onto the curvilinear grid. Given the curved grid lines, the optimal projection on the coarse curvilinear mesh tends to show weaker similarity to the exact solution as compared to the projection on the orthogonal mesh with similar refinement. Nevertheless, the solution shown in \Cref{fig:advDiff2D_u_bar_ex_p3_skew} is, by definition, the best representation of the exact solution on that mesh in the energy norm. Considering the Galerkin solution shown in \Cref{fig:advDiff2D_u_Glk_skew}, we observe that is it equally poor as the Galerkin solution on the orthogonal mesh. Following that, we turn to \Cref{fig:advDiff2D_u_bar_VMS_p2_k4_skew} and \Cref{fig:advDiff2D_u_prime_VMS_p3_skew} where the VMS solution and the computed unresolved scales for $k = 4$ are shown. These figures clearly show that the VMS solution and the computed unresolved scales closely coincides with their exact counterparts.

Next, we consider the convergence tests where we compute the error using a general form of the norms from \eqref{eq:err_norm_mix_1}-\eqref{eq:err_norm_mix_3}, namely 
\begin{align}
    error \: wrt \: the  \: exact \: solution : \: &\lVert \nabla \cdot \bar{\underline{q}} - (\underline{c} \kappa^{-1} \underline{q} - f) \rVert_{L^2} \label{eq:err_gen_norm_mix_1} \\
    error \: wrt \: the  \: exact \: projection : \: &\left(\lVert \bar{\phi} - \mathcal{P}\phi \rVert_{L^2}^2 + \lVert {\underline{\bar{q}}} - \mathcal{P} \underline{q} \rVert_{L^2}^2\right)^{\frac{1}{2}} \label{eq:err_gen_norm_mix_2} \\
    error \: wrt \: the  \: exact \: fine \: scales : \: &\left(\lVert \phi^{\prime}_k - \phi^{\prime} \rVert_{L^2}^2 + \lVert \underline{c} \kappa^{-1} (\underline{q}^{\prime}_k - \underline{q}^{\prime}) \rVert_{L^2}^2\right)^{\frac{1}{2}}, \label{eq:err_gen_norm_mix_3}
\end{align}
where we again use over-integration with a degree of precision of 25 to evaluate the integrals. We perform the initial convergence tests at a low Peclet number by setting $\nu = 0.08$ yielding a Peclet number of $\alpha = \frac{1}{\nu}$ to be used in \eqref{eq:manufactured_soln}. Solutions for large Peclet numbers require a lot of mesh resolution to capture, thereby making it quite expensive to arrive at the asymptotic region where the error decreases monotonically. Setting a low Peclet number ensures that we can cheaply reach the asymptotic region to assess the rate of convergence. In the latter part of this section, we will also assess the errors for larger Peclet numbers in the non-asymptotic region where the effect of employing VMS is most pronounced. 

\begin{figure}[htp]
\begin{subfigure}{0.49\linewidth}
    \centering
    \includegraphics[width = \linewidth]{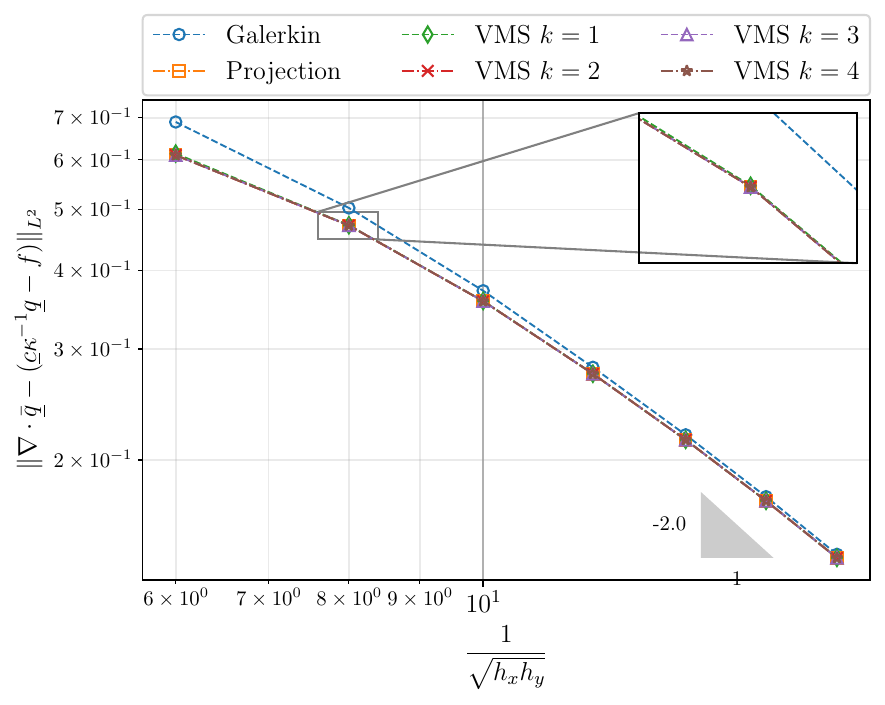}
    \caption{$p = 2$ orthogonal mesh}
    \label{fig:advDiff2D_convergence_surf_p1}
\end{subfigure}
\begin{subfigure}{0.49\linewidth}
    \centering
    \includegraphics[width = \linewidth]{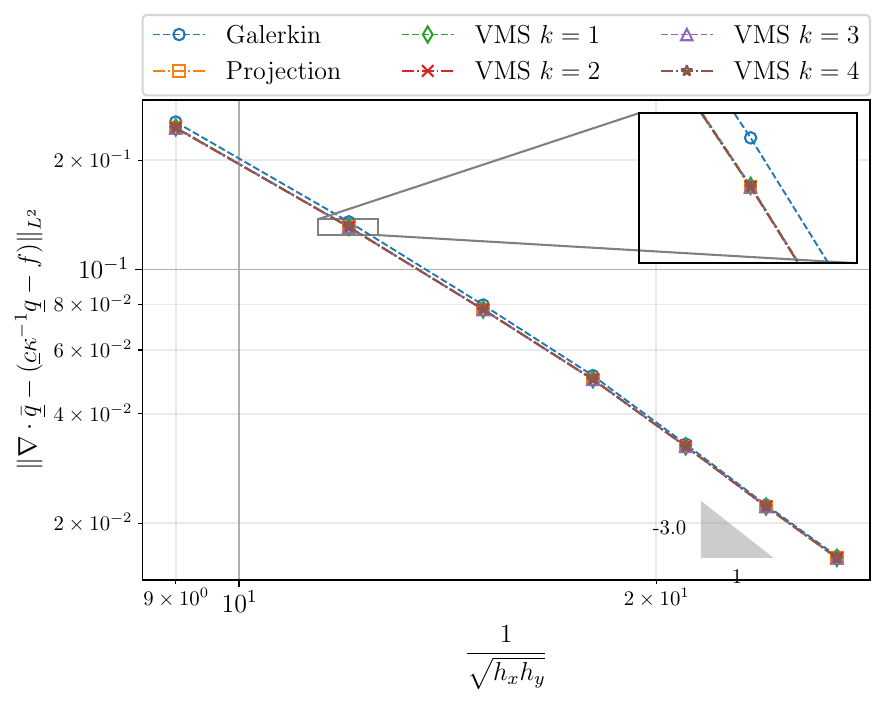}
    \caption{$p = 3$ orthogonal mesh}
    \label{fig:advDiff2D_convergence_surf_p2}
\end{subfigure}
\caption{$h$-$p$ convergence for the mixed formulation of the 2D steady advection-diffusion equation on orthogonal meshes}
\label{fig:advDiff2D_convergence_surf_ortho}
\end{figure}

\begin{figure}[htp]
\begin{subfigure}{0.49\linewidth}
    \centering
    \includegraphics[width = \linewidth]{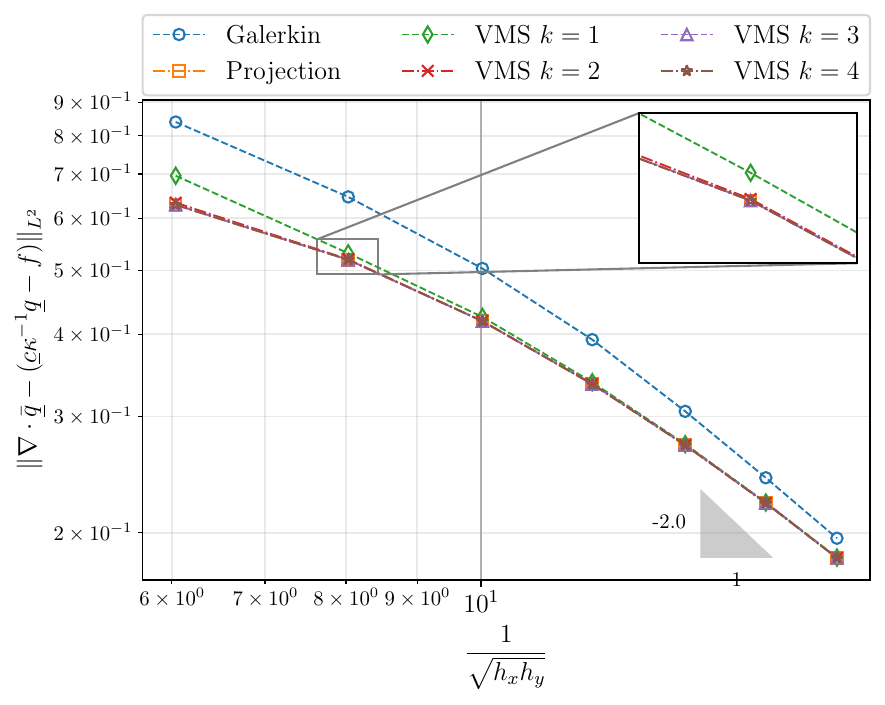}
    \caption{$p = 2$ curvilinear mesh}
    \label{fig:advDiff2D_convergence_surf_p3}
\end{subfigure} 
\begin{subfigure}{0.49\linewidth}
    \centering
    \includegraphics[width = \linewidth]{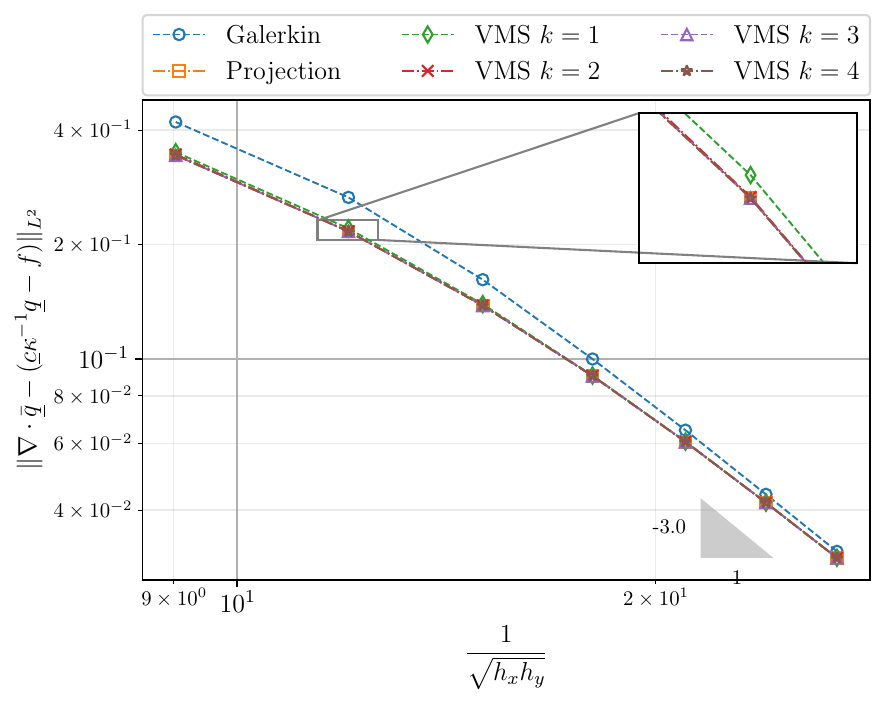}
    \caption{$p = 3$ curvilinear mesh}
    \label{fig:advDiff2D_convergence_surf_p4}
\end{subfigure}
\caption{$h$-$p$ convergence for the mixed formulation of the 2D steady advection-diffusion equation on curvilinear meshes}
\label{fig:advDiff2D_convergence_surf_skew}
\end{figure}

\begin{figure}[htp]
\begin{subfigure}{0.49\linewidth}
    \centering
    \includegraphics[width = \linewidth]{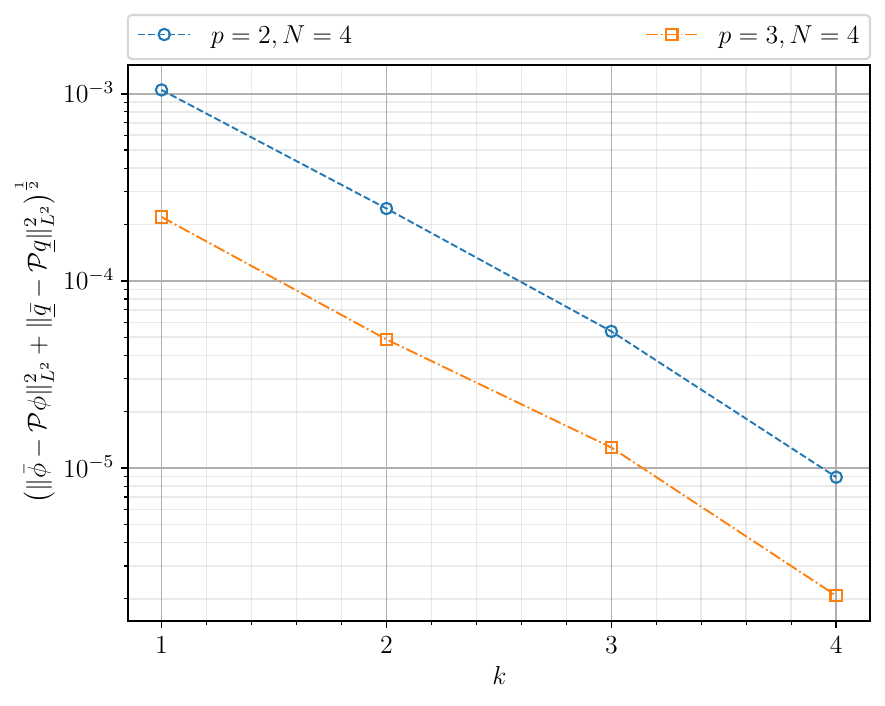}
    \caption{$k$-convergence of the VMS solution}
    \label{fig:advDiff2D_convergence_surf_u_bar_k}
\end{subfigure}
\begin{subfigure}{0.49\linewidth}
    \centering
    \includegraphics[width = \linewidth]{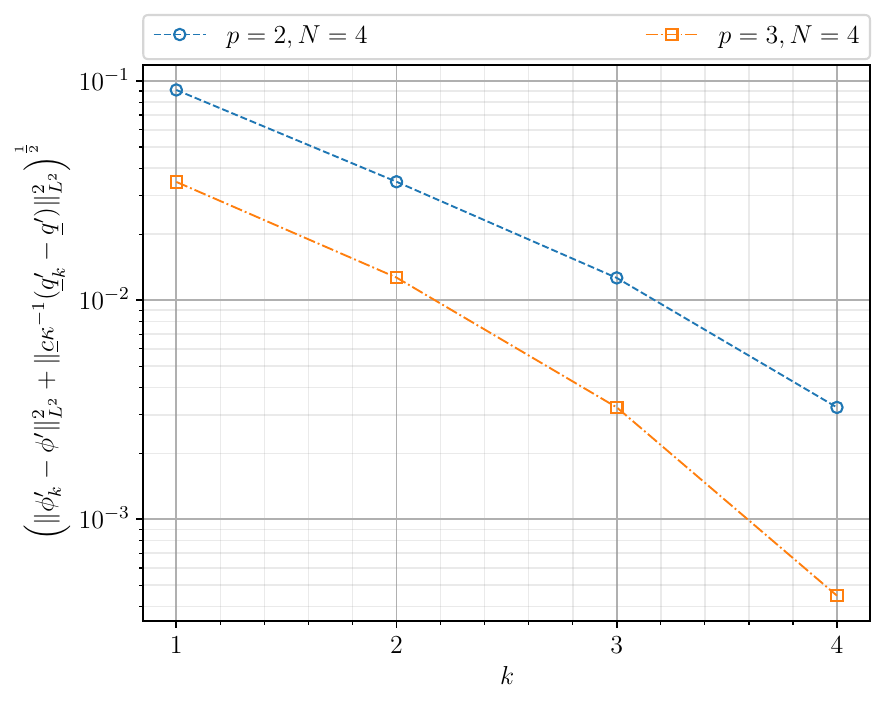}
    \caption{$k$-convergence of the unresolved scales}
    \label{fig:advDiff2D_convergence_surf_u_prime_k}
\end{subfigure}
\caption{$k$-convergence VMS solution and unresolved scales for the mixed formulation of the 2D steady advection-diffusion equation on coarse orthogonal meshes with varying polynomial degree}
\label{fig:advDiff2D_convergence_surf_u_bar_k_u_prime_k}
\end{figure}

\begin{figure}[htp]
\begin{subfigure}{0.49\linewidth}
    \centering
    \includegraphics[width = \linewidth]{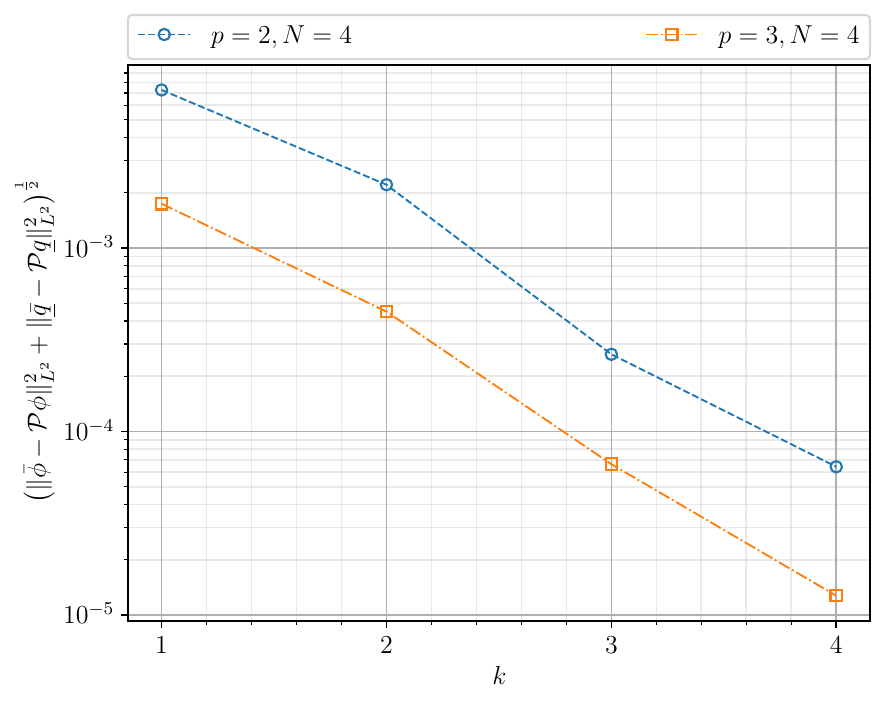}
    \caption{$k$-convergence of the VMS solution}
    \label{fig:advDiff2D_convergence_surf_u_bar_k_skew}
\end{subfigure}
\begin{subfigure}{0.49\linewidth}
    \centering
    \includegraphics[width = \linewidth]{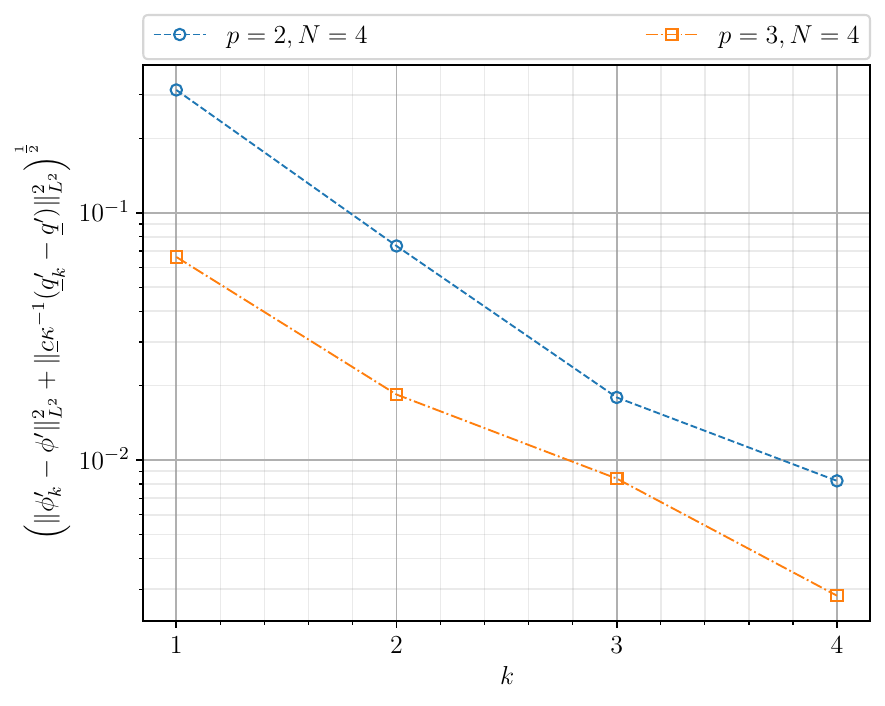}
    \caption{$k$-convergence of the unresolved scales}
    \label{fig:advDiff2D_convergence_surf_u_prime_k_skew}
\end{subfigure}
\caption{$k$-convergence VMS solution and unresolved scales for the mixed formulation of the 2D steady advection-diffusion equation on coarse curvilinear meshes with varying polynomial degree}
\label{fig:advDiff2D_convergence_surf_u_bar_k_u_prime_k_skew}
\end{figure}

Considering the plots in \Cref{fig:advDiff2D_convergence_surf_ortho} and \Cref{fig:advDiff2D_convergence_surf_skew}, we note the same behaviour as for the 1D case. The Galerkin solution has the highest error, while the projection has the lowest by definition and the VMS solutions rapidly approach the optimal projection. Moreover, we see that all schemes converge at the expected rate of $p$. Once again, we observe exponential convergence with $k$ for the convergence of the VMS solution to the exact projection and the unresolved scales as highlighted through \Cref{fig:advDiff2D_convergence_surf_u_bar_k_u_prime_k} and \Cref{fig:advDiff2D_convergence_surf_u_bar_k_u_prime_k_skew}. These observations are consistent for both the orthogonal and curvilinear mesh where the effect of the curved mesh tends to shift the error curves vertically without affecting the slope. Moreover, we observe a greater improvement with VMS over the Galerkin solution on the curvilinear mesh as highlighted by the larger gap between the error curves of the Galerkin and VMS solutions seen through \Cref{fig:advDiff2D_convergence_surf_skew}.

\begin{figure}[htp]
\begin{subfigure}{0.49\linewidth}
    \centering
    \includegraphics[width = \linewidth]{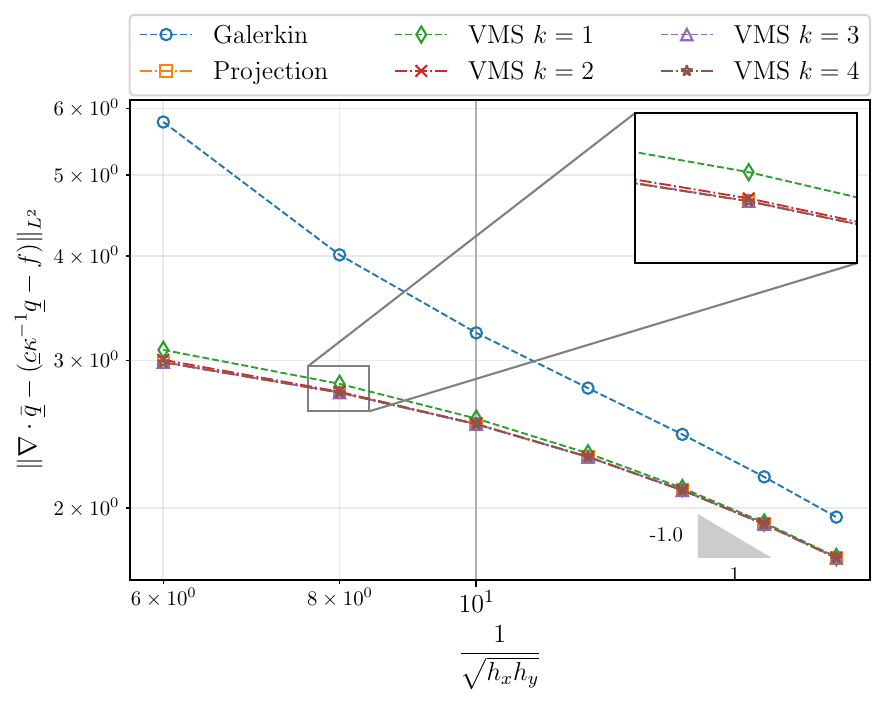}
    \caption{$\alpha = 50$ ($\nu = 0.02, \varepsilon = 1$)}
    \label{fig:convergence_nu_0.02}
\end{subfigure}
\begin{subfigure}{0.49\linewidth}
    \centering
    \includegraphics[width = \linewidth]{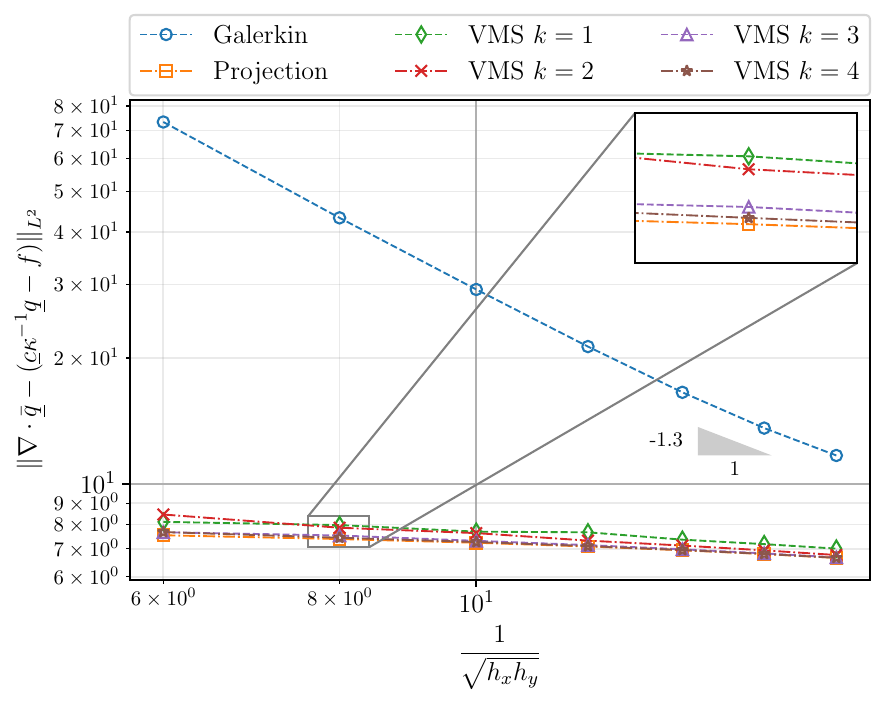}
    \caption{$\alpha = 200$ ($\nu = 0.005, \varepsilon = 1$)}
    \label{fig:convergence_nu_0.005}
\end{subfigure}
\caption{$h$-convergence for the mixed formulation of the 2D steady advection-diffusion equation for different Peclet numbers with polynomial degree $p = 2$}
\label{fig:convergence_nu}
\end{figure}

We now move on to assessing the effect of the Peclet number where we consider two values, namely $\alpha = 50$ and $\alpha = 200$. We achieve this uniform Peclet numbers by setting $\nu = 0.02$ and $\nu = 0.005$ for the respective cases and we set $\varepsilon = 1$ in \eqref{eq:diff_tensor}. The $h$-convergence plots for these two cases for a polynomial degree of $p = 2$ are shown in \Cref{fig:convergence_nu}. The first thing to note in the plots in \Cref{fig:convergence_nu} is that we are not in the asymptotic region hence the rate of convergence is lower than expected. Naturally, we would require much finer meshes to capture the sharp layer of the exact solution at high Peclet numbers. Nevertheless, the effect of employing VMS analysis is the most prominent in the non-asymptotic region and we clearly see the expected improvement in the solution. It is evident that the VMS solution possesses a significantly smaller error than the Galerkin solution and again approaches. The higher Peclet number does influence the rate at which the VMS solution converges to the projection, however, as seen through \Cref{fig:convergence_nu_0.02} and \Cref{fig:convergence_nu_0.005}. 

\begin{figure}[htp]
\begin{subfigure}{0.33\linewidth}
    \centering
    \includegraphics[width = \linewidth]{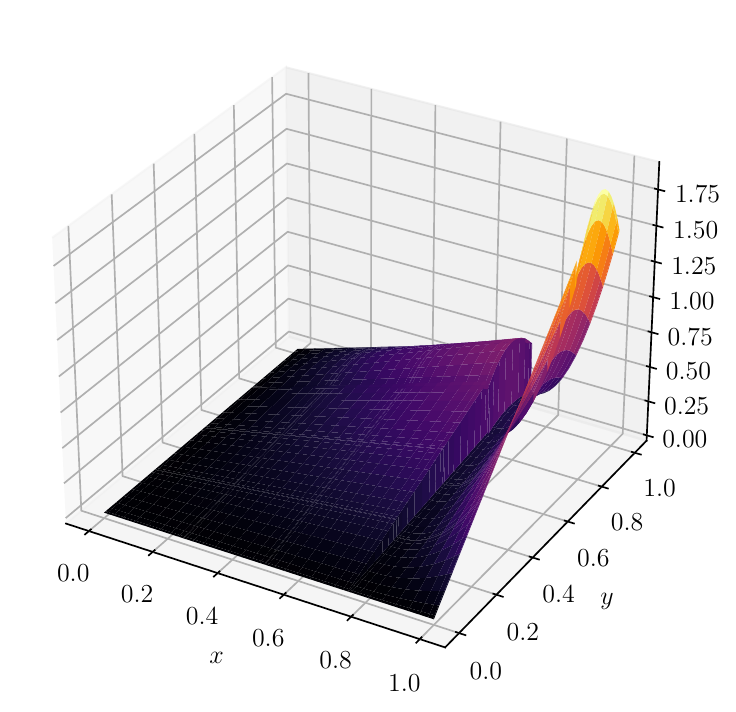}
    \caption{Projection}
    \label{fig:advDiff2D_u_bar_ex_p3_ansio}
\end{subfigure}
\begin{subfigure}{0.33\linewidth}
    \centering
    \includegraphics[width = \linewidth]{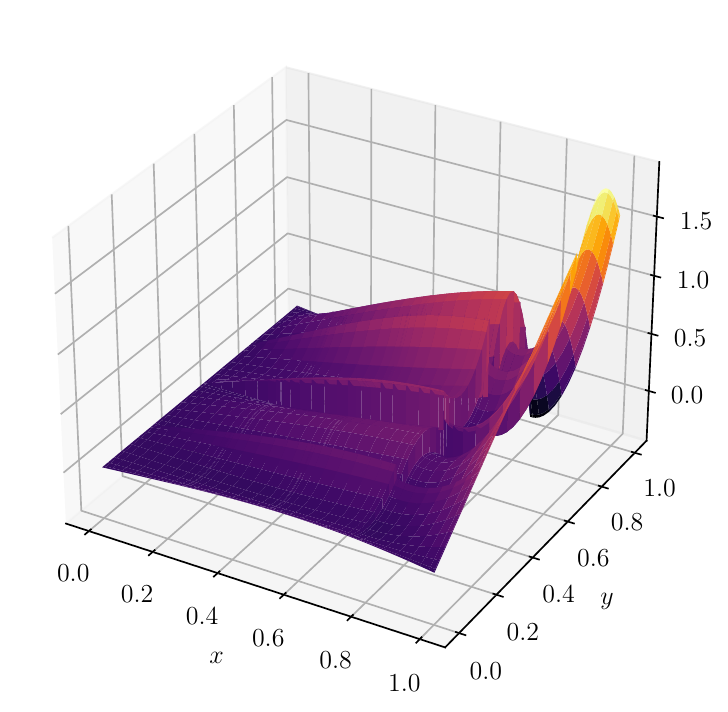}
    \caption{Galerkin solution}
    \label{fig:advDiff2D_u_Glk_ansio}
\end{subfigure}
\begin{subfigure}{0.33\linewidth}
    \centering
    \includegraphics[width = \linewidth]{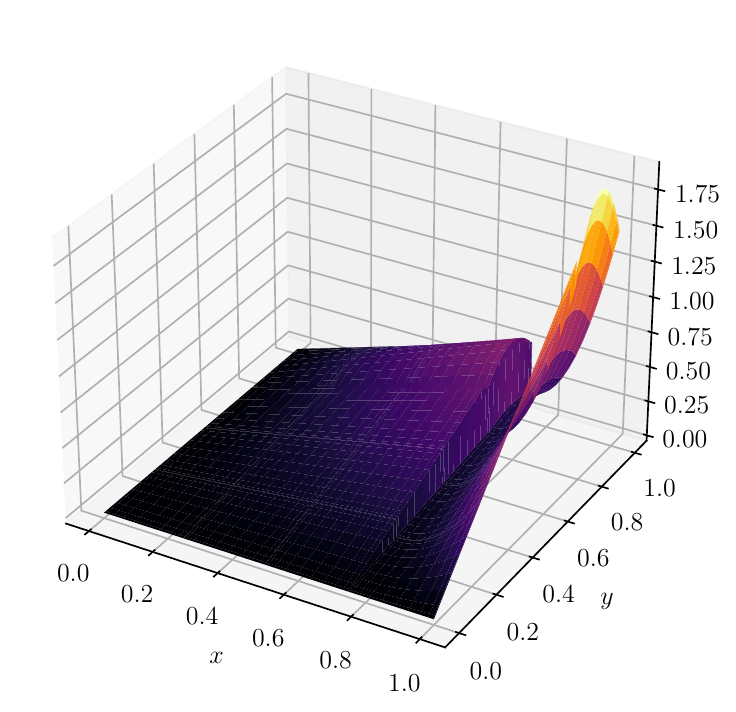}
    \caption{VMS solution with $k = 5$}
    \label{fig:advDiff2D_u_bar_VMS_p2_k4_ansio}
\end{subfigure}
\caption{Optimal projection, Galerkin solution, and VMS solution on a coarse orthogonal mesh with polynomial degree $p = 3$ and
$4 \times 4$ elements}
\end{figure}

\begin{figure}[htp]
\begin{subfigure}{0.33\linewidth}
        \centering
        \includegraphics[width=\linewidth]{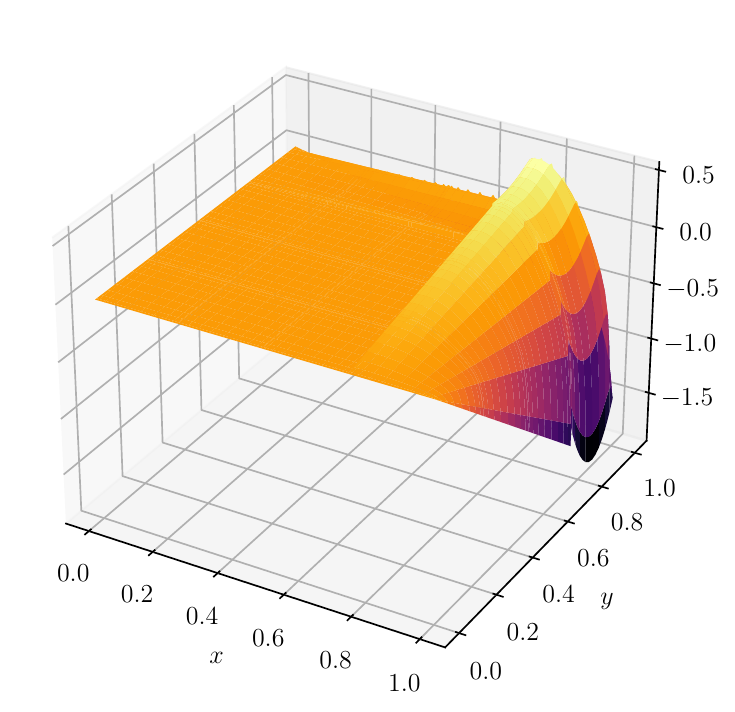}
        \caption{Exact unresolved scales}
        \label{fig:dvDiff2D_u_prime_ex_p2_k4_ansio}
\end{subfigure}
\begin{subfigure}{0.33\linewidth}
        \centering
        \includegraphics[width=\linewidth]{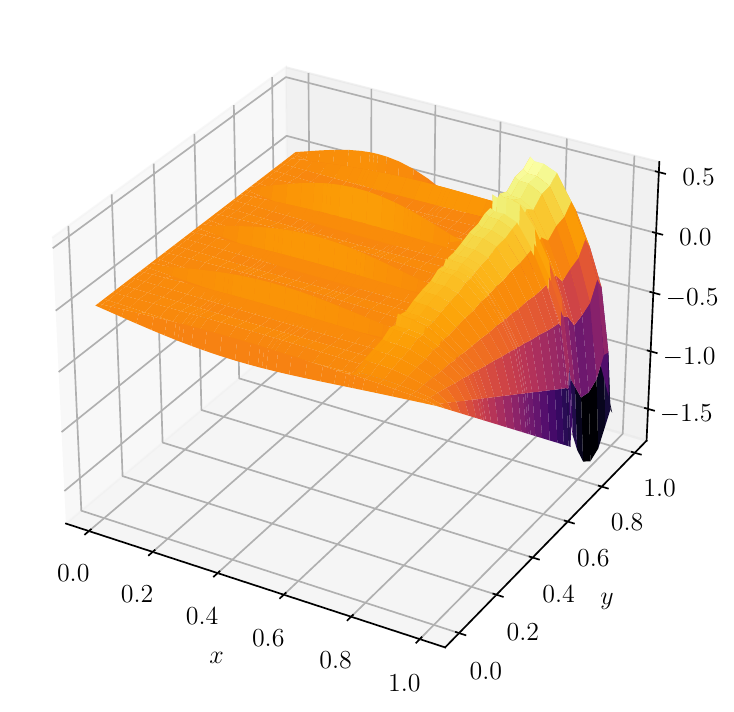}
        \caption{$\mathcal{G}^{\prime}(\mathbb{I} - \sigma_{SG}^{\mathcal{C}})\mathscr{R}\bar{\phi}$ with $k = 5$}
        \label{fig:dvDiff2D_u_prime_VMS_p3_k5_ansio}
\end{subfigure}
\begin{subfigure}{0.33\linewidth}
        \centering
        \includegraphics[width=\linewidth]{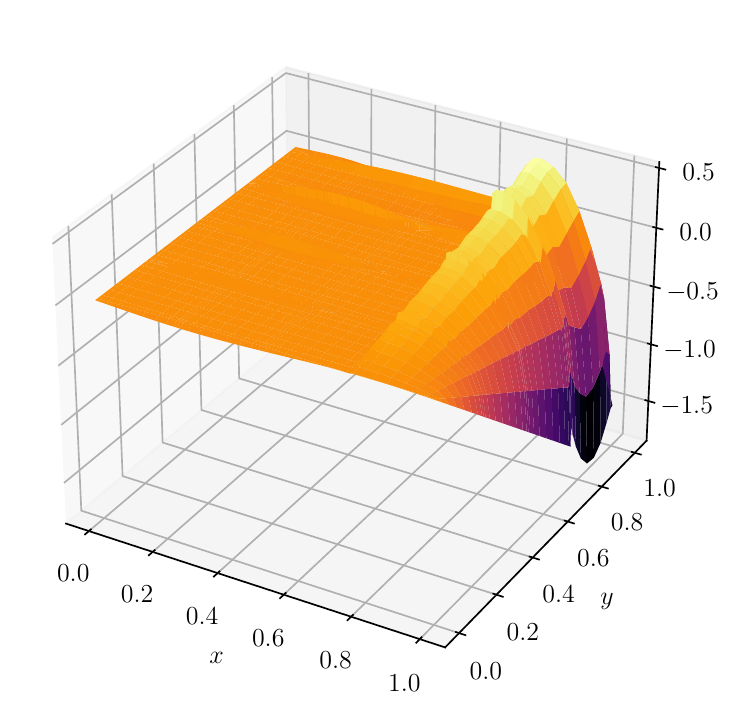}
        \caption{$\mathcal{G}^{\prime}(\mathbb{I} - \sigma_{SG}^{\mathcal{C}})\mathscr{R}\bar{\phi}$ with $k = 6$}
        \label{fig:dvDiff2D_u_prime_VMS_p3_k6_ansio}
\end{subfigure}
\caption{Exact and computed unresolved scales with varying k on a coarse orthogonal mesh with polynomial degree $p = 3$ and $4 \times 4$ elements}
\label{fig:dvDiff2D_u_prime_VMS_ansio}
\end{figure}

To further assess the effect of the Peclet number, we consider an additional test similar to the one presented in \cite{Giraldo2023ANorms}, by considering the following parameters for the advection speed and the diffusion tensor
\begin{gather}
    \underline{c} = \boldsymbol{1} = [1, 1]^T, \quad \quad \kappa = \left[\begin{array}{cc}
        1 & 0 \\
        0 & 10^{-6}
    \end{array} \right].
\end{gather}
This highly anisotropic diffusion tensor makes the problem quite challenging as the problem becomes heavily convection-dominated. We take the same manufactured solution as shown in \eqref{eq:manufactured_soln} with $\alpha = 100$ which generates a sharp layer in the domain (shown in \Cref{fig:advDiff2D_u_ex}). \\ \\
\Cref{fig:advDiff2D_u_bar_ex_p3_ansio} shows the optimal projection of the exact solution onto an orthogonal mesh with polynomial degree 3 and \Cref{fig:advDiff2D_u_bar_VMS_p2_k4_ansio} shows the Galerkin solution on that same mesh. Contrary to the prior tests, we can note that we only observe slight discrepancies between the Galerkin solution and the optimal projection. The difference between the two is most prominent in the interior of the domain where the Galerkin solution exhibits oscillations, whereas the projection remains smooth. Besides that difference, we see that the two solutions behave very similarly near the sharp layer where both solutions exhibit a sharp increase near the boundary. Though the distinctions between the two solutions may not be as significant as the preceding test cases, the projected solution is, by definition, the optimal solution in the energy norm and we see in \Cref{fig:advDiff2D_u_bar_VMS_p2_k4_ansio} that the VMS approach successfully yields a solution that closely resembles the optimal projection. Furthermore, \Cref{fig:dvDiff2D_u_prime_VMS_ansio} shows that the Fine-Scale Greens' function consistently computes the unresolved scales. We do note that the highly anisotropic diffusion tensor does pose difficulties for the computations of the fine scales as apparent through the oscillations observed for the $k = 5$ case in \Cref{fig:advDiff2D_u_bar_VMS_p2_k4_ansio}. Despite that, however, the VMS solution for the same $k$ value shown closely resembles the exact projection seen through \Cref{fig:advDiff2D_u_bar_VMS_p2_k4_ansio}. In line with expectations, we see from \Cref{fig:dvDiff2D_u_prime_VMS_p3_k6_ansio} that increasing $k$ eliminates these oscillations as the approximation of the classic Greens' function is improved.

Lastly, we close off the discussion on the numerical experiments by noting an important fact regarding the orthogonality between the computed unresolved scales and the resolved space. Throughout the derivations presented in \Cref{sec:optmal_skew_sym}, we use orthogonality conditions between the resolved and unresolved subspaces invoked by the projector to simplify the variational form. These orthogonality conditions are numerically satisfied through the construction of the Fine-Scale Greens' function and are thus fully independent of the polynomial degree of the coarse space where we solve the full PDE as well as the fine space where the classic Greens' function is approximated. To highlight this fact, we present \Cref{tab:ortho1}, \Cref{tab:ortho2}, and \Cref{tab:ortho3} showing the numerically computed values of the corresponding inner products between the computed unresolved scales and the resolved space. We see that all the values in the tables are in the order of machine precision, thereby confirming that the orthogonality conditions are exactly satisfied in the numerical setting.

\begin{table}[H]
    \centering
    \caption{Inner product from \eqref{eq:ortho_dir} between the computed unresolved scales and the resolved scales for the advection-diffusion problem in the direct formulation with varying $p$ and $k$}
    \label{tab:ortho1}
    \begin{tabular}{c|c|c|c|c}
         $\bilnear{\nabla v^h}{\nabla \phi^{\prime}}$ & $k = 1$  & $k = 2$ & $k = 3$ & $k = 4$\\ \hline \hline
         $p = 1$ & $7.66\times 10^{-15}$ & $3.82\times 10^{-14}$ & $8.88\times 10^{-16}$& $-1.31\times 10^{-14}$ \\ \hline
         $p = 2$ & $4.08\times 10^{-14}$ & $-2.66\times 10^{-15}$ & $-1.59\times 10^{-14}$ & $-1.2434\times 10^{-14}$ \\\hline
         $p = 4$ & $-6.57\times 10^{-14}$ & $-1.77\times 10^{-14}$ & $-1.24\times 10^{-14}$ & $3.81\times 10^{-14}$ \\ 
    \end{tabular}
\end{table}

\begin{table}[H]
    \centering
    \caption{Inner product from \eqref{eq:mix_ortho_1} between the computed unresolved scales and the resolved scales for the advection-diffusion problem in the mixed formulation with varying $p$ and $k$}
    \label{tab:ortho2}
    \begin{tabular}{c|c|c|c|c}
         $\bilnear{v^h}{\underline{q}^{\prime}} + \bilnear{\nabla \cdot v^h}{\phi^{\prime}}$ & $k = 1$  & $k = 2$ & $k = 3$ & $k = 4$\\ \hline \hline
         $p = 1$ & $2.39\times 10^{-14}$ & $-6.70\times 10^{-15}$ & $6.85\times 10^{-15}$& $-1.62\times 10^{-14}$ \\ \hline
         $p = 2$ & $4.38\times 10^{-15}$ & $-9.21\times 10^{-15}$ & $1.11\times 10^{-16}$ & $-2.77\times 10^{-16}$ \\\hline
         $p = 4$ & $-3.34\times 10^{-15}$ & $6.96\times 10^{-16}$ & $-6.95\times 10^{-16}$ & $-4.28\times 10^{-16}$ \\ 
    \end{tabular}
\end{table}

\begin{table}[H]
    \centering
    \caption{Inner product from \eqref{eq:mix_ortho_2} between the computed unresolved scales and the resolved scales for the advection-diffusion problem in the mixed formulation with varying $p$ and $k$}
    \label{tab:ortho3}
    \begin{tabular}{c|c|c|c|c}
         $\bilnear{\eta^h}{\nabla \cdot \underline{q}^{\prime}}$ & $k = 1$  & $k = 2$ & $k = 3$ & $k = 4$\\ \hline \hline
         $p = 1$ & $7.95\times 10^{-13}$ & $-3.41\times 10^{-13}$ & $-1.13\times 10^{-12}$& $-1.59\times 10^{-12}$ \\ \hline
         $p = 2$ & $-5.11\times 10^{-13}$ & $-8.52\times 10^{-14}$ & $-5.68\times 10^{-13}$ & $2.27\times 10^{-13}$ \\\hline
         $p = 4$ & $7.10\times 10^{-15}$ & $-4.97\times 10^{-14}$ & $2.48\times 10^{-14}$ & $-2.48\times 10^{-14}$ \\ 
    \end{tabular}
\end{table}
\section{Summary}


\label{sec:summary}

In the present study, we formulated an algebraic VMS approach employing the Fine-Scale Greens' function. The proposed VMS approach uses two meshes with varying refinements where the full problem is solved on a coarse mesh while a fine mesh is used to compute the approximate Greens' function, but only for the underlying symmetric problem. We considered the 1D advection-diffusion problem in both direct and mixed formulations along with the 2D advection-diffusion problem in a mixed formulation for the numerical tests. The proposed approach can also be naturally extended to 3D, see \cite[\S 5]{Jain2021ConstructionMeshes}. For the aforementioned cases, the traditional Galerkin method, characterised by its central-like scheme, produces solutions with significant oscillations, particularly on coarse meshes, resulting in the largest errors. On the other hand, the optimal projection exhibits the lowest errors by construction. The newly proposed VMS method shows promise, with its errors lying close to those of the optimal projection. Additionally, the VMS solution demonstrates rapid convergence towards the optimal projection as the Greens' function for the symmetric problem is refined by increasing the polynomial degree increment $k$. The convergence analysis further reveals that all methods converge at the expected rate upon mesh refinement, and the VMS solution and the unresolved scales achieve exponential convergence with $k$ refinement. Notably, the error in unresolved scales is consistently larger than in resolved scales, implying that a highly accurate estimation of fine scales is not crucial for a good approximation of the optimal projection. These findings highlight the potential of the VMS method to provide more precise solutions in numerical simulations of advection-diffusion problems. We emphasise that the VMS method is not a stabilisation technique as the base Galerkin discretisation scheme is already stable. As such, VMS acts to supplement the base scheme to improve the solution by yielding a solution close to the desired projection. The subsequent part of this study will involve applying the presented linear theory to more complex problems including time-dependent problems and the extension to non-linear problems such as Burgers' and Navier-Stokes equations.

\section*{Acknowledgments}
Suyash Shrestha and Esteban Ferrer acknowledge the funding received by Clean Aviation Joint Undertaking under the European Union’s Horizon Europe research and innovation programme under Grant Agreement HERA (Hybrid-Electric Regional Architecture) no. 101102007. Views and opinions expressed are,
however, those of the author(s) only and do not necessarily reflect those of the European Union or CAJU. Neither the European Union nor the granting authority can be held responsible for them. Esteban Ferrer and Gonzalo Rubio acknowledge the funding received by the Grant DeepCFD (Project No. PID2022-137899OB-I00) funded by MICIU/AEI/10.13039/501100011033 and by ERDF, EU. 



\bibliographystyle{elsarticle-num} 
\bibliography{references_mendeley}

\newpage
\appendix
\section{Greens' function plots}
\label{app:A}
In this appendix, we present plots of the exact classic Greens' function for the Poisson problem and its approximation computed using the methodology described in \Cref{subsec:greens_approx}. 

The 1D Poisson problem with homogeneous Dirichlet boundary conditions reads
\begin{gather}
    -\parddxdx{u}{x} = f, \quad \text{in } \Omega = [0, 1] \label{eq:poisson_eq}\\
    u(x) = 0, \quad \text{on } \partial \Omega. \label{eq:poisson_bc}
\end{gather}
The classic Greens' function for this problem defined in $\Omega$ is given by
\begin{equation}
    g(x, s) = \begin{cases}
        (1 - s)x, \quad x \leq s \\
        s (1 - x), \quad x > s
    \end{cases}.
\end{equation}
Similarly, the 2D problem reads
\begin{gather}
    -\nabla \cdot \nabla \phi = f, \quad \text{in } \Omega = [0, 1]^2 \\
    \phi(\boldsymbol{x}) = 0, \quad \text{on } \partial \Omega,
\end{gather}
and the Greens' function associated with this problem can be expressed as an eigenfunction expansion as follows \cite[\S 8]{RichardHaberman2019HabermanBoundary}
\begin{gather}
    g(\boldsymbol{x}, \boldsymbol{s}) = g(x, y, s_1, s_2) =
    \displaystyle \sum_{n = 1}^{\infty} -\dfrac{(2\sin{\left(n \pi s_1 \right)} \sin{(n \pi x)})}{(n \pi \sinh{(n \pi)})} 
    \begin{cases}
        \sinh{(n \pi (s_2 - 1))} \sinh{(n \pi y)}, \quad & y < s_2 \\
        \sinh{(n \pi (y - 1))} \sinh{(n \pi s_2)}, \quad & y \geq s_2.
    \end{cases}
    \label{eq:2dPoisson_Greens}
\end{gather}

\subsection{1D direct formulaion}
\begin{figure}[H]
    \centering
    \includegraphics[width=0.49\linewidth]{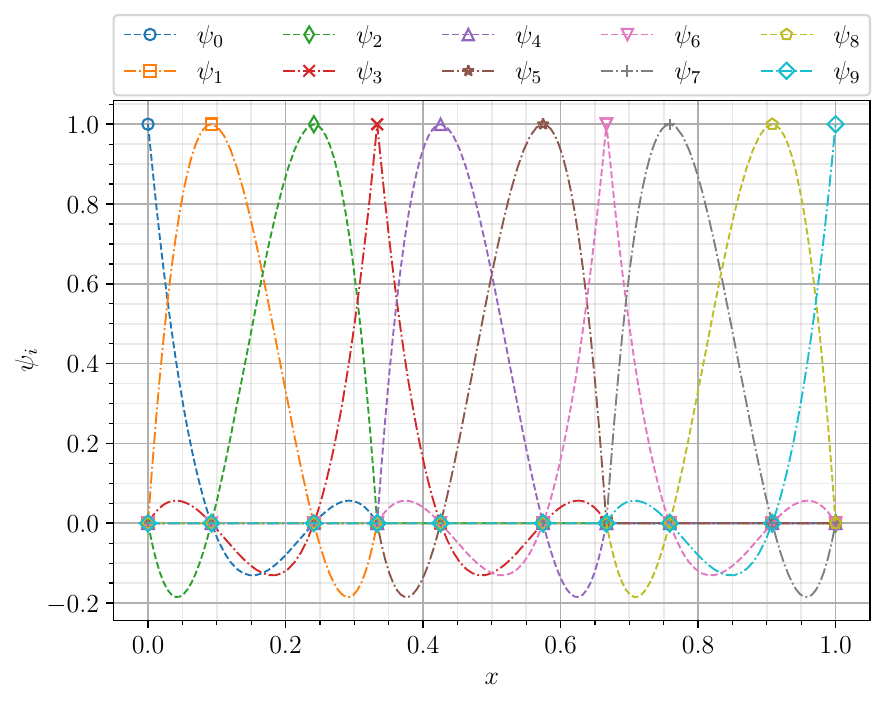}
    \caption{Polynomials of degree 3 over 3 elements spanning the $H^1(\Omega)$ subspace}
    \label{fig:lagrange_basis}
\end{figure}

\begin{figure}[H]
\begin{subfigure}{0.49\linewidth}
    \centering
    \includegraphics[width=\linewidth]{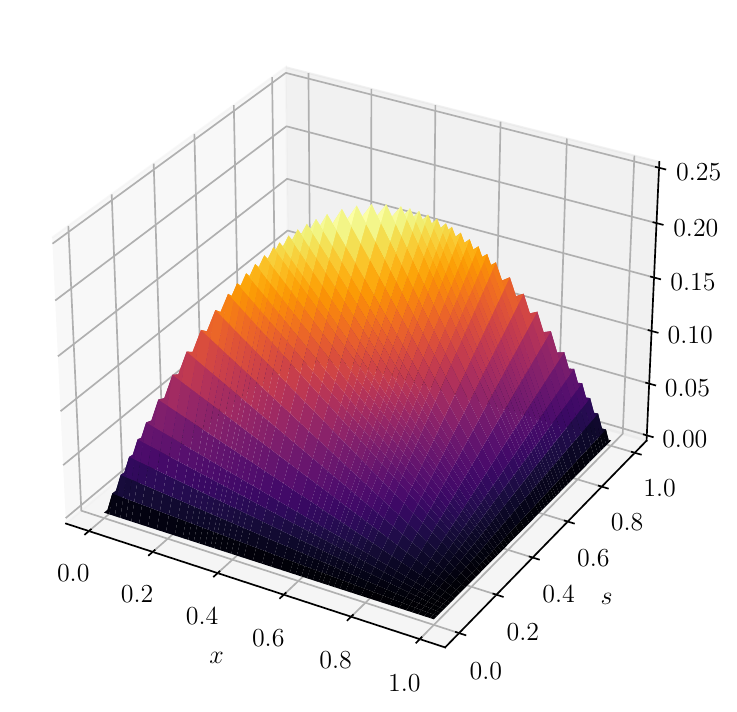}
    \caption{Exact}
\end{subfigure}
\begin{subfigure}{0.49\linewidth}
    \centering
    \includegraphics[width=\linewidth]{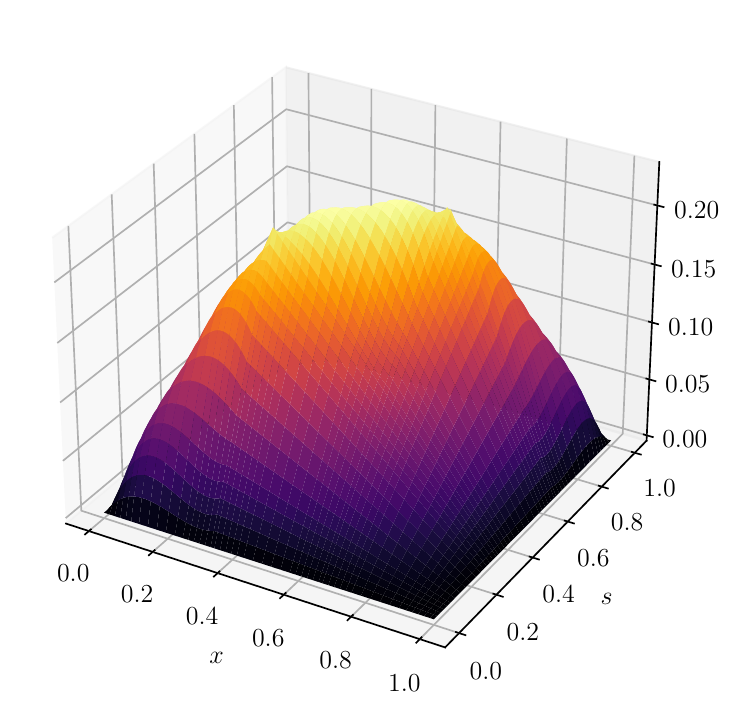}
    \caption{Approximate}
\end{subfigure}
\caption{Exact and approximate Greens' function for 1D Poisson equation in the direct formulation with the approximation computed using polynomials of degree 3 over 3 elements}
\label{fig:greens_approx_direct}
\end{figure}

\begin{figure}[H]
\begin{subfigure}{0.49\linewidth}
    \centering
    \includegraphics[width=\linewidth]{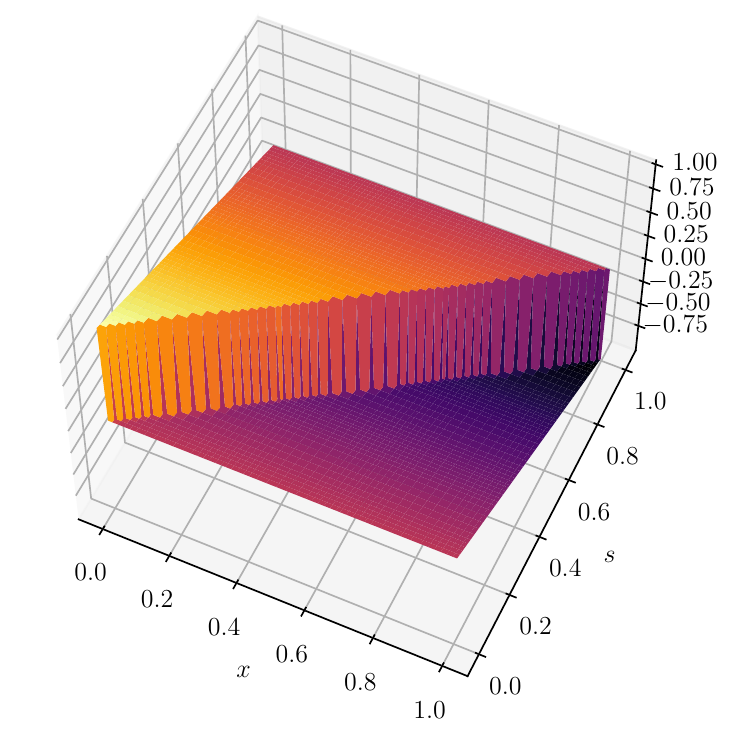}
    \caption{Exact}
\end{subfigure}
\begin{subfigure}{0.49\linewidth}
    \centering
    \includegraphics[width=\linewidth]{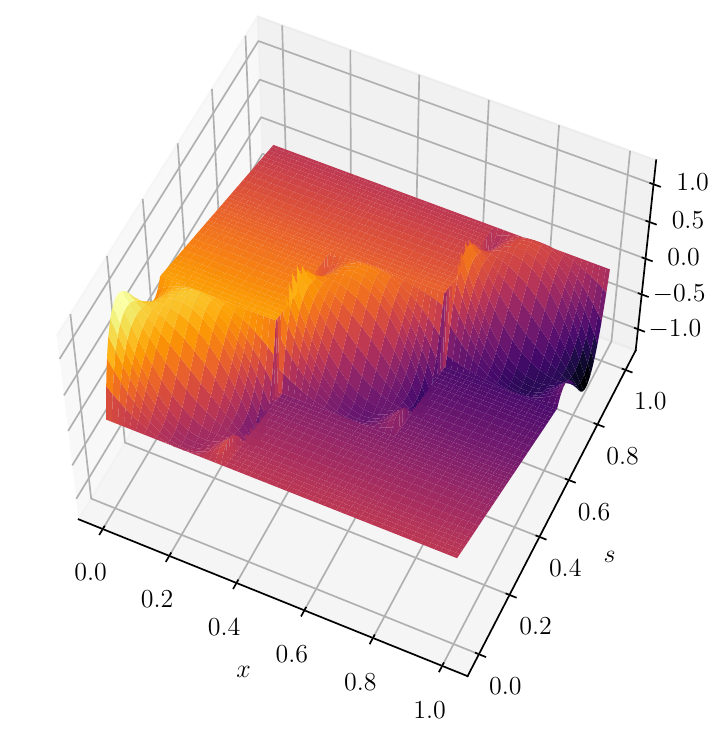}
    \caption{Approximate}
\end{subfigure}
\caption{Derivative of the exact and approximate Greens' function for the 1D Poisson equation in the direct formulation with the approximation computed using polynomials of degree 3 over 3 elements}
\label{fig:dgreens_approx_direct}
\end{figure}

For computing the approximate Greens' function for the direct formulation, we employ \eqref{eq:greens_approx} where we construct the discrete Laplacian operator using Lagrange polynomials as the basis for our discrete $H^1(\Omega)$ subspace. \Cref{fig:lagrange_basis} shows an example of these basis functions along with \Cref{fig:greens_approx_direct} and \Cref{fig:dgreens_approx_direct} which show the corresponding computed Greens' function and its derivative alongside the exact one. 


\subsection{1D mixed formulaion}
\begin{figure}[H]
\begin{minipage}{0.49\linewidth}
    \centering
    \includegraphics[width=\linewidth]{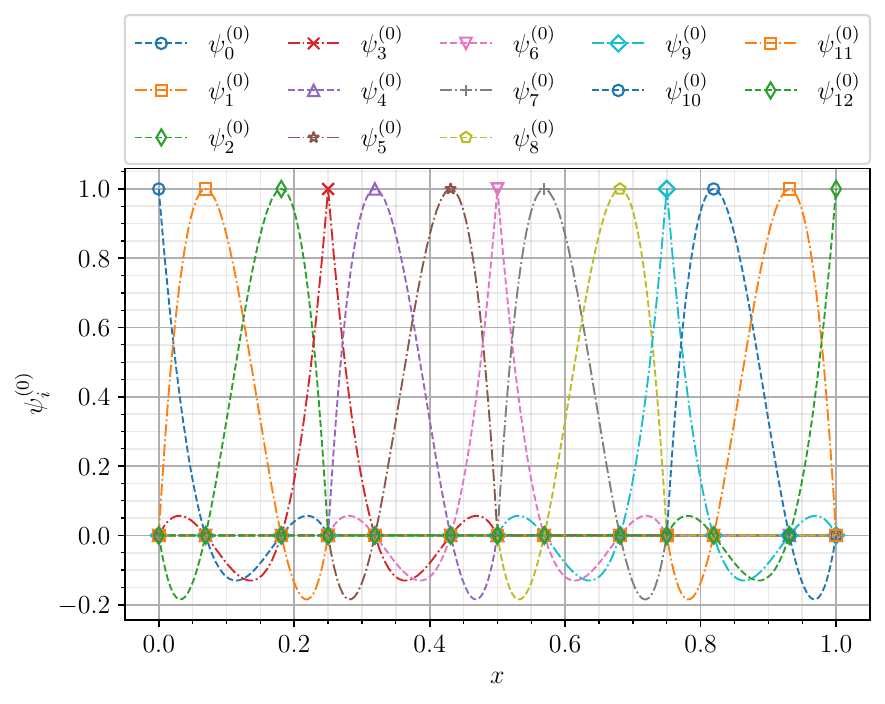}
    \caption{Polynomials of degree 3 over 4 elements spanning the $H^1(\Omega)$ subspace}
    \label{fig:basis_H1}
\end{minipage}
\begin{minipage}{0.49\linewidth}
    \centering
    \includegraphics[width=\linewidth]{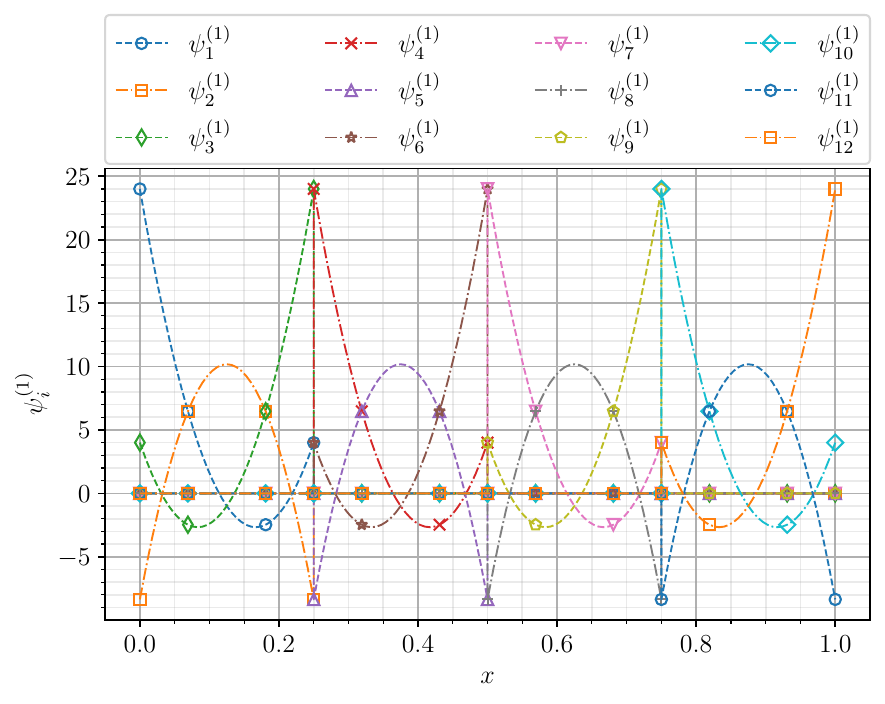}
    \caption{Polynomials of degree 2 over 4 elements spanning the $L^2(\Omega)$ subspace}
    \label{fig:basis_L2}
\end{minipage}
\end{figure}

\begin{figure}[H]
\begin{subfigure}{0.49\linewidth}
    \centering
    \includegraphics[width=\linewidth]{Images/advDiff1D/greens_exact.pdf}
    \caption{Exact}
\end{subfigure}
\begin{subfigure}{0.49\linewidth}
    \centering
    \includegraphics[width=\linewidth]{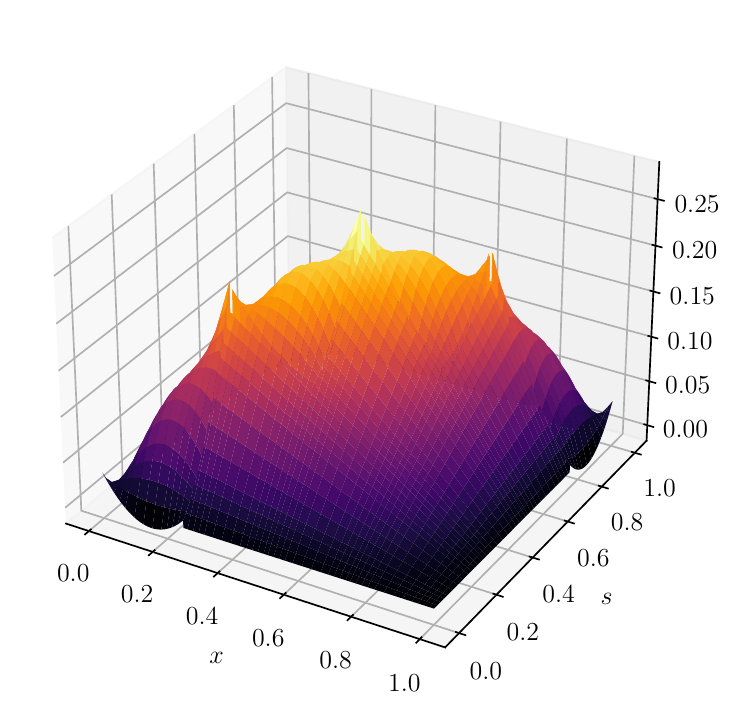}
    \caption{Approximate}
\end{subfigure}
\caption{Exact and approximate Greens' function for the 1D Poisson equation in the mixed formulation with the approximation computed using polynomials of degree 2 over 4 elements}
\label{fig:greens_approx_mix}
\end{figure}

\begin{figure}[H]
\begin{subfigure}{0.49\linewidth}
    \centering
    \includegraphics[width=\linewidth]{Images/advDiff1D/dgreens_exact.pdf}
    \caption{Exact}
\end{subfigure}
\begin{subfigure}{0.49\linewidth}
    \centering
    \includegraphics[width=\linewidth]{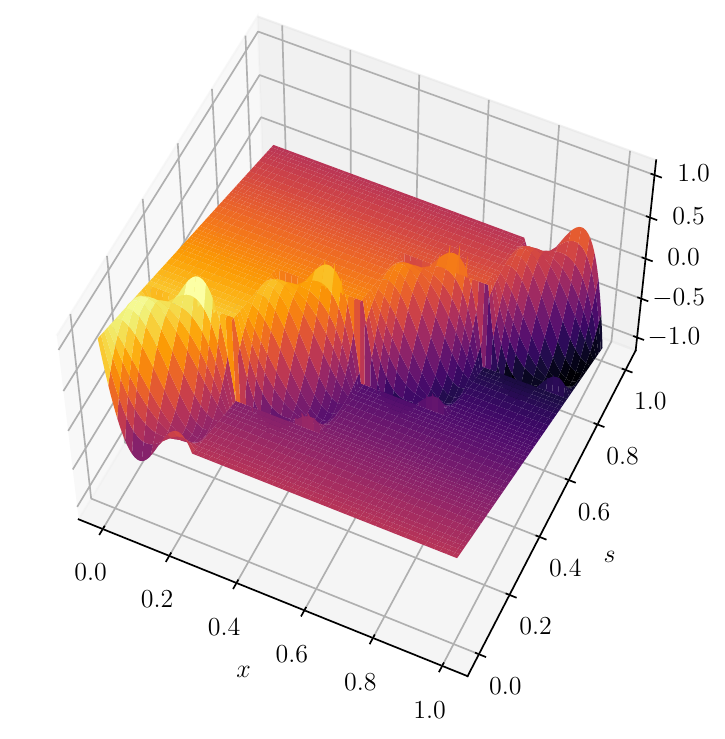}
    \caption{Approximate}
\end{subfigure}
\caption{Derivative of the exact and approximate Greens' function for the 1D Poisson equation in the mixed formulation with the approximation computed using polynomials of degree 2 over 4 elements}
\label{fig:dgreens_approx_mix}
\end{figure}

The same concept is used for the mixed formulation, with the addition of edge polynomials \cite{MarcGerritsma2010EDGEMETHODS} to construct the discrete $L^2(\Omega)$ subspace. Plots of the polynomials spanning the discrete $H^1(\Omega)$ and $L^2(\Omega)$ subspaces are shown in \Cref{fig:basis_H1} and \Cref{fig:basis_L2}, with the corresponding plots of the Greens' function and its derivative show in \Cref{fig:greens_approx_mix} and \Cref{fig:dgreens_approx_mix} respectively. 

\subsection{2D mixed formulaion}
\begin{figure}[H]
\begin{subfigure}{0.33\linewidth}
    \centering
    \includegraphics[width = \linewidth]{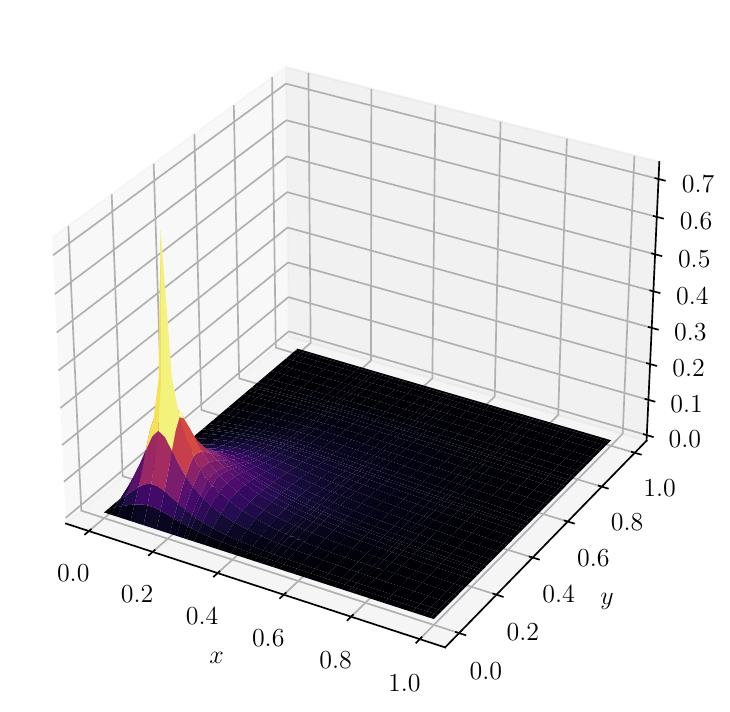}
    \caption{Exact $(s_1, s_2) = (0.125, 0.125)$}
\end{subfigure}
\begin{subfigure}{0.33\linewidth}
    \centering
    \includegraphics[width = \linewidth]{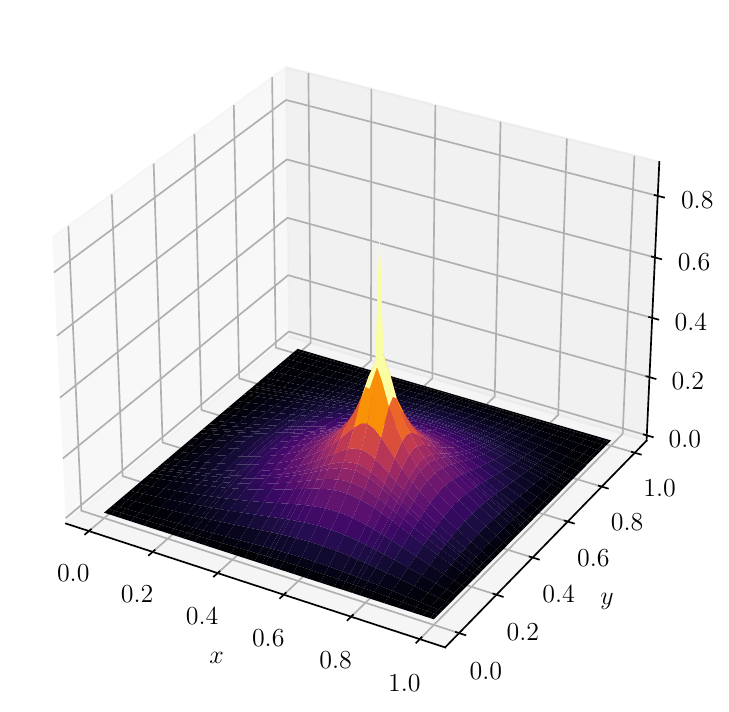}
    \caption{Exact $(s_1, s_2) = (0.625, 0.375)$}
\end{subfigure}
\begin{subfigure}{0.33\linewidth}
    \centering
    \includegraphics[width = \linewidth]{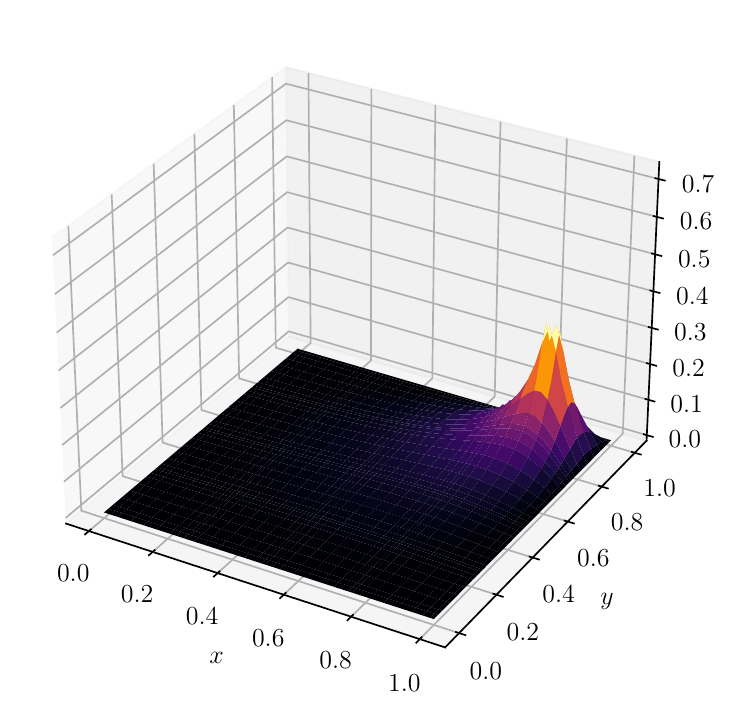}
    \caption{Exact $(s_1, s_2) = (0.875, 0.875)$}
\end{subfigure} \\
\begin{subfigure}{0.33\linewidth}
    \centering
    \includegraphics[width = \linewidth]{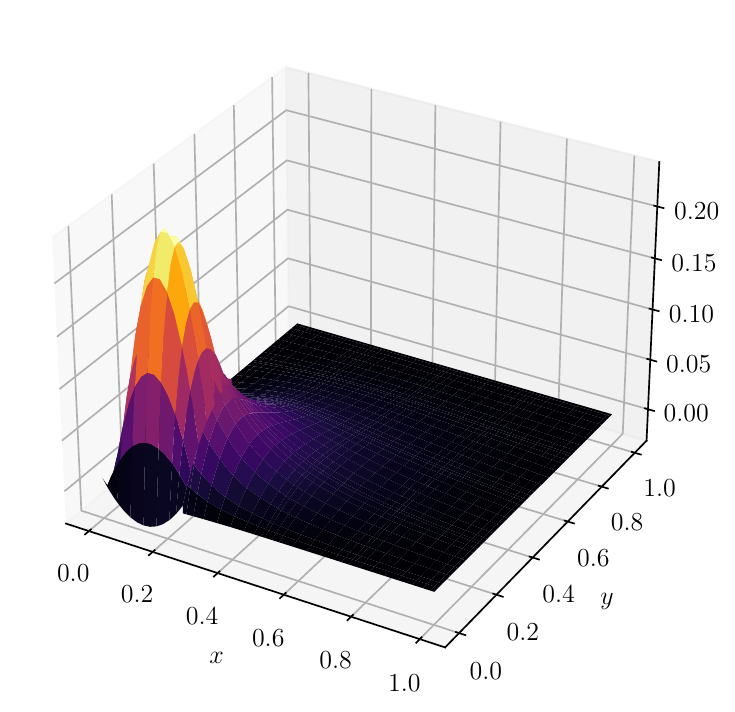}
    \caption{Approximate $(s_1, s_2) = (0.125, 0.125)$}
\end{subfigure}
\begin{subfigure}{0.33\linewidth}
    \centering
    \includegraphics[width = \linewidth]{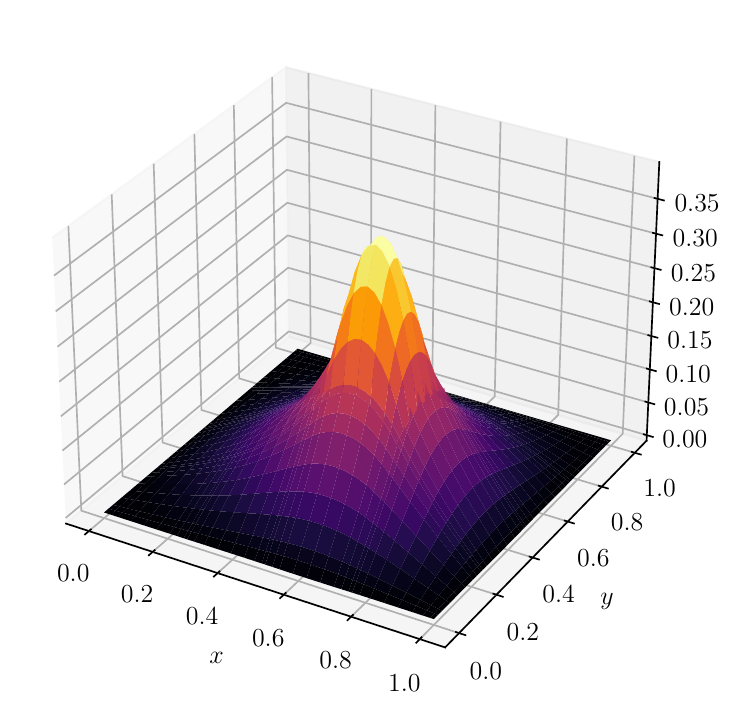}
    \caption{Approximate $(s_1, s_2) = (0.625, 0.375)$}
\end{subfigure}
\begin{subfigure}{0.33\linewidth}
    \centering
    \includegraphics[width = \linewidth]{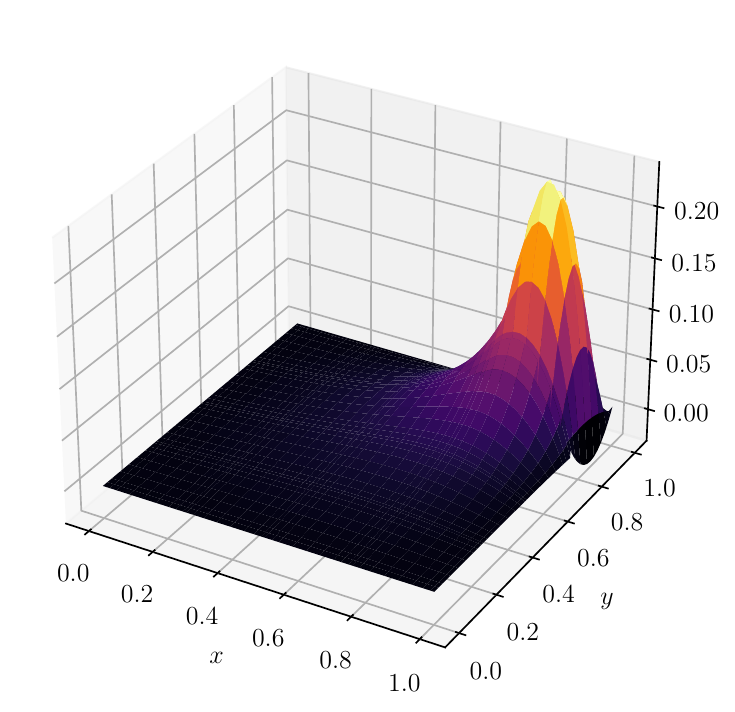}
    \caption{Approximate $(s_1, s_2) = (0.875, 0.875)$}
\end{subfigure}
\caption{Exact and approximate Greens' function for the 2D Poisson equation in the mixed formulation with the approximation computed using polynomials of degree 4 over $4 \times 4$ elements}
\label{fig:greens_approx_2D}
\end{figure}

\begin{figure}[H]
\begin{subfigure}{0.33\linewidth}
    \centering
    \includegraphics[width = \linewidth]{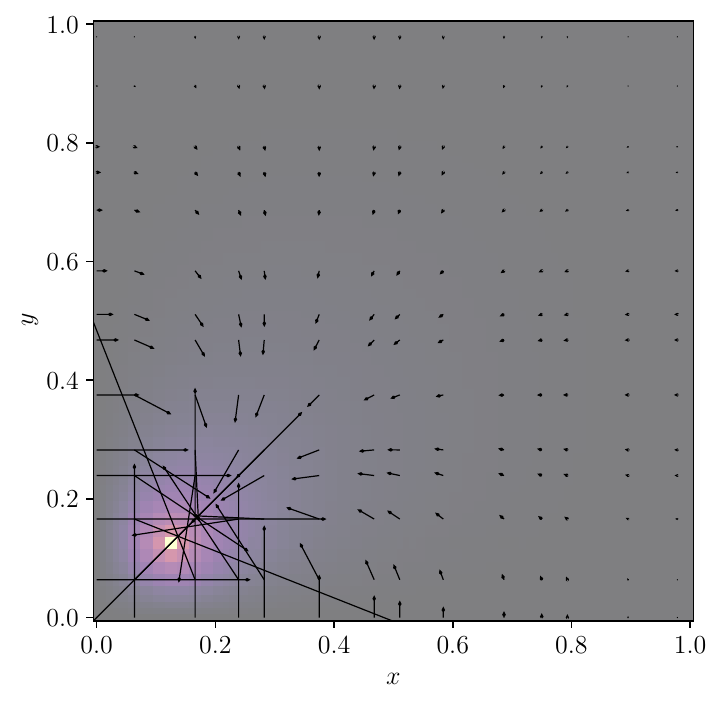}
    \caption{Exact $(s_1, s_2) = (0.125, 0.125)$}
\end{subfigure}
\begin{subfigure}{0.33\linewidth}
    \centering
    \includegraphics[width = \linewidth]{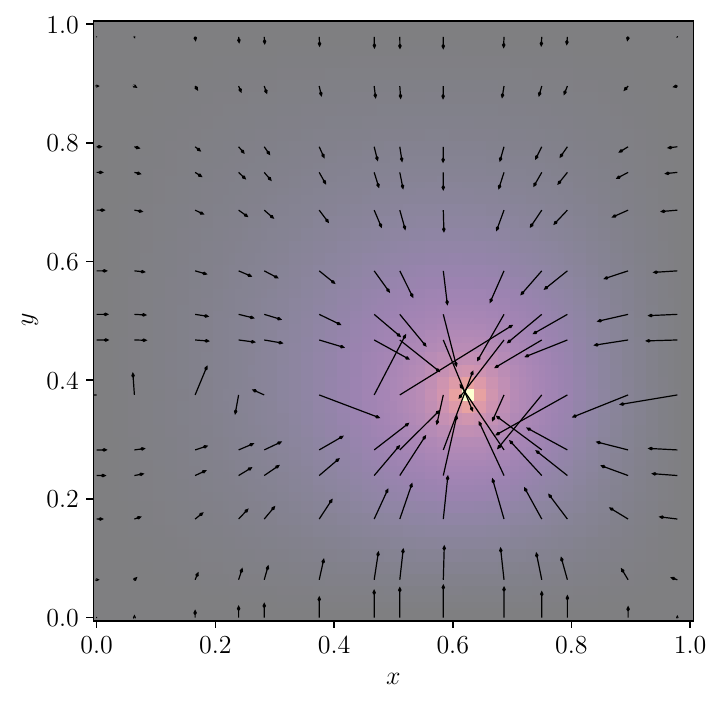}
    \caption{Exact $(s_1, s_2) = (0.625, 0.375)$}
\end{subfigure}
\begin{subfigure}{0.33\linewidth}
    \centering
    \includegraphics[width = \linewidth]{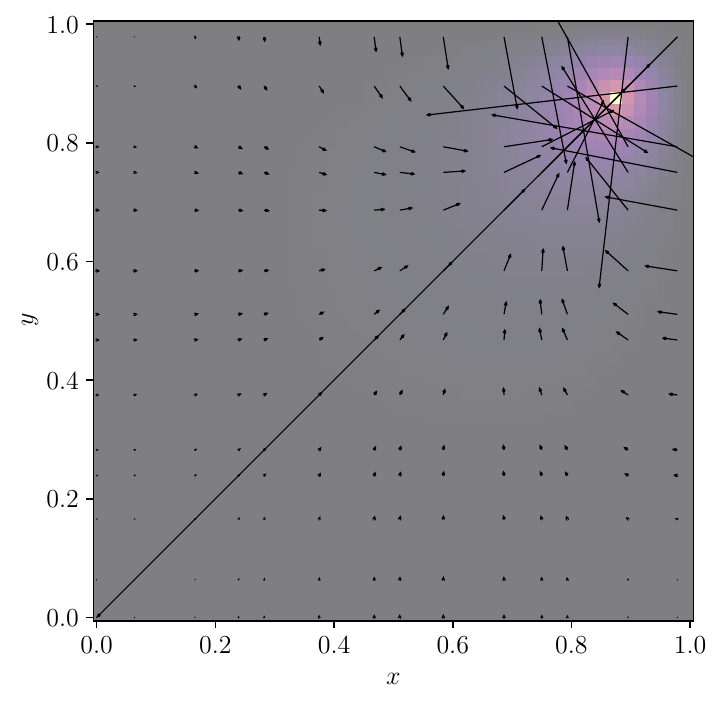}
    \caption{Exact $(s_1, s_2) = (0.875, 0.875)$}
\end{subfigure} \\
\begin{subfigure}{0.33\linewidth}
    \centering
    \includegraphics[width = \linewidth]{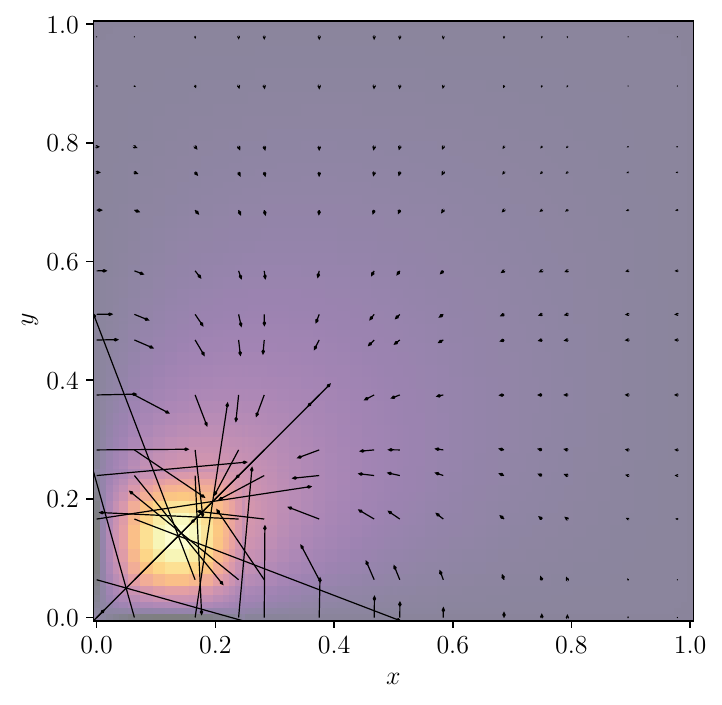}
    \caption{Approximate $(s_1, s_2) = (0.125, 0.125)$}
\end{subfigure}
\begin{subfigure}{0.33\linewidth}
    \centering
    \includegraphics[width = \linewidth]{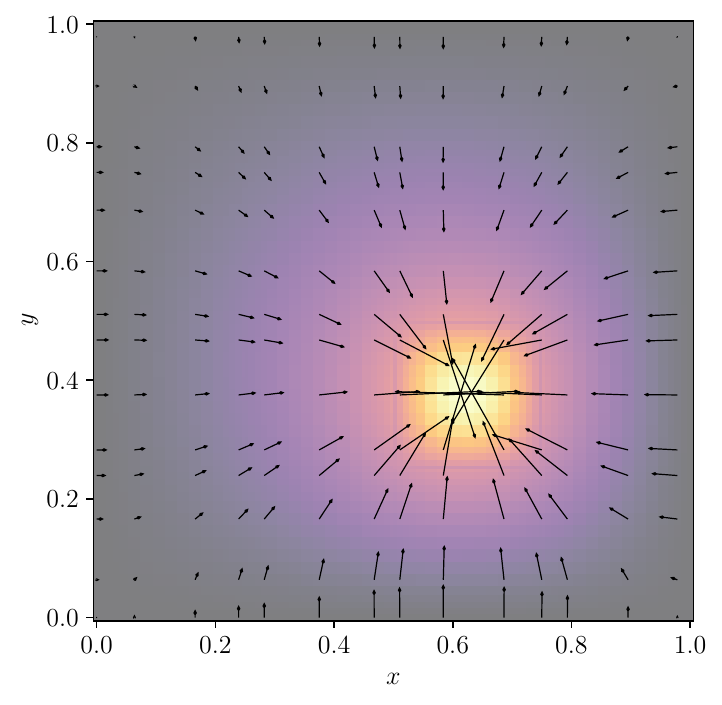}
    \caption{Approximate $(s_1, s_2) = (0.625, 0.375)$}
\end{subfigure}
\begin{subfigure}{0.33\linewidth}
    \centering
    \includegraphics[width = \linewidth]{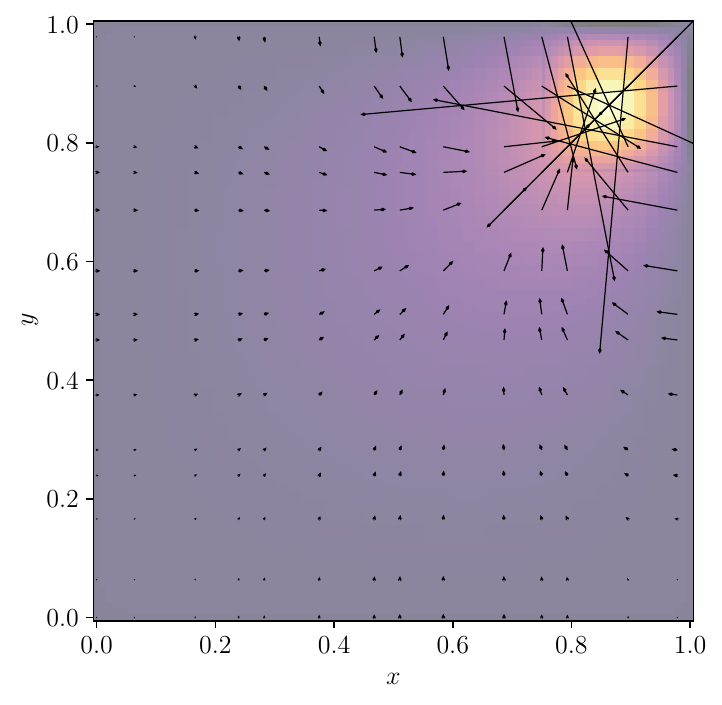}
    \caption{Approximate $(s_1, s_2) = (0.875, 0.875)$}
\end{subfigure}
\caption{Gradient of the exact and approximate Greens' function for the 2D Poisson equation in the mixed formulation with the approximation computed using polynomials of degree 4 over $4 \times 4$ elements}
\label{fig:dgreens_approx_2D}
\end{figure}

The extension to 2D is naturally done by taking the appropriate tensor products of the 1D basis functions to create the discrete $H(\mathrm{div}, \Omega)$ and $L^2(\Omega)$ subspaces \cite{Jain2021ConstructionMeshes}. \Cref{fig:greens_approx_2D} and \Cref{fig:dgreens_approx_2D} show the plots of the computed Greens' function and its gradient alongside the exact Greens' function from \eqref{eq:2dPoisson_Greens} with the infinite sum truncated to include the first 100 terms.

\end{document}